\documentclass[3p,12pt]{elsarticle}
\usepackage{amsthm}
\usepackage{amssymb}
\usepackage{amsmath}
\usepackage{float}
\usepackage{epstopdf}
\usepackage{graphicx}
\usepackage{multirow}
\usepackage{subcaption}
\usepackage{listings}
\usepackage{array}
\usepackage{flafter}
\usepackage[section]{placeins}
\usepackage{epsfig}
\usepackage{booktabs}
\usepackage{xcolor}
\usepackage[ruled]{algorithm2e}
\usepackage{color}
\usepackage{tikz}
\usetikzlibrary{arrows.meta, positioning, calc}
\usepackage{algpseudocode}
\numberwithin{equation}{section}

\newproof{pf}{Proof}[section]
\newdefinition{rmk}{Remark}
\definecolor{darkgreen}{rgb}{0,0.5,0}

\def\p{\partial}

\biboptions{sort&compress}

\journal{}
\begin{document}

\begin{frontmatter}

\title{A hybrid crossover kangaroo escape optimization framework for engineering optimization and UAV path planning}

\author[A]{Han-Bin Liu}
\author[A,B]{Huayi Wei\corref{cor1}}
\ead{weihuayi@xtu.edu.cn}
\cortext[cor1]{Corresponding author.}

\address[A]{School of Mathematics and Computational Science, Xiangtan University, Xiangtan 411105, China}
\address[B]{National Center of Applied Mathematics in Hunan, Hunan Key Laboratory for Computation
 and Simulation in Science and Engineering, Xiangtan 411105, China}

\begin{abstract}
Complex engineering optimization problems are often characterized by multimodality, high dimensionality, and nonlinear constraints, posing significant challenges for efficient and reliable computation.
To address these challenges, this paper develops a hybrid crossover-based optimization framework that enhances population interaction and improves search efficiency. The proposed framework integrates two complementary mechanisms, namely a Lévy long jump crossover strategy for global exploration and a horizontal–vertical crossover strategy for effective information exchange and local refinement, thereby improving convergence behavior and robustness.
To support scalable computation, the method is implemented within a unified multi-backend computational framework based on FEALPy, enabling consistent and efficient execution across heterogeneous platforms, including NumPy and PyTorch on both CPU and GPU. This design enhances portability, reproducibility, and computational efficiency in large-scale optimization tasks.
Extensive experiments on the IEEE CEC2022 benchmark suite demonstrate that the proposed framework achieves competitive or superior performance compared with several representative metaheuristic algorithms, as validated by Wilcoxon rank-sum and Friedman statistical tests. In addition, the method shows strong performance on constrained engineering design problems.
Finally, the proposed framework is applied to UAV path planning, formulated as a constrained optimization problem, demonstrating its effectiveness and scalability in complex engineering scenarios.
\end{abstract}


\end{frontmatter}

\section{Introduction}

Complex optimization problems in engineering and autonomous systems are often characterized by high dimensionality, multimodality, and nonlinear constraints, posing significant challenges for efficient computation.
These challenges become more critical in large-scale or real-time scenarios where both solution quality and computational efficiency are required.
In particular, problems such as UAV path planning can be formulated as constrained optimization tasks involving obstacle avoidance, trajectory smoothness, and dynamic feasibility requirements. Solving such problems requires not only robust optimization strategies but also scalable computational frameworks capable of operating efficiently across heterogeneous computing environments.

Path planning in unmanned aerial vehicle (UAV) systems is a representative application of constrained optimization that directly affects mission efficiency and operational safety.
It aims to determine an optimal trajectory under multiple constraints, including obstacle avoidance, path smoothness, and vehicle dynamic feasibility.
Various approaches have been proposed for UAV path planning, including graph-based methods such as Voronoi diagrams and probabilistic roadmap (PRM)~\cite{Beard2002,Mclain2005,Pettersson2006}, sampling-based techniques such as rapidly exploring random trees (RRT)~\cite{Lin2017}, grid-based search methods such as A*~\cite{Penin2019,Kwak2018}, and potential field approaches~\cite{Heidari2021,Di2015}.
Despite their effectiveness in specific scenarios, these methods often face challenges in scalability, handling complex constraints, or avoiding local optima in high-dimensional environments.
As a result, optimization-based approaches, particularly metaheuristic algorithms, have become an important direction for solving UAV path planning problems in complex and dynamic scenarios.

In recent years, metaheuristic algorithms have been widely adopted for solving complex engineering optimization problems due to their flexibility and ability to handle high-dimensional and non-differentiable search spaces.
Representative methods include Particle Swarm Optimization (PSO)~\cite{Manh2021}, Genetic Algorithm (GA)~\cite{Roberge2013}, and Differential Evolution (DE)~\cite{Sun2016}.
Among recent developments, the Kangaroo Escape Optimizer (KEO), proposed by Almutairi et al.~\cite{Almutairi2025}, is a metaheuristic algorithm that models escape behaviors to balance exploration and exploitation.
Although KEO has demonstrated promising performance across various optimization problems~\cite{Sulaiman2025,Abdullah}, it may still suffer from premature convergence and insufficient exploration capability when dealing with high-dimensional or complex optimization landscapes.
From a methodological perspective, existing KEO variants primarily rely on movement-based search mechanisms, which limits the diversity of information exchange among individuals.
In particular, the potential of crossover-based interactions for enhancing population diversity and accelerating convergence has not been fully explored.
In addition, most existing implementations are tightly coupled to specific programming environments, restricting scalability and portability in heterogeneous computing settings.

To address these limitations, this study develops a hybrid crossover-based optimization framework, referred to as HCKEO.
The proposed framework introduces two complementary crossover strategies to enhance information exchange among individuals and improve the balance between global exploration and local exploitation, thereby mitigating premature convergence and insufficient exploration in complex search spaces.
Furthermore, a unified multi-backend computational framework is constructed based on FEALPy, enabling efficient and consistent execution across heterogeneous platforms.
The effectiveness of the proposed framework is validated through experiments on the IEEE CEC2022 benchmark suite and constrained engineering optimization problems.
In addition, it is applied to UAV path planning formulated as a constrained optimization problem, demonstrating its applicability in complex real-world scenarios.
The main contributions of this study are summarized as follows:
\begin{itemize}
	\item A hybrid crossover-based optimization framework is developed to address complex engineering optimization problems characterized by multimodality, high dimensionality, and nonlinear constraints. The framework enhances population interaction and improves overall search efficiency.

	\item Two complementary crossover strategies, namely the L\'{e}vy long jump crossover and the horizontal–vertical crossover, are integrated within the framework to achieve a better balance between global exploration and local exploitation, leading to improved convergence behavior and robustness.

	\item A unified multi-backend computational framework is established based on FEALPy, enabling consistent and efficient execution across heterogeneous platforms (NumPy and PyTorch on CPU and GPU), thereby improving scalability, portability, and reproducibility.

	\item Comprehensive evaluations on the IEEE CEC2022 benchmark suite and constrained engineering optimization problems demonstrate the effectiveness and robustness of the proposed framework compared with representative metaheuristic methods.

	\item The proposed framework is applied to UAV path planning formulated as a constrained optimization problem, demonstrating its capability in handling complex real-world engineering scenarios.
\end{itemize}

The structure of the remainder of this paper is summarized as follows.
Section~\ref{RW} provides a brief review of related studies.
Section~\ref{HCKEO} describes the basic KEO framework and details the proposed HCKEO algorithm.
Section~\ref{EA} reports the experimental evaluation on benchmark functions and engineering optimization problems.
Section~\ref{UAV} discusses the UAV path planning application.
Section~\ref{fealpy} presents the multi-backend implementation of HCKEO based on the FEALPy framework and analyzes its computational performance across different computing platforms.
Finally, Section~\ref{conclusion} concludes the paper.

\section{Related Work}\label{RW}

\subsection{Metaheuristic algorithms in optimization}

In recent years, many real-world engineering and scientific problems have been formulated as optimization problems whose complexity continues to increase.
To address such challenges, MHAs have gained widespread attention in recent years because of their flexibility and strong global optimization capability.
MHAs are population-based optimization methods inspired by natural, biological, or social phenomena~\cite{Ali2021}, and have become effective tools for solving complex optimization problems.

\begin{figure}[!ht]
  \centering
  \includegraphics[width=0.55\linewidth]{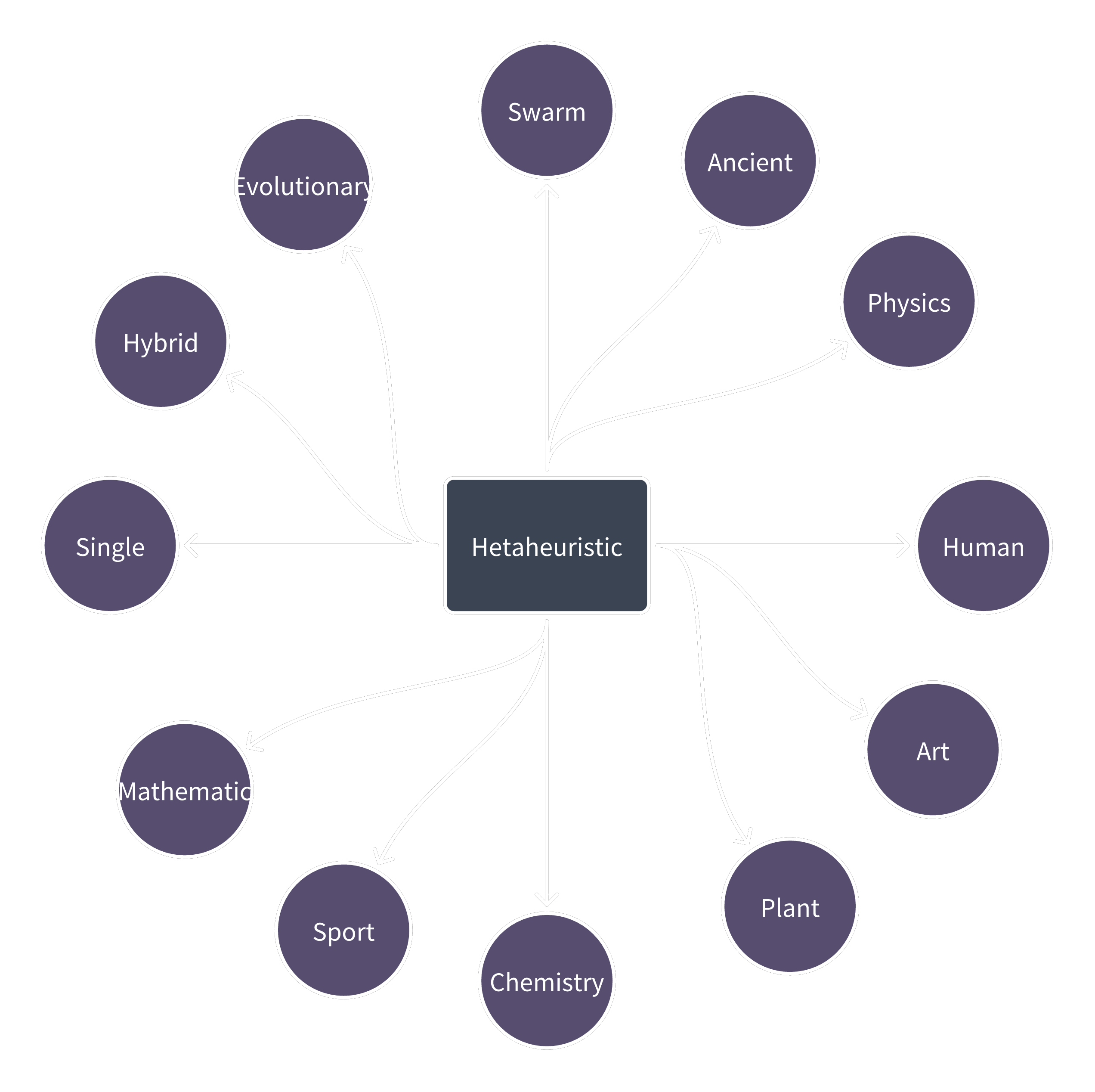}
  \caption{Classification of metaheuristic algorithms.}
  \label{fig:Metaheuristics classification}
\end{figure}

Over the past decades, a large number of MHAs have been proposed and can be broadly classified into several categories~\cite{Peraza2024}, including evolutionary algorithms (Genetic Algorithm \cite{Roberge2013}), swarm-based algorithms (Crayfish optimization algorithm \cite{COA}), physics-based algorithms (Snow ablation optimizer \cite{Deng2023}), chemistry-based algorithms (Chemical Reaction Optimization \cite{Lam2012}), plant-based algorithms (Strawberry Algorithm \cite{Khan2018}), human-based algorithms (Society Civilization Algorithm \cite{Ray2003}), art-based algorithms (Musical Composition Algorithm \cite{Mora2014}), sport-based algorithms (Running City game optimizer \cite{Ma2023}), mathematics-based algorithms (Stochastic Fractal Search \cite{Salimi2015}), single-solution-based algorithms (Large Neighbourhood Search \cite{Pisinger2019}), and hybrid algorithms (Cuckoo Search and Quantum-Behavior Particle Swarm Optimization \cite{Xiongfa2023}, Butterfly Optimization Algorithm Combined with Black Widow Optimization \cite{Xu2022}), as illustrated in Fig.~\ref{fig:Metaheuristics classification}.
%

The popularity of MHAs can be attributed to several advantages, including simple implementation, derivative-free search capability, and strong adaptability to complex and nonlinear optimization problems~\cite{Kumar2021}.
Unlike gradient-based optimization methods, MHAs rely on stochastic population evolution to explore the search space, which helps avoid premature convergence and improves robustness.
However, according to the No Free Lunch (NFL) theorem \cite{NFL}, no single algorithm performs best for all optimization problems.
Therefore, continuous development of new metaheuristics and hybrid strategies remains an important research direction in optimization.
Among recently proposed MHAs, the Kangaroo Escape Optimizer has attracted increasing attention due to its adaptive escape-inspired search mechanisms and promising optimization performance.


\subsection{Kangaroo Escape Optimizer and its variants}

Recently, the Kangaroo Escape Optimizer (KEO) has been proposed as a nature-inspired metaheuristic algorithm that simulates the escape strategies of kangaroos when facing potential threats.
In KEO, each kangaroo represents a candidate solution that dynamically moves within the search space to locate a safer region corresponding to a better objective value.
The algorithm models three main behavioral mechanisms, including zigzag motion, long-jump escape, and decoy drop strategies, which collectively balance global exploration and local exploitation during the optimization process~\cite{Abdullah}.

Due to its simple structure and adaptive search mechanism, KEO has recently attracted increasing attention in engineering optimization problems.
For example, Aljumah and Abdullah~\cite{Abdullah} applied KEO to parameter extraction of photovoltaic models based on different equivalent circuit configurations, demonstrating improved accuracy and convergence efficiency compared with several state-of-the-art optimizers.
Similarly, Sulaiman \emph{et al.}~\cite{Sulaiman2025} introduced a Kangaroo Escape Algorithm (KEA) for solving combined heat and power economic dispatch problems, showing promising performance in large-scale energy management scenarios.
In addition, Almutairi and Shaheen~\cite{Almutairi2025} proposed a Kangaroo Escape Optimization to optimize fractional-order PID controller in load frequency control of interconnected power systems.

Despite the encouraging results reported in engineering applications, existing research on KEO mainly focuses on practical implementations and parameter tuning rather than algorithmic enhancement.
The potential of integrating crossover-based information exchange mechanisms to enhance population diversity and search capability has received limited attention.
However, the original KEO primarily relies on movement-based search strategies such as zigzag motion and long-jump escape, which may limit information exchange among kangaroos and reduce population diversity in complex optimization problems.
Motivated by this observation, the present work introduces two crossover-based strategies to improve the exploration–exploitation balance of the original KEO framework.

\subsection{Hybrid and crossover-based strategies}

Hybridization strategies have been widely adopted to enhance the performance of MHAs by integrating complementary search mechanisms.
By combining different operators or search behaviors, hybrid metaheuristics can effectively improve exploration capability, convergence speed, and solution robustness in complex optimization problems \cite{Jialing2025,Heming2024}.

Among various hybridization techniques, crossover-based operators play a crucial role in promoting information exchange among agents and maintaining population diversity.
Originally introduced in evolutionary algorithms such as GA, crossover operations enable candidate solutions to share useful information, thereby improving the global search ability of population-based optimization methods.
Recent studies have demonstrated that incorporating crossover mechanisms into swarm intelligence and other metaheuristic frameworks can significantly enhance convergence stability and reduce the risk of premature convergence \cite{Yanyun2023}.

In addition, L\'{e}vy flight has been widely employed in MHAs to improve global exploration performance \cite{Xize2026,Zhenyu2022}.
Due to its heavy-tailed step-length distribution governed by a power law, L\'{e}vy flight enables a combination of frequent short moves and occasional long-distance jumps, which helps algorithms escape from local optima and explore distant promising regions.
Therefore, incorporating L\'{e}vy-based search steps into existing metaheuristic frameworks has become an effective approach for enhancing global exploration capability.

Therefore, integrating crossover mechanisms and L\'{e}vy-based exploration into the KEO framework is expected to enhance both population diversity and global search capability.
Motivated by these observations, this study introduces two crossover-based enhancement strategies into the original KEO framework, namely the L\'{e}vy Long Jump Crossover (LLJC) strategy and the Horizontal–Vertical Crossover (HVC) strategy.
These mechanisms aim to improve information exchange among individuals and achieve a better balance between global exploration and local exploitation.

\section{The Proposed Algorithm}\label{HCKEO}

In this section, we first present a brief overview of the original KEO to clarify its core mechanisms and optimization principles.
Then, we introduce the proposed HCKEO, which improves the exploration–exploitation balance of the standard KEO by integrating two enhanced strategies: the L\'{e}vy Long Jump Crossover (LLJC) strategy and the Horizontal–Vertical Crossover (HVC) strategy.
The overall framework, algorithmic improvements, and computational procedure of HCKEO are described in detail in the following subsections.

\subsection{Kangaroo Escape Optimizer}

The KEO \cite{Almutairi2025} is a nature-inspired MHA that emulates the escape behavior of kangaroos when facing predators.
In nature, kangaroos exhibit adaptive and cooperative movement strategies that enable them to efficiently evade threats.
Motivated by this behavior, KEO models the optimization process as a group of kangaroos that dynamically adjust their positions in the search space through three coordinated mechanisms: Energy-Based Decision to adaptive search strategy, zigzag motion to avoid direct pursuit, long jumps to traverse distant regions, and decoy drops to mislead predators.
These mechanisms collectively achieve a balance between global exploration and local exploitation during the search process.

\subsubsection{Energy-Based Decision}

In the original KEO algorithm, the chaotic logistic map is used to model the energy of each kangaroo  to represent the dynamic energy variation of kangaroos during the search process.
The energy value at \(t\) iteration is defined as:
\begin{equation}\label{energy}
energy(t) = \left(1 - r \cdot \frac{t}{T}\right) \cdot \left(0.95 + 5 \cdot \frac{\lambda(t)}{100}\right),
\end{equation}
where \( t \) and \( T \) denote the current and maximum iteration numbers, respectively; \( r \in [0,1] \) is a random number; and \( \lambda(t) \) is updated according to the logistic map:
\begin{equation}\label{lambda}
\lambda(t+1) = 4 \cdot \lambda(t) \cdot (1 - \lambda(t)),
\end{equation}
with the initial value \( \lambda(0) = 0.7 \).
This mechanism enables kangaroos to exhibit varying activity levels throughout the optimization process, maintaining a controlled degree of randomness while gradually reducing their energy as the iterations proceed.

\subsubsection{Zigzag Motion}

Let $X_i^t$ denote the position of the $i$th kangaroo at iteration $t$.
The zigzag motion is mathematically modeled as:
\begin{equation}
X_i^{t+1} = X_i^t + \sin(\theta)\, \beta\, v_{rotate}\, N_1,
\label{zigzag motion}
\end{equation}
where the second term on the right-hand side represents the zigzag step, $\beta = 0.5$ is a scaling coefficient, $N_1$ is a random vector drawn from a standard normal distribution, and $v_{rotate}$ is the rotated direction vector defined as:
\begin{equation}
v_{rotate} = \sin(\theta)\, U\, \|V_i\| + \cos(\theta)\, V_i,
\label{rotated direction vector}
\end{equation}
where $U$ is an orthogonal unit vector perpendicular to $V_i$, the direction vector from $X_i$ to the global best position ($gbest$).
The rotation angle $\theta$ is given by:
\begin{equation}\label{random angle}
\theta = \theta^{\max}(2r_1 - 1),
\end{equation}
where $\theta^{\max} = \pi/6$ denotes the maximum allowable rotation angle, and the random variable $r_1$ sampled from the uniform distribution over $(0,1)$.
This mechanism introduces directional perturbations to maintain diversity while guiding the population toward promising regions.

\subsubsection{Decoy Drops}\label{DD}

When no immediate threat is perceived, kangaroos tend to move toward safer and resource-rich regions.
This behavior is modeled as:
\begin{equation}\label{decoy drops}
X_i^{t+1} = X_{safe} + decoy_{drop}\, N_3\, (X_i^t - X_{safe}),
\end{equation}
where $X_{safe}$ denotes a position considered safe, determined under three movement modes:
\textit{(i) exploration} – selecting a random kangaroo position;
\textit{(ii) exploitation} – adopting a locally best position within a small group; and
\textit{(iii) deep exploitation} – moving toward the global best position.
The term $decoy_{drop}$ characterizes the decoy-dropping behavior and is defined as:
\begin{equation}\label{dd}
decoy_{drop} =
\begin{cases}
    \text{round}(r_2 \times r_3), & \text{if } r_5 \ge 2/3, \\
    \text{round}(r_4), & \text{if } 2/3 > r_5 \ge 1/3, \\
    [1,1,\dots,1]_{1 \times D}, & \text{otherwise,}
\end{cases}
\end{equation}
where $r_2$, $r_3$, and $r_4$ are random vectors within $[0,1]$, and $r_5 \in [0,1]$ is a scalar random variable.
This mechanism enhances stochasticity and enables individuals to exhibit diverse adaptive behaviors, thereby mitigating the risk of premature convergence.

\subsubsection{Long Jump}

Kangaroos may search for distant areas through long-distance jump, which can be formulated as:
\begin{equation}\label{long jump}
X_i^{t+1} = X_i^t + 2\, X_i^t \cdot N_2 \cdot decoy_{drop},
\end{equation}
where $N_2$ denotes a random vector sampled from a standard normal distribution.
The inclusion of $decoy_{drop}$ introduces controlled randomness, enhancing the algorithm’s global exploration ability and enabling individuals to effectively escape from local optima.

\subsection{Hybrid Crossover Kangaroo Escape Optimizer}

Although KEO adopts three strategies to balance exploration and exploitation, its convergence speed and local refinement capability can still be further improved, particularly in high-dimensional or complex search spaces.
To address these limitations, we propose the Hybrid Crossover Kangaroo Escape Optimizer (HCKEO), which integrates two key mechanisms: the L\'{e}vy Long Jump Crossover (LLJC) strategy and the Horizontal–Vertical Crossover (HVC) strategy.
By combining these crossover-based improvements with the original KEO framework, HCKEO achieves a more robust and adaptive optimization process.
The details of these strategies and their computational procedures are elaborated in the following subsections.

\subsubsection{L\'{e}vy Long Jump Crossover Strategy}

To address the exploration–exploitation imbalance and reduce premature convergence in KEO, a L\'{e}vy Long Jump Crossover (LLJC) strategy is introduced.
This strategy leverages the statistical characteristics of L\'{e}vy flights, which generate a combination of frequent short moves and occasional long-distance jumps.
Such long-tailed step-size distributions, governed by a power-law behavior, have been demonstrated to significantly enhance global search ability and improve convergence robustness in swarm intelligence and evolutionary algorithms \cite{Hanbin2024}.
Furthermore, recent studies suggest that combining L\'{e}vy-based movements with other search strategies can effectively balance diversification and intensification, thereby improving overall optimization efficiency \cite{Xiongfa2023}.
By embedding the L\'{e}vy-based crossover operator within KEO, the proposed approach enables non-Gaussian stochastic movements that facilitate escaping local optima and exploring distant promising regions, ultimately enhancing adaptability and convergence performance.

The mathematical model of the LLJC strategy is formulated as:
\begin{equation}\label{lcs}
X_i^{t+1} = X_i^t + 2 \cdot N_2 \cdot decoy_{drop} \cdot (X_{r_1}^t - gbest) + k \cdot l(\lambda) \cdot (X_{r_2}^t - X_i^t),
\end{equation}
where $X_{r_1}^t$ and $X_{r_2}^t$ are two randomly selected kangaroos at iteration $t$, and $k \cdot l(\lambda)$ represents the L\'{e}vy crossover factor defined as:
\begin{eqnarray}
k &=& 2\, r_6 \cdot \frac{T - t}{T}, \\
l(\lambda) &=& \frac{r_a}{\|r_b\|^{1/\alpha}} \cdot \left(\frac{\Gamma(1+\alpha)\sin(\pi\alpha/2)}{\Gamma((1+\alpha)/2)\,\alpha\,2^{(\alpha-1)/2}}\right)^{1/\alpha},
\end{eqnarray}
where $T$ is the maximum number of iterations, $\Gamma(\cdot)$ denotes the Gamma function, $r_6 \in [0,1]$ is a random variable, $r_a$ and $r_b$ represent Gaussian random variables with zero mean and unit variance, and $\alpha$ is a constant coefficient whose value is fixed at 1.5 in HCKEO.
This mechanism integrates the exploratory characteristics of L\'{e}vy flight–like movements with the information-sharing capability introduced by crossover operations.
Consequently, individuals are able to perform adaptive large-step movements guided by high-quality solutions, while preserving stochastic perturbations that help the population avoid premature convergence.

\begin{figure}[!ht]
  \centering
  \includegraphics[width=0.55\linewidth]{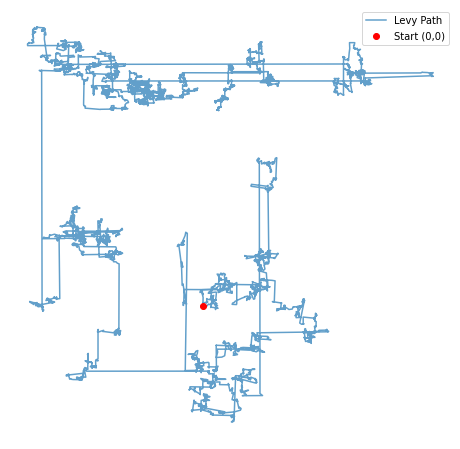}
  \caption{Illustration of L\'{e}vy flight step-length distribution at different iteration stages.
  Larger steps occur with higher probability during early iterations, facilitating global exploration,
  while the step size gradually decreases as the iteration proceeds, encouraging local exploitation and convergence.}
  \label{fig:levy_step}
\end{figure}

It is worth noting that the effective step length induced by the L\'{e}vy flight is iteration-dependent.
At early stages, larger step sizes dominate the search process to enhance global exploration, whereas the step length gradually shrinks as the iteration progresses, promoting local exploitation.
This dynamic behavior is illustrated in Fig.~\ref{fig:levy_step}.

\subsubsection{Horizontal–Vertical Crossover Strategy}

To further enhance information exchange within the population, the Horizontal–Vertical Crossover (HVC) strategy \cite{Xiongbin2025} is incorporated into the proposed HCKEO framework.
The HVC operator introduces two complementary crossover mechanisms that operate from different perspectives of the search space, namely horizontal crossover and vertical crossover.

The horizontal crossover mechanism focuses on interactions between different individuals.
In this operation, two individuals are randomly selected from the population and their position vectors are combined to generate new candidate solutions.
By allowing individuals to exchange useful search information, horizontal crossover helps improve the global exploration capability of the algorithm and accelerates convergence toward promising regions.

Let $X_{i1}$ and $X_{i2}$ denote two kangaroos selected from the population.
The offspring generated through horizontal crossover can be expressed as:
\begin{equation}
X^{hc}{i1} = r_7 X{i1} + (1-r_7)X_{i2} + c_1 (X_{i1}-X_{i2}),
\end{equation}
\begin{equation}
X^{hc}{i2} = r_8 X{i2} + (1-r_8)X_{i1} + c_2 (X_{i2}-X_{i1}),
\end{equation}
where $r_7$ and $r_8$ are uniformly distributed random numbers in $(0,1)$, and $c_1$ and $c_2$ are random coefficients sampled from $(-1,1)$.

Different from horizontal crossover, the vertical crossover mechanism operates across different dimensions of the same individual.
During the search process, certain dimensions may stagnate and limit the algorithm’s ability to further improve the solution.
Vertical crossover alleviates this issue by introducing interactions among dimensions within an individual.

Specifically, two dimensions of an individual are randomly selected and combined to update one target dimension, while the remaining dimensions remain unchanged.
This dimension-level recombination introduces additional variation into the population and helps the algorithm escape from locally trapped regions of the search space.
The vertical crossover operation is defined as:
\begin{equation}\label{vc}
X^{vc}_{i,j_1} = r_9 \cdot X_{i,j_1} + (1 - r_9) \cdot X_{i,j_2}
\end{equation}
where $j_1$ and $j_2$ represent two randomly selected dimensions and $r_9 \in (0,1)$ is a uniformly distributed random coefficient.
By jointly employing horizontal and vertical crossover mechanisms, the HVC strategy promotes both inter-individual information exchange and intra-individual dimensional interaction.
This dual interaction mechanism enhances population diversity, improves the ability to escape premature convergence, and ultimately strengthens the overall search performance of HCKEO.

\subsubsection{HCKEO Algorithmic Process}

\begin{figure}[!ht]
\centering
\includegraphics[width=\linewidth]{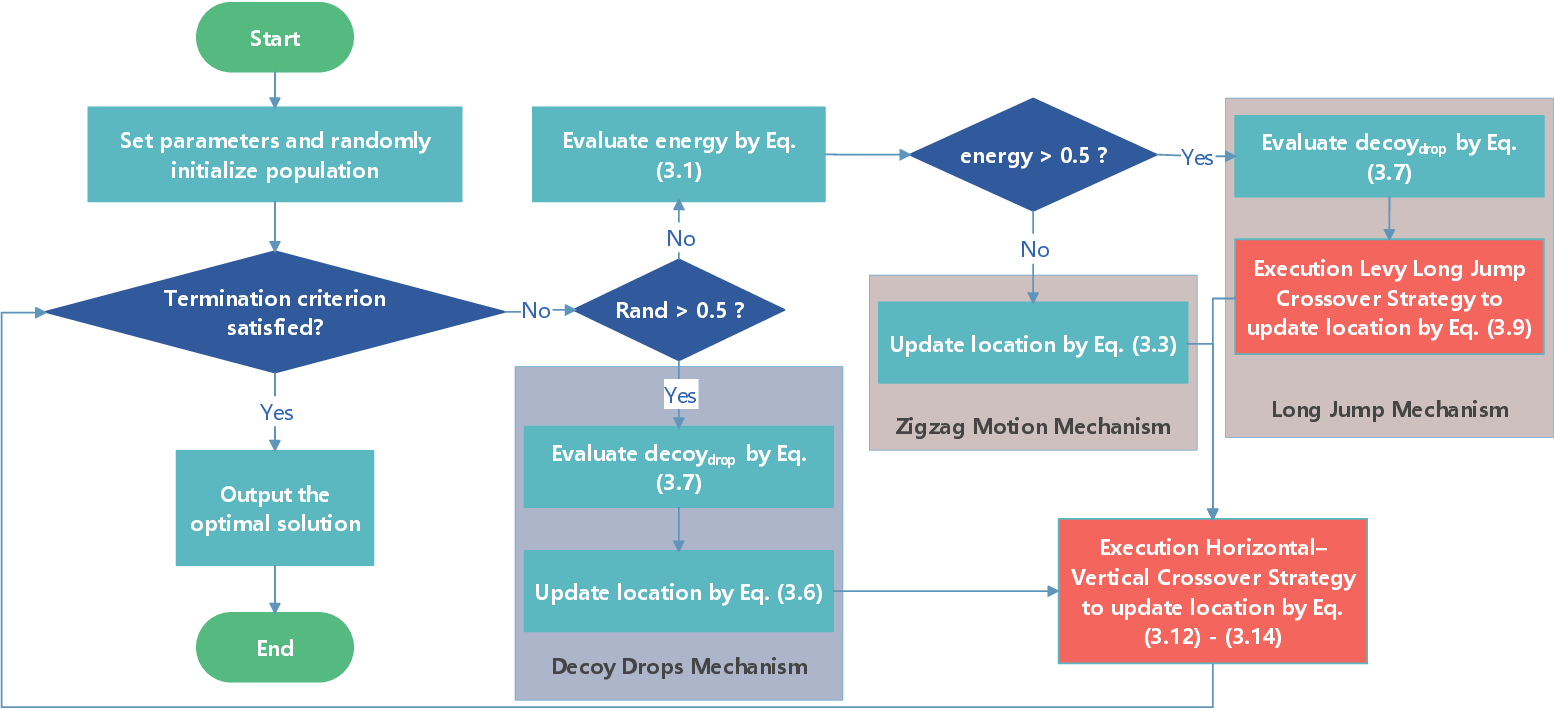}
\caption{Flowchart of the proposed HCKEO algorithm.}
\label{flowchat}
\end{figure}

The flowchart of the proposed HCKEO algorithm is shown in Fig.~\ref{flowchat}.
The corresponding pseudo-code is presented in Algorithm~\ref{pseudo}.

\begin{algorithm}[!ht]
	\caption{Hybrid Crossover Kangaroo Escape Optimizer (HCKEO)}\label{pseudo}
	\LinesNumbered
	\KwIn{Population size, problem dimension, and maximum number of iterations}
	\KwOut{Optimal solution}
	
	Initialize the population\;
	
	\While{the termination criterion is not satisfied}{
        \For{each kangaroo}{
            \uIf{rand $>$ 0.5}{
                \% \textbf{Decoy Drops Mechanism}\;
                Compute $decoy_{drop}$ using Eq.~(\ref{dd})\;
                Update position using Eq.~(\ref{decoy drops})\;
    		}
            \Else{
                Evaluate energy using Eq.~(\ref{energy})\;
                \uIf{energy $>$ 0.5}{
                    \% \textbf{Long Jump Mechanism}\;
                    Compute $decoy_{drop}$ using Eq.~(\ref{dd})\;
                    Execute the L\'{e}vy Long Jump Crossover Strategy using Eq.~(\ref{lcs})\;
        		}
                \Else{
                    \% \textbf{Zigzag Motion Mechanism}\;
                    Update position using Eq.~(\ref{zigzag motion})\;
                }
            }
        }
        \% \textbf{Population Update}\;
        Apply the Horizontal-Vertical Crossover Strategy to update kangaroo states\;
        Update the population by retaining the best current positions\;
	}
\end{algorithm}

\subsubsection{HCKEO Performance Analysis}

This section analyzes the computational complexity and memory requirements of HCKEO, corresponding to the algorithm’s time and space complexity, respectively.

\paragraph{Computational Complexity}

The computational cost of the proposed HCKEO algorithm mainly arises from population initialization, fitness evaluation, position updating, and the HVC operations.
Before discussing the computational complexity of the algorithm, it is agreed that $N $represents the population size, $D $represents the dimension, and $T $represents the maximum number of iterations.
During initialization, the population is randomly generated within the search space, producing an ($N \times D$) matrix of candidate solutions.
This step requires ($O(ND)$) time.
For each iteration of the optimization process, the algorithm evaluates the fitness of all individuals and updates their positions according to the movement mechanisms of KEO together with the proposed crossover operators.
Both the fitness computation and the position update involve operations across all individuals and dimensions, leading to a complexity of ($O(ND)$) per iteration.
In addition, the HVC strategy performs arithmetic operations on solution vectors with the same order of computational cost, i.e., ($O(ND)$).

Consequently, the dominant cost of the algorithm occurs during the iterative search stage, and the omputational complexity of HCKEO can be calculated as:
\[
O(ND) + O(ND) + O(TND) + O(TND) = O(TND).
\]
Since the original KEO algorithm follows a similar population update mechanism, it shares the same asymptotic time complexity ($O(TND)$).
Therefore, incorporating the crossover operators does not change the order of computational complexity, although it introduces additional constant-level arithmetic operations within each iteration.

\paragraph{Memory Efficiency}

Space complexity measures the amount of memory resources required by an algorithm.
It measures the amount of additional memory required during the execution of the algorithm and how this requirement grows with the problem size.
In the analysis of space complexity, attention is mainly given to the memory usage associated with data structures, variables, pointers, and potential recursive operations.
This analysis provides insight into memory utilization, facilitates the selection of appropriate data structures, and supports improvements in memory efficiency.

The parameter settings of HCKEO occupy a constant amount of memory, which is independent of the problem size.
The overall space complexity of HCKEO is $O(N)$, where $N$ denotes the population size.
This memory requirement mainly originates from maintaining and updating the population throughout the iterative optimization process.
Since the space complexity of KEO is also $O(N)$, the two algorithms share the same linear space complexity property.

In summary, HCKEO maintains the same order of time and space complexity as KEO, i.e., $O(TND)$ and $O(N)$, respectively.
This indicates that the additional mechanisms in HCKEO enhance optimization performance without increasing the algorithm’s computational or memory burden.

\section{Experimental Evaluation and Statistical Analysis}\label{EA}

In this section, the performance of the proposed HCKEO algorithm is comprehensively evaluated through a series of benchmark experiments.
The evaluation includes four aspects:
(1) effectiveness verification of the proposed components,
(2) exploration–exploitation analysis,
(3) testing on the CEC2022 benchmark suite, and
(4) validation on real-world engineering optimization problems.

\paragraph{Experimental Environment and Settings}

The experimental environment is summarized in Table~\ref{tab:environment}.
Unless otherwise specified, all algorithms were run using the same maximum number of iterations and population size to guarantee consistency across comparisons.
The performance was evaluated in terms of convergence accuracy, robustness, and computational efficiency.

\begin{table}[!ht]
\centering
\caption{Experimental environment.}
\begin{tabular}{ll}
\hline
\textbf{Component} & \textbf{Specification} \\
\hline
CPU & 11th Gen Intel(R) Core(TM) i7-11800H (2.30\,GHz) \\
RAM & 16\,GB (2\,×\,8\,GB) \\
Operating System & Ubuntu 24.04 \\
Implementation & Python 3.10.14 \\
\hline
\end{tabular}
\label{tab:environment}
\end{table}

\subsection{Effectiveness Verification of Different Strategies}

To enhance the optimization performance of KEO, two improvement strategies are introduced: the LLJC Strategy and the HVC Strategy.
To verify the effectiveness of each component, two algorithmic variants, namely LLJC-KEO and HVC-KEO, are designed, as summarized in Table~\ref{tab:cross_variants}.
In this table, “+” indicates that the corresponding strategy is incorporated, while “–” means it is not applied.

\begin{table}[!ht]
  \centering
  \caption{Different KEO variants with incorporated improvement strategies.}
  \label{tab:cross_variants}
  \begin{tabular}{lccc}
    \hline\noalign{\smallskip}
    Parameters & LLJC-KEO & HVC-KEO & HCKEO \\
    \noalign{\smallskip}\hline\noalign{\smallskip}
    L\'{e}vy Long Jump Crossover Strategy & + & – & + \\
    Horizontal–Vertical Crossover Strategy & – & + & + \\
    \noalign{\smallskip}\hline
  \end{tabular}
\end{table}

All algorithms were tested on the 10-dimensional CEC2022 benchmark suite \cite{CEC2022} with a population size of 50 and a maximum of 500 iterations, which is given in Table \ref{tab:cec2022}.
\begin{table}[!ht]
\centering
\caption{CEC2022 Benchmark Functions}
\resizebox{\textwidth}{!}{
\begin{tabular}{lllll}
\toprule
\textbf{Category} & \textbf{Function} & \textbf{Description} & \textbf{Dimension} & \textbf{Optimum} \\
\midrule
Unimodal & F1 & Shifted and rotated Zakharov function & 10/20 & 300 \\
Multimodal & F2 & Shifted and rotated Rosenbrock function & 10/20 & 400 \\
 & F3 & Shifted and rotated Rastrigin function & 10/20 & 600 \\
 & F4 & Shifted and rotated non-continuous Rastrigin function & 10/20 & 800 \\
 & F5 & Shifted and rotated Levy function & 10/20 & 900 \\
Hybrid & F6 & Hybrid function (N = 3) & 10/20 & 1800 \\
 & F7 & Hybrid function (N = 6) & 10/20 & 2000 \\
 & F8 & Hybrid function (N = 5) & 10/20 & 2200 \\
Composition & F9 & Composition function (N = 5) & 10/20 & 2300 \\
 & F10 & Composition function (N = 4) & 10/20 & 2400 \\
 & F11 & Composition function (N = 5) & 10/20 & 2600 \\
 & F12 & Composition function (N = 6) & 10/20 & 2700 \\
\bottomrule
\end{tabular}}
\label{tab:cec2022}
\end{table}

The same parameter settings are maintained for the subsequent experiments for consistency.
Friedman’s non-parametric test was employed as the performance evaluation criterion.
The comparative results are presented in Table~\ref{tab:cec2022_variants}.

\begin{table}[!ht]
  \centering
  \caption{Performance comparison of KEO variants on the 10-dimensional CEC2022 benchmark functions.}
  \label{tab:cec2022_variants}
  \begin{tabular}{llccc}
    \hline\noalign{\smallskip}
    Function & Metric & LLJC-KEO & HVC-KEO & HCKEO \\
    \noalign{\smallskip}\hline\noalign{\smallskip}
    F1  & Ave & 3.00E+02 & 3.00E+02 & \textbf{3.00E+02} \\
    F2  & Ave & \textbf{402.9083} & 403.6444 & 405.3619 \\
    F3  & Ave & 600.0512 & 600.0005 & \textbf{600.0000} \\
    F4  & Ave & 814.6590 & 813.8962 & \textbf{810.5417} \\
    F5  & Ave & 901.2022 & 900.1083 & \textbf{900.0392} \\
    F6  & Ave & \textbf{1808.0229} & 2053.7615 & 1811.2306 \\
    F7  & Ave & 2012.5819 & \textbf{2008.3568} & 2015.3281 \\
    F8  & Ave & 2220.3490 & 2218.4813 & \textbf{2215.1911} \\
    F9  & Ave & \textbf{2497.7294} & 2501.2772 & 2529.2844 \\
    F10 & Ave & 2522.9880 & 2504.1033 & \textbf{2500.3522} \\
    F11 & Ave & 2646.7467 & \textbf{2615.0442} & 2645.1511 \\
    F12 & Ave & \textbf{2856.1626} & 2857.7005 & 2862.1113 \\
    \noalign{\smallskip}\hline\noalign{\smallskip}
    Friedman average rank & & 2.1667 & 2.0000 & \textbf{1.8333} \\
    Overall rank & & 3 & 2 & \textbf{1} \\
    \noalign{\smallskip}\hline
  \end{tabular}
\end{table}

\paragraph{Result Analysis}

As observed from Table~\ref{tab:cec2022_variants}, HCKEO exhibits the best overall performance, obtaining the lowest Friedman ranking \cite{Thomas2016} value of 1.833, indicating that the combination of both strategies provides complementary advantages.
Although the individual variants (LLJC-KEO and HVC-KEO) exhibit competitive results in certain benchmark functions, the hybrid version (HCKEO) consistently maintains superior average performance across most test functions.
The LLJC component enhances the algorithm’s global exploration ability by introducing stochastic long jumps, while the HVC mechanism promotes efficient local exploitation through structured crossover in both horizontal and vertical directions.
Therefore, the integration of these two strategies significantly improves the balance between exploration and exploitation, leading to enhanced optimization robustness and convergence efficiency.

\subsection{Exploration–Exploitation Analysis}

To investigate the search dynamics of the proposed HCKEO algorithm, the balance between exploration and exploitation is examined on the CEC2022 benchmark suite with 20-dimensional problems.

In population-based optimization algorithms, the search process typically alternates between two complementary behaviors.
One aims to explore new areas of the search space in order to discover potential candidate solutions, while the other concentrates on improving the quality of solutions within promising regions.
An appropriate balance between these two behaviors is essential for achieving both effective global search and reliable convergence.

To quantitatively evaluate this behavior, the dimensional diversity metric proposed by Kashif \emph{et al.}~\cite{Kashif2019} is adopted.
This metric reflects the dispersion of individuals across different dimensions of the population.
Based on the diversity value at each iteration, the proportions of exploration and exploitation during the search process can be estimated as:
\begin{equation}\label{Exploration}
  Exploration(\%) = \frac{{Div(t)}}{{Di{v_{\max }}}} \times 100,
\end{equation}
\begin{equation}\label{Exploitation}
  Exploitation(\%) = \frac{{|Div(t) - Di{v_{\max }}|}}{{Di{v_{\max }}}} \times 100,
\end{equation}
\begin{equation}\label{caldiv}
  Div(t) = \frac{1}{{Dim}}\sum_{d=1}^{Dim} \left( \frac{1}{N}\sum_{i = 1}^N \big| \mathrm{median}(x_d(t)) - x_{id}(t) \big| \right),
\end{equation}
where $x_{id}(t)$ denotes the value of the (d)-th dimension of the $i^{th}$ individual at iteration $t$, $Dim$ is the problem dimension, and $\mathrm{median}(x_d(t)))$represents the median value of the $d^{th}$ dimension across the population.
This formulation measures the average deviation of individuals from the population median, thereby characterizing the distribution of solutions in the search space.

\begin{figure}[!ht]
\centering
\subfloat[F1]{\includegraphics[width=0.32\textwidth]{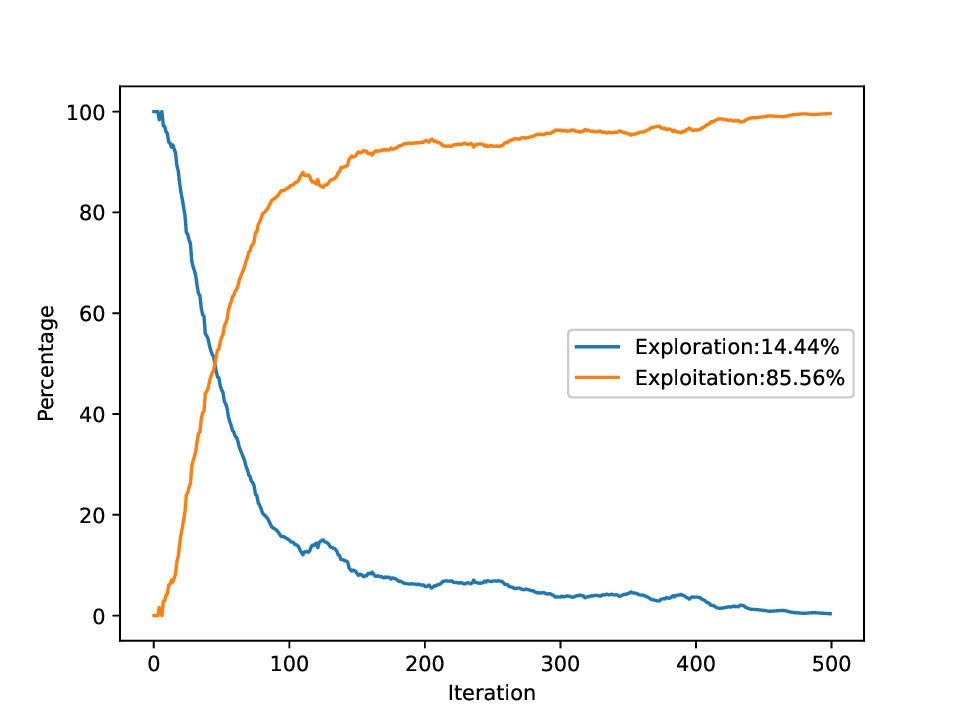}}
\hfill
\subfloat[F2]{\includegraphics[width=0.32\textwidth]{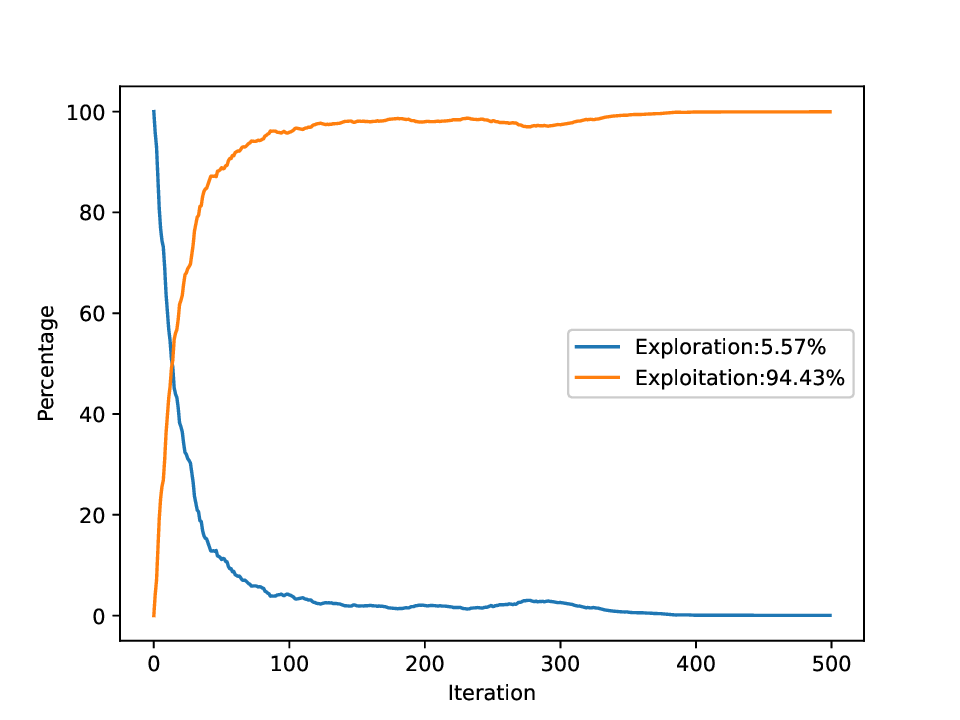}}
\hfill
\subfloat[F3]{\includegraphics[width=0.32\textwidth]{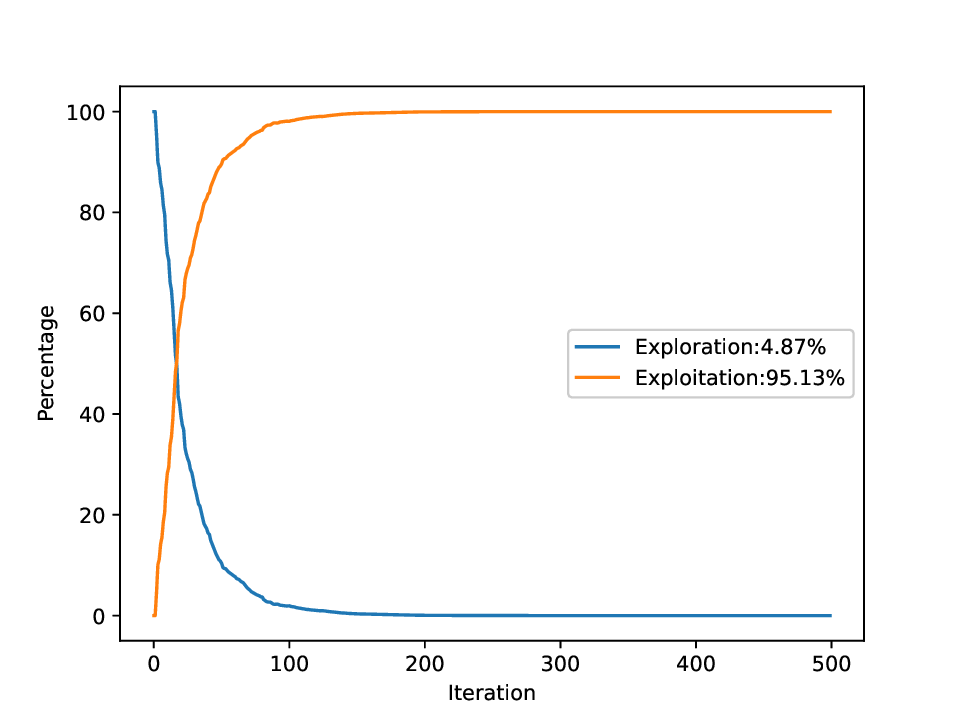}}\\
\subfloat[F4]{\includegraphics[width=0.32\textwidth]{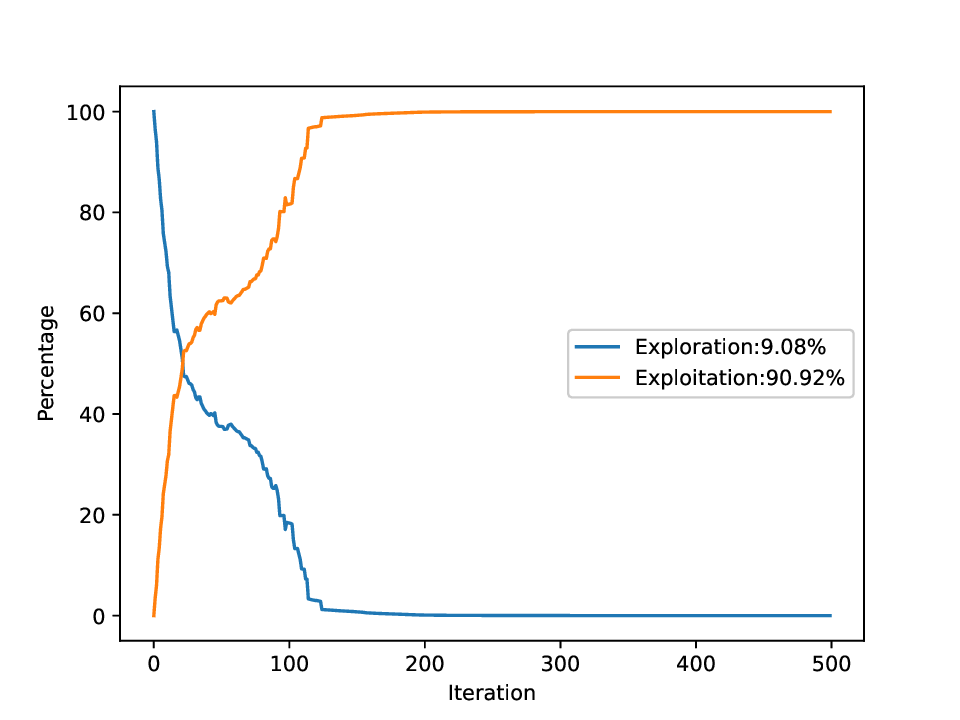}}
\hfill
\subfloat[F5]{\includegraphics[width=0.32\textwidth]{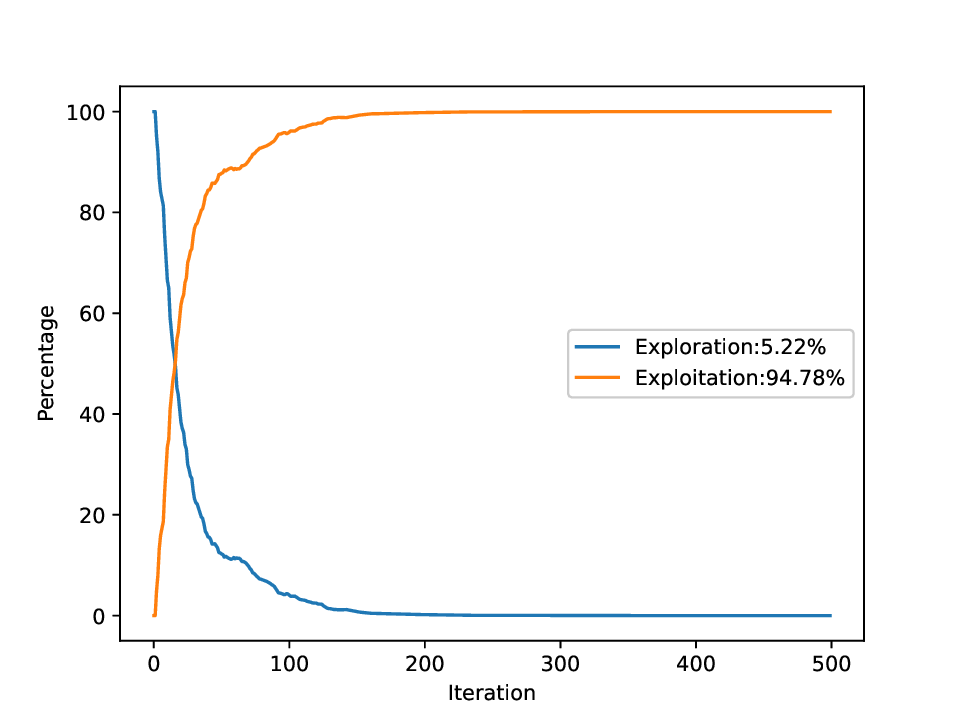}}
\hfill
\subfloat[F6]{\includegraphics[width=0.32\textwidth]{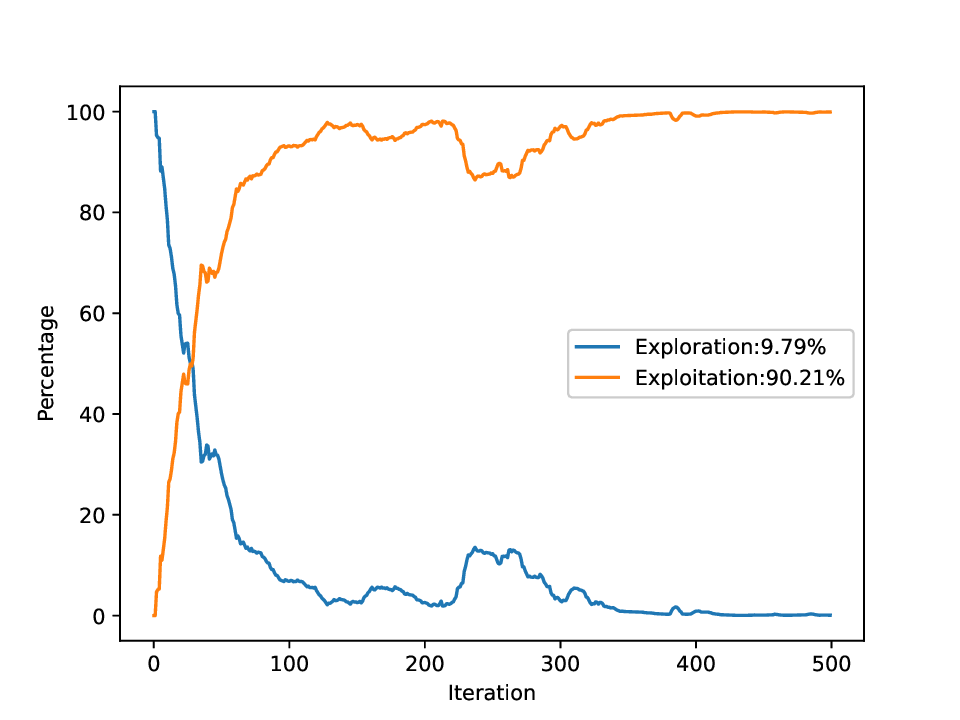}}\\
\subfloat[F7]{\includegraphics[width=0.32\textwidth]{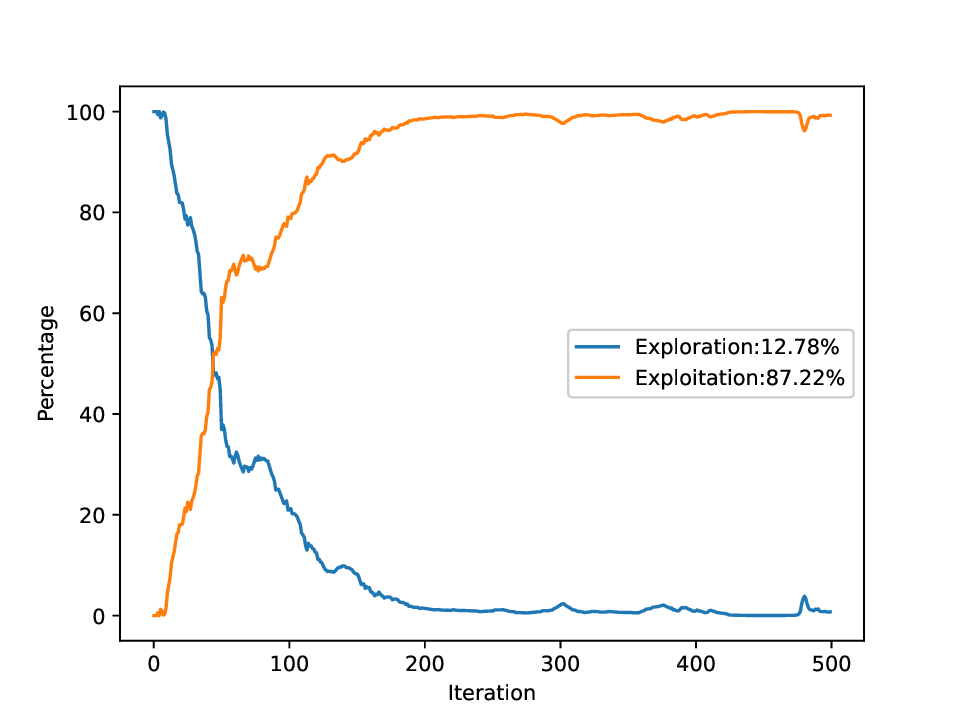}}
\hfill
\subfloat[F8]{\includegraphics[width=0.32\textwidth]{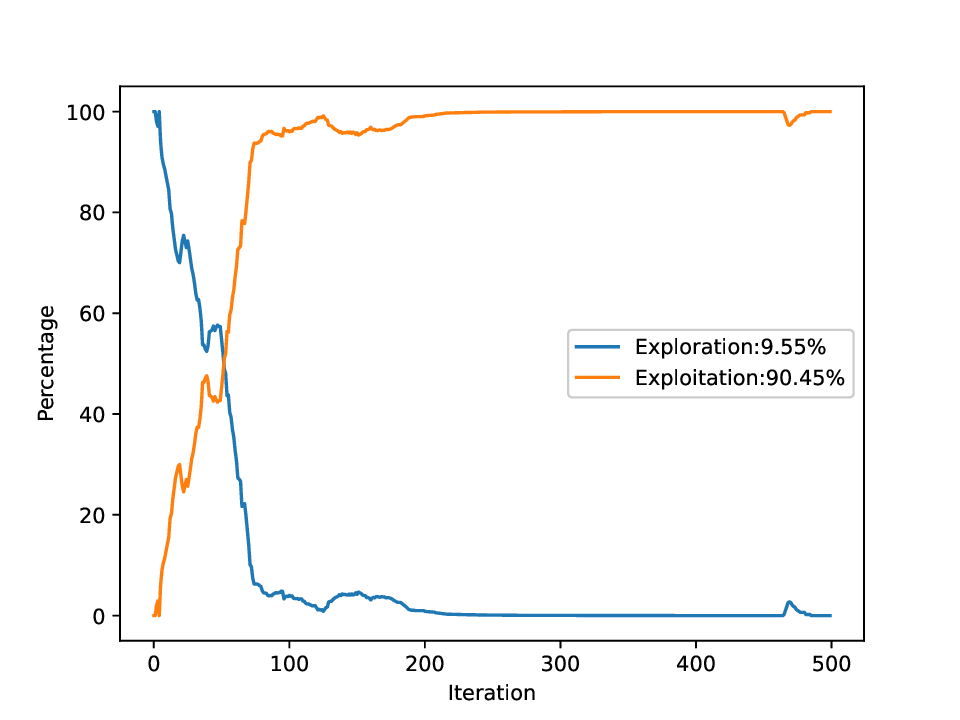}}
\hfill
\subfloat[F9]{\includegraphics[width=0.32\textwidth]{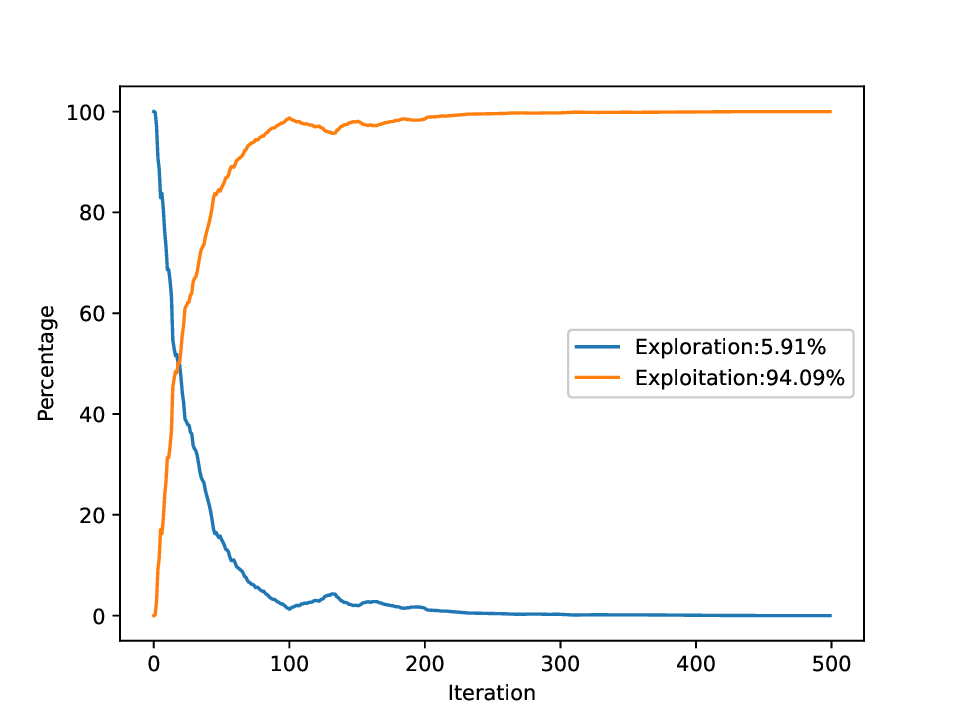}}\\
\subfloat[F10]{\includegraphics[width=0.32\textwidth]{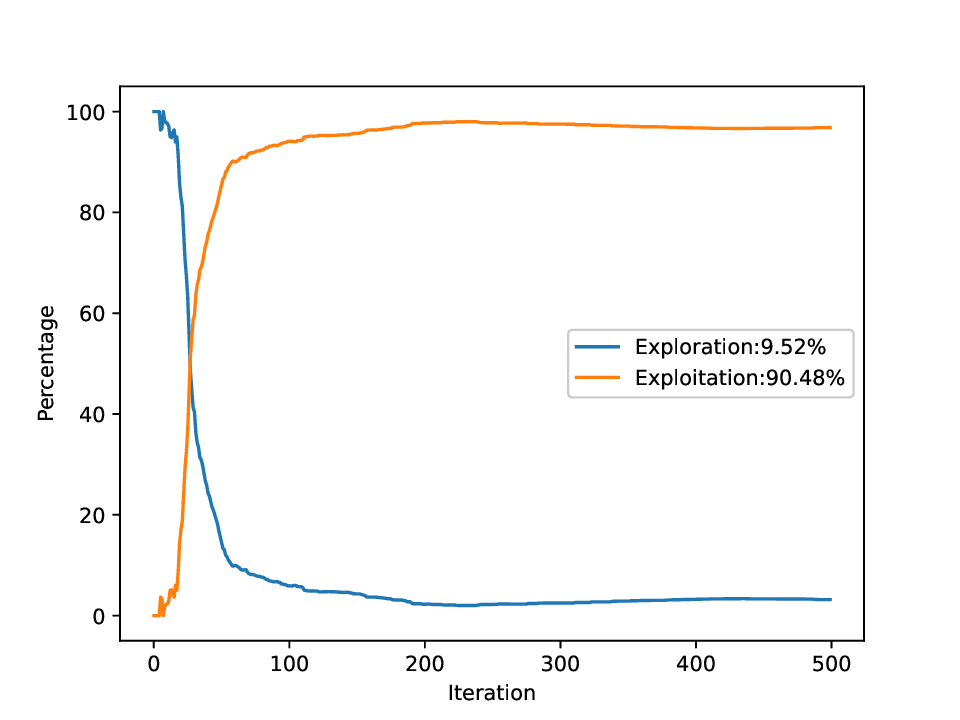}}
\hfill
\subfloat[F11]{\includegraphics[width=0.32\textwidth]{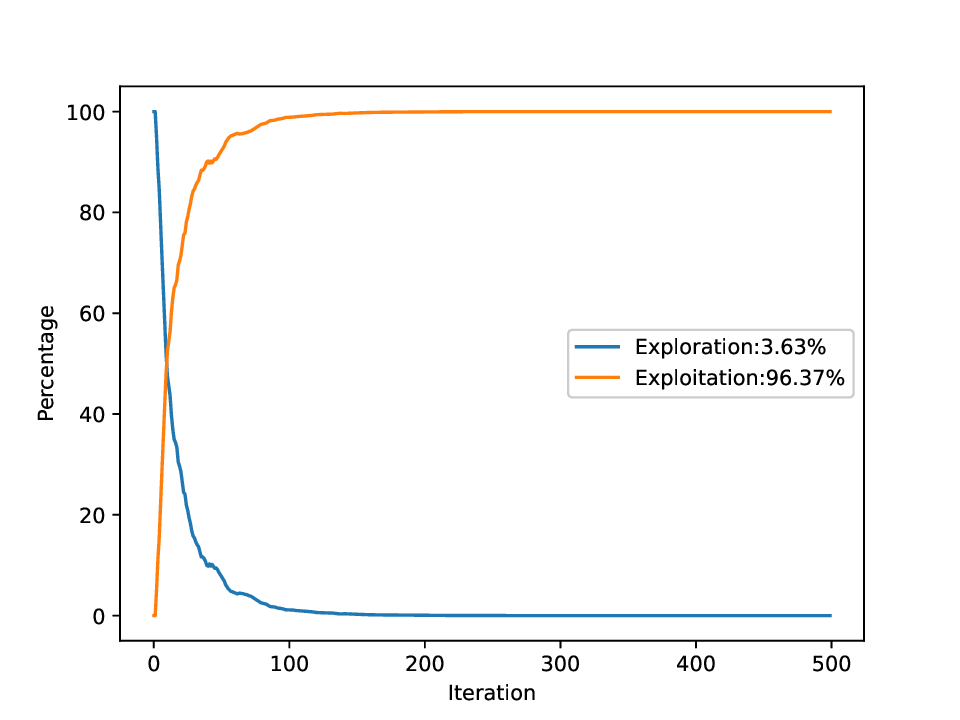}}
\hfill
\subfloat[F12]{\includegraphics[width=0.32\textwidth]{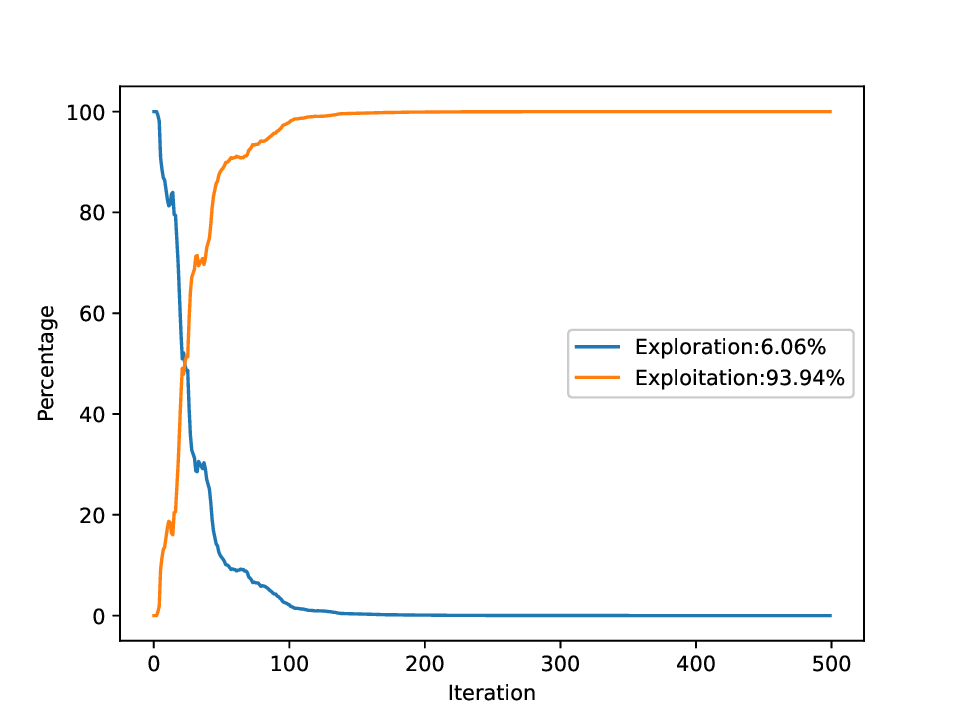}}\\
\caption{Exploration–exploitation balance of HCKEO on selected CEC2022 functions (Dim=20).}
\label{fig:plpt}
\end{figure}

Fig.~\ref{fig:plpt} presents the variation of exploration and exploitation during the optimization process in HCKEO for all benchmark functions.
According to Bernardo \emph{et al.}~\cite{Bernardo2020}, an ideal balance during the search process is often associated with approximately 10\% exploration and 90\% exploitation.
However, rather than maintaining a fixed ratio between exploration and exploitation, an effective optimization process typically involves a gradual transition from exploration-dominated search in the early stages to exploitation-oriented refinement in the later stages.
As can be observed from the figure, the exploration and exploitation curves of HCKEO intersect mainly during the early stages of the optimization, indicating strong exploration ability at the beginning.
In later iterations, the exploration level gradually decreases while exploitation becomes dominant, enabling the algorithm to refine solutions and improve convergence accuracy.
The results indicate that HCKEO is able to maintain a balanced and smooth transition between exploration and exploitation throughout the entire search procedure.

\subsection{Experiments on the CEC2022 Benchmark Suite}

This paper used the CEC2022 benchmark test suite and conducts comparative experiments to verify the optimization capability of HCKEO.
The algorithm was evaluated against several representative optimization methods, including classical algorithms, recently proposed metaheuristic techniques, and the original KEO framework.
In addition to KEO in 2025, the comparison algorithms also encompass classical methods such as Particle Swarm Optimization (PSO)~\cite{Manh2021} and Genetic Algorithm (GA)~\cite{Roberge2013}.
5
State-of-the-art algorithms utilized for comparison include the Starling Murmuration Optimizer (SMO)~\cite{Hoda2022} introduced in 2021, the White Shark Optimizer (WSO)~\cite{Braik2022} developed in 2022, the Snow Ablation Optimizer (SAO)~\cite{Deng2023} from 2023, and the Black-winged Kite Algorithm (BKA)~\cite{Wang2024} released in 2024.
Furthermore, advanced optimizers employed for comparison consist of the Linear Population Size Reduction Differential Evolution with Success-History based Adaptation (LSHADE)~\cite{Tanabe2013} and the Quantum Squirrel Search Algorithm (QSSA)~\cite{Yu2022}.
These algorithms were selected because they represent widely used optimization approaches and recent developments in population-based metaheuristics.

To ensure a fair comparison, the parameter settings of all algorithms follow the configurations reported in their original studies.
The detailed parameter values are summarized in Table~\ref{parameter settings}:
\begin{table}[!ht]
\centering
\small
\caption{Control parameters of the compared algorithms.}
\setlength{\tabcolsep}{4mm}
\begin{tabular}{lll}
\hline
Algorithm & Parameter & Value \\
\hline
GA & Crossover probability ($P_c$) & 0.8 \\
   & Mutation probability ($P_m$) & 0.1 \\

PSO & Cognitive coefficient ($C_1$) & 2 \\
    & Social coefficient ($C_2$) & 2 \\

SMO & $\lambda$ & 20 \\
    & $\mu$ & 0.5 \\

WSO & $f$ & $(0.07,0.075)$ \\
    & $\tau$ & 4.11 \\
    & $p$ & $(0.5,1.5)$ \\
    & $a_0,a_1,a_2$ & $6.25,100,0.0005$ \\

SAO & Parameter setting & Default \\

BKA & $P$ & 0.9 \\

QSSA & $\omega_1,\omega_2$ & $0.6,0.9$ \\
     & $C_1,C_2$ & $0.6,0.9$ \\
     & $P_{dp}$ & 0.2 \\
     & $P$ & 0.5 \\
     & $sf$ & 32.8 \\

LSHADE & $P_{rate}$ & 0.11 \\
       & $Arc_{rate}$ & 1.4 \\

KEO & Parameter setting & Default \\

HCKEO & Parameter setting & Default \\
\hline
\end{tabular}
\label{parameter settings}
\end{table}
The experiments were carried out on the CEC2022 benchmark suite using two commonly adopted problem dimensions, namely (D=10) and (D=20), which allow the algorithms to be evaluated under different levels of search complexity.

\begin{table}[!ht]
\centering
\caption{Statistical performance of the compared algorithms on the CEC2022 benchmark functions (10D).}
\label{tab:cec2022_10dim}
\scriptsize
\setlength{\tabcolsep}{1.5pt}
\begin{tabular}{@{}llcccccccccc@{}}
\toprule
Function & Metric & GA & PSO & SMO & WSO & SAO & BKA & QSSA & LSHADE & KEO & HCKEO \\
\midrule
F1 & Mean & 2683.492 & 300.0011 & 300.1723 & 328.4874 & 301.6204 & 300.4817 & 9769.695 & 300.0803 & 300.0028 & 300.0000 \\
   & Best & 479.3674 & 300.0000 & 300.0001 & 300.0000 & 300.0000 & 300.0211 & 628.2560 & 300.0004 & 300.0000 & 300.0000 \\
   & Worst & 8077.836 & 300.0075 & 302.1188 & 598.0290 & 347.8756 & 303.6512 & 23454.85 & 300.7740 & 300.0325 & 300.0000 \\
   & Std & 1699.888 & 0.001893 & 0.406482 & 73.63505 & 8.736884 & 0.729447 & 6499.785 & 0.148412 & 0.007318 & 6.06E-14 \\
F2 & Mean & 430.2422 & 412.8810 & 407.8281 & 403.2120 & 406.9057 & 407.8432 & 596.3323 & 406.3384 & 405.9879 & 405.9595 \\
   & Best & 400.7336 & 400.0045 & 400.0002 & 400.0000 & 400.0019 & 400.0033 & 402.9028 & 400.0533 & 400.0109 & 400.0000 \\
   & Worst & 475.7380 & 472.2589 & 491.0057 & 470.9449 & 408.9161 & 471.8784 & 1041.769 & 408.9161 & 408.9161 & 408.9161 \\
   & Std & 31.10451 & 20.45451 & 16.10447 & 12.89332 & 3.377306 & 13.23624 & 159.1483 & 2.818537 & 3.240954 & 3.277361 \\
F3 & Mean & 600.4334 & 600.7363 & 600.0376 & 601.0015 & 600.0003 & 619.9046 & 641.2472 & 600.0000 & 600.0808 & 600.0000 \\
   & Best & 600.1861 & 600.0000 & 600.0000 & 600.0010 & 600.0000 & 603.8774 & 614.1057 & 600.0000 & 600.0000 & 600.0000 \\
   & Worst & 600.8509 & 612.1389 & 600.8613 & 603.7398 & 600.0078 & 637.1418 & 670.3612 & 600.0000 & 600.6555 & 600.0000 \\
   & Std & 0.170175 & 2.307255 & 0.160261 & 1.117901 & 0.001418 & 10.80939 & 12.75370 & 1.87E-09 & 0.162208 & 2.59E-07 \\
F4 & Mean & 828.5288 & 813.2329 & 816.7459 & 815.0684 & 812.0968 & 815.3374 & 835.3702 & 817.9605 & 821.1263 & 812.2380 \\
   & Best & 809.0986 & 804.9748 & 806.9647 & 803.2638 & 803.9798 & 805.0406 & 812.8771 & 809.7094 & 809.9496 & 803.9798 \\
   & Worst & 857.8537 & 823.8790 & 833.8285 & 834.2472 & 825.8689 & 828.8661 & 874.6216 & 824.5849 & 836.8134 & 823.8790 \\
   & Std & 11.66958 & 4.006752 & 7.241736 & 9.371571 & 5.070975 & 5.860810 & 16.18016 & 3.925013 & 7.077403 & 5.400205 \\
F5 & Mean & 1109.352 & 900.1100 & 900.8610 & 901.2764 & 900.0422 & 1012.120 & 1376.006 & 900.0000 & 902.0522 & 900.0000 \\
   & Best & 904.5327 & 900.0000 & 900.0000 & 900.0000 & 900.0000 & 905.7471 & 973.5272 & 900.0000 & 900.0000 & 900.0000 \\
   & Worst & 3062.927 & 900.5439 & 912.7252 & 907.8182 & 900.5439 & 1353.681 & 2489.293 & 900.0000 & 915.1903 & 900.0000 \\
   & Std & 413.7847 & 0.158087 & 2.394438 & 2.042119 & 0.127677 & 86.29854 & 293.4764 & 0.000000 & 3.674741 & 0.000000 \\
F6 & Mean & 11284.77 & 3268.917 & 7.53E+05 & 1820.721 & 4788.522 & 1963.865 & 15947.52 & 1820.805 & 1969.261 & 1806.114 \\
   & Best & 2240.471 & 1822.002 & 1812.222 & 1801.518 & 1819.183 & 1859.317 & 1882.410 & 1805.393 & 1804.588 & 1800.101 \\
   & Worst & 50302.40 & 8065.210 & 1.76E+07 & 1850.667 & 7897.379 & 2559.443 & 395048.7 & 1856.456 & 3750.878 & 1830.040 \\
   & Std & 10381.10 & 1840.758 & 3.23E+06 & 13.70603 & 2251.593 & 128.7054 & 71614.91 & 13.03945 & 405.6685 & 6.348536 \\
F7 & Mean & 2014.197 & 2018.783 & 2018.504 & 2017.327 & 2021.465 & 2034.748 & 2074.128 & 2008.942 & 2017.466 & 2012.490 \\
   & Best & 2001.090 & 2001.990 & 2000.176 & 2001.105 & 2000.996 & 2020.641 & 2026.000 & 2001.640 & 2000.000 & 2000.000 \\
   & Worst & 2022.480 & 2030.698 & 2023.982 & 2030.877 & 2031.936 & 2073.558 & 2163.616 & 2021.366 & 2024.598 & 2021.618 \\
   & Std & 8.807140 & 9.799195 & 6.821683 & 10.60948 & 4.476896 & 13.81353 & 33.41550 & 5.752513 & 7.873875 & 9.580332 \\
F8 & Mean & 2219.411 & 2231.809 & 2219.827 & 2222.833 & 2220.939 & 2226.493 & 2232.135 & 2216.039 & 2220.973 & 2215.220 \\
   & Best & 2202.232 & 2200.094 & 2201.952 & 2204.817 & 2218.695 & 2212.321 & 2224.189 & 2206.800 & 2218.579 & 2200.067 \\
   & Worst & 2222.457 & 2341.862 & 2226.940 & 2233.699 & 2224.800 & 2244.714 & 2259.128 & 2223.924 & 2222.349 & 2221.617 \\
   & Std & 4.486108 & 37.07464 & 7.538526 & 6.837663 & 1.022231 & 5.625041 & 8.493838 & 5.969250 & 0.763433 & 8.837725 \\
F9 & Mean & 2534.926 & 2529.321 & 2505.369 & 2529.445 & 2529.284 & 2539.081 & 2670.452 & 2529.284 & 2529.284 & 2529.284 \\
   & Best & 2529.461 & 2529.284 & 2493.347 & 2529.287 & 2529.284 & 2529.284 & 2529.821 & 2529.284 & 2529.284 & 2529.284 \\
   & Worst & 2548.099 & 2529.840 & 2541.847 & 2530.506 & 2529.284 & 2676.218 & 2787.655 & 2529.284 & 2529.284 & 2529.284 \\
   & Std & 4.706392 & 0.141063 & 8.688296 & 0.270915 & 0.000000 & 37.27778 & 59.73417 & 0.000000 & 8.33E-12 & 0.000000 \\
F10 & Mean & 2525.437 & 2553.617 & 2523.155 & 2525.672 & 2538.155 & 2571.479 & 2604.392 & 2526.836 & 2523.993 & 2500.359 \\
   & Best & 2400.540 & 2500.162 & 2500.261 & 2500.210 & 2500.240 & 2500.500 & 2500.621 & 2500.269 & 2500.342 & 2500.235 \\
   & Worst & 2636.947 & 2626.143 & 2621.868 & 2609.631 & 2621.358 & 3260.697 & 3373.651 & 2621.104 & 2624.679 & 2500.549 \\
   & Std & 80.08988 & 58.04875 & 46.31690 & 46.45799 & 54.44259 & 143.4459 & 164.0188 & 48.86223 & 47.92314 & 0.071315 \\
F11 & Mean & 2746.662 & 2706.223 & 2691.884 & 2689.163 & 2670.137 & 2685.907 & 3025.557 & 2615.014 & 2663.522 & 2640.137 \\
   & Best & 2604.554 & 2600.000 & 2600.000 & 2600.004 & 2600.000 & 2600.462 & 2757.086 & 2600.000 & 2600.000 & 2600.000 \\
   & Worst & 3184.174 & 3176.586 & 3000.100 & 3000.001 & 2900.000 & 3183.986 & 3562.037 & 2900.000 & 3000.000 & 2751.075 \\
   & Std & 152.1320 & 183.0775 & 111.6070 & 130.3398 & 102.3279 & 156.8790 & 267.5875 & 60.41966 & 95.70967 & 67.69738 \\
F12 & Mean & 2868.155 & 2872.268 & 2856.592 & 2868.734 & 2864.155 & 2865.067 & 2904.812 & 2862.674 & 2863.780 & 2858.618 \\
   & Best & 2863.267 & 2863.851 & 2852.572 & 2862.700 & 2861.435 & 2862.287 & 2865.133 & 2858.618 & 2859.682 & 2751.075 \\
   & Worst & 3184.174 & 3176.586 & 2866.985 & 3000.001 & 2900.000 & 3183.986 & 3102.055 & 2900.000 & 3000.000 & 2862.328 \\
   & Std & 2.654030 & 13.14142 & 3.395383 & 11.20498 & 1.411175 & 1.952986 & 55.47136 & 1.657132 & 1.396138 & 1.401152 \\
   & (L/T/W) & (0/0/12) & (0/0/12) & (0/0/12) & (2/0/10) & (2/0/10) & (0/0/12) & (0/2/10) & (1/2/9) & (0/0/12) &  \\
\bottomrule
\end{tabular}
\end{table}

\begin{table}[!ht]
\centering
\caption{Statistical performance of the compared algorithms on the CEC2022 benchmark functions (20D).}
\label{tab:cec2022_20dim}
\scriptsize
\setlength{\tabcolsep}{2pt}
\begin{tabular}{llcccccccccc}
\toprule
Function & Metric & GA & PSO & SMO & WSO & SAO & BKA & QSSA & LSHADE & KEO & HCKEO \\
\midrule
F1 & Mean & 10945.54 & 971.2756 & 10914.11 & 9060.090 & 17432.06 & 1626.722 & 40847.59 & 9803.456 & 1785.563 & 301.7422 \\
   & Best & 3451.392 & 464.3323 & 1862.353 & 2617.096 & 6822.308 & 525.0722 & 24187.42 & 5044.454 & 794.6222 & 300.0485 \\
   & Worst & 22373.28 & 2017.952 & 34185.71 & 20239.26 & 33876.94 & 3544.859 & 59426.11 & 17704.83 & 4151.324 & 306.7600 \\
   & Std & 4748.447 & 379.0005 & 7763.654 & 4567.043 & 7259.164 & 764.3025 & 10514.24 & 3097.586 & 933.6172 & 2.055797 \\
F2 & Mean & 467.8963 & 456.6098 & 466.8967 & 520.2429 & 450.5928 & 484.9500 & 1453.752 & 448.6656 & 455.6984 & 456.5544 \\
   & Best & 439.0728 & 411.0092 & 434.7607 & 473.6702 & 429.0133 & 425.7750 & 824.3681 & 444.8955 & 402.9115 & 444.8955 \\
   & Worst & 570.8365 & 553.5791 & 553.5909 & 608.3047 & 472.1735 & 596.4861 & 2272.146 & 449.0845 & 504.5142 & 475.0196 \\
   & Std & 25.98132 & 24.98771 & 28.40320 & 38.59635 & 7.794733 & 37.54182 & 388.9540 & 1.278186 & 19.34183 & 11.49435 \\
F3 & Mean & 601.4026 & 602.9482 & 601.8135 & 627.7275 & 600.0392 & 646.2505 & 665.9026 & 600.0000 & 602.6907 & 600.0056 \\
   & Best & 600.6009 & 600.0296 & 600.1947 & 608.4489 & 600.0000 & 632.8501 & 639.3851 & 600.0000 & 600.2740 & 600.0001 \\
   & Worst & 602.0737 & 627.9980 & 609.2663 & 650.7990 & 600.5770 & 660.7498 & 695.4861 & 600.0000 & 613.0082 & 600.0623 \\
   & Std & 0.410557 & 5.299602 & 1.847030 & 9.058023 & 0.144835 & 7.531078 & 13.25208 & 7.76E-06 & 2.723341 & 0.014175 \\
F4 & Mean & 873.2619 & 850.7535 & 845.1869 & 848.8149 & 842.7550 & 871.9826 & 930.3689 & 888.7568 & 852.6664 & 831.5187 \\
   & Best & 833.5051 & 826.8639 & 828.9117 & 817.3690 & 814.9244 & 838.3329 & 887.5014 & 865.9185 & 828.8538 & 812.9345 \\
   & Worst & 935.1136 & 898.5010 & 883.5987 & 913.8315 & 934.8497 & 895.1019 & 1024.855 & 904.9961 & 884.5712 & 852.7327 \\
   & Std & 22.41393 & 15.73614 & 11.68080 & 21.11926 & 32.60343 & 15.51227 & 30.11453 & 9.609513 & 15.35370 & 10.40263 \\
F5 & Mean & 1943.649 & 907.0274 & 1049.664 & 1521.799 & 900.6166 & 1959.114 & 3461.775 & 900.0030 & 1210.476 & 904.1335 \\
   & Best & 1038.961 & 900.1004 & 908.5853 & 995.7836 & 900.0000 & 1477.613 & 1911.945 & 900.0000 & 911.7124 & 900.0000 \\
   & Worst & 2786.740 & 932.6981 & 1457.254 & 2973.184 & 904.2687 & 2500.972 & 5425.717 & 900.0895 & 2324.640 & 936.1301 \\
   & Std & 536.9347 & 6.355236 & 144.6286 & 422.7810 & 1.074520 & 240.7493 & 919.3704 & 0.016346 & 313.8326 & 6.691075 \\
F6 & Mean & 150870.7 & 4344.141 & 21277.57 & 2030.931 & 3707.810 & 7277.250 & 3.69E+08 & 153121.9 & 4913.264 & 3669.385 \\
   & Best & 7350.626 & 1853.130 & 1995.864 & 1854.048 & 2002.609 & 2057.705 & 2261.017 & 40520.20 & 1892.537 & 1854.100 \\
   & Worst & 673596.5 & 13828.79 & 332578.8 & 2462.086 & 9077.567 & 25353.86 & 2.23E+09 & 424102.1 & 18105.32 & 14160.21 \\
   & Std & 143865.3 & 3308.016 & 60321.66 & 160.8989 & 1579.844 & 6664.396 & 4.92E+08 & 91386.23 & 3880.783 & 2972.937 \\
F7 & Mean & 2070.231 & 2077.439 & 2050.973 & 2068.845 & 2051.642 & 2115.401 & 2213.394 & 2059.661 & 2056.175 & 2049.516 \\
   & Best & 2029.754 & 2031.725 & 2023.187 & 2028.620 & 2025.290 & 2053.918 & 2108.151 & 2039.931 & 2026.049 & 2021.309 \\
   & Worst & 2186.075 & 2232.582 & 2120.014 & 2138.518 & 2175.875 & 2172.842 & 2384.581 & 2083.133 & 2145.447 & 2166.803 \\
   & Std & 33.40676 & 47.21851 & 23.11200 & 30.85958 & 27.86789 & 32.96045 & 60.61582 & 9.180605 & 24.13267 & 28.04217 \\
F8 & Mean & 2227.508 & 2274.137 & 2229.003 & 2236.369 & 2228.441 & 2282.436 & 2302.004 & 2234.832 & 2233.329 & 2221.915 \\
   & Best & 2222.861 & 2221.597 & 2221.750 & 2227.058 & 2220.416 & 2226.364 & 2229.994 & 2228.851 & 2221.073 & 2220.537 \\
   & Worst & 2241.329 & 2462.379 & 2242.439 & 2347.568 & 2339.294 & 2473.332 & 2462.100 & 2238.894 & 2356.458 & 2237.549 \\
   & Std & 6.718385 & 73.96457 & 6.000294 & 21.27199 & 21.82776 & 74.70743 & 74.60914 & 2.233858 & 31.92022 & 2.990461 \\
F9 & Mean & 2486.541 & 2493.727 & 2482.437 & 2508.725 & 2480.783 & 2490.756 & 2896.628 & 2480.781 & 2480.806 & 2480.781 \\
   & Best & 2482.223 & 2480.782 & 2473.808 & 2486.606 & 2480.781 & 2481.328 & 2669.696 & 2480.781 & 2480.781 & 2480.781 \\
   & Worst & 2493.840 & 2549.062 & 2508.599 & 2547.784 & 2480.803 & 2533.682 & 3252.109 & 2480.781 & 2480.948 & 2480.781 \\
   & Std & 3.681514 & 20.96695 & 8.475978 & 14.02898 & 0.004557 & 11.27802 & 157.8898 & 2.87E-11 & 0.041477 & 5.97E-07 \\
F10 & Mean & 2437.508 & 3062.900 & 2872.348 & 2815.575 & 2993.435 & 4077.370 & 5016.947 & 2666.355 & 2655.501 & 2530.986 \\
   & Best & 2404.446 & 2500.405 & 2500.468 & 2500.889 & 2500.528 & 2501.061 & 2509.998 & 2500.380 & 2500.473 & 2403.891 \\
   & Worst & 2669.909 & 3775.970 & 3723.513 & 6886.930 & 4421.093 & 5548.350 & 6596.323 & 4535.105 & 3146.860 & 3857.017 \\
   & Std & 57.18316 & 462.9929 & 401.3427 & 997.7101 & 632.4629 & 935.0873 & 1439.041 & 492.3086 & 231.0757 & 254.5261 \\
F11 & Mean & 3168.026 & 3055.325 & 2979.358 & 3457.643 & 2893.345 & 3176.766 & 6359.105 & 2906.679 & 2941.204 & 2916.667 \\
   & Best & 2663.937 & 2600.107 & 2604.762 & 3050.253 & 2600.000 & 2630.716 & 4747.253 & 2900.000 & 2600.120 & 2900.000 \\
   & Worst & 5921.992 & 3565.915 & 3388.977 & 3977.275 & 3000.351 & 3842.202 & 8492.389 & 3000.365 & 3360.572 & 3000.000 \\
   & Std & 711.0581 & 257.0948 & 191.0623 & 207.1374 & 86.84942 & 290.9062 & 928.5483 & 25.41715 & 112.9617 & 37.90491 \\
F12 & Mean & 2973.267 & 2983.370 & 2900.005 & 3028.017 & 2949.316 & 3027.433 & 3162.235 & 2939.420 & 2956.337 & 2943.224 \\
   & Best & 2945.637 & 2945.623 & 2900.004 & 2965.101 & 2932.627 & 2955.125 & 3008.565 & 2932.008 & 2939.866 & 2933.974 \\
   & Worst & 3006.194 & 3077.369 & 2900.005 & 3113.293 & 2993.027 & 3243.460 & 3487.975 & 2950.508 & 3002.472 & 2965.384 \\
   & Std & 16.53742 & 31.67861 & 0.000137 & 47.34230 & 12.31073 & 65.33425 & 126.1453 & 5.607964 & 14.29076 & 6.067805 \\
   & (L/T/W) & (1/0/11) & (0/0/12) & (1/0/11) & (1/0/11) & (2/0/10) & (0/0/12) & (0/0/12) & (5/1/6) & (1/0/11) &  \\
\bottomrule
\end{tabular}
\end{table}

\begin{figure}[!ht]
\centering
\subfloat[F1(D=10)]{\includegraphics[width=0.32\textwidth]{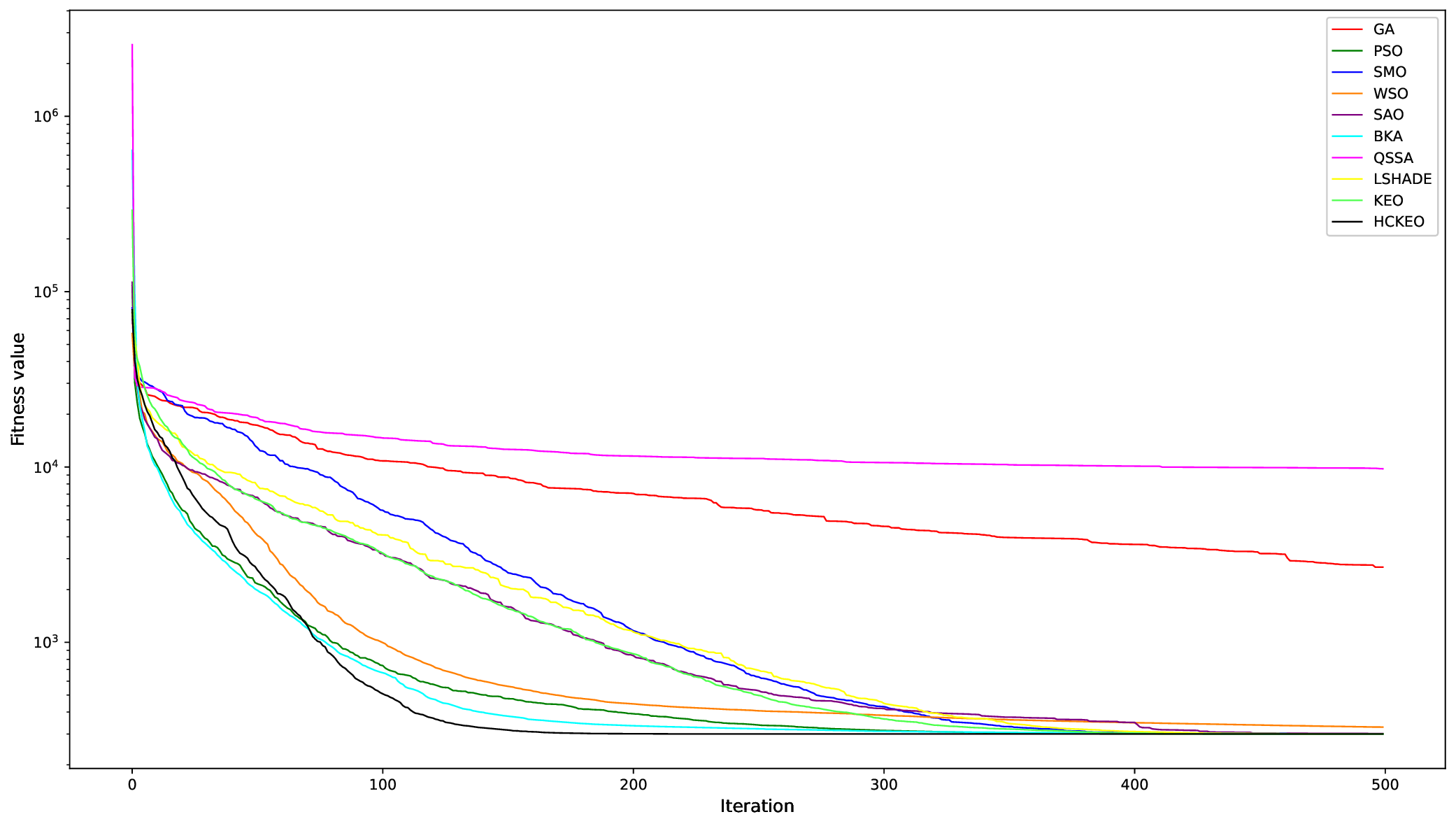}}
\subfloat[F2(D=10)]{\includegraphics[width=0.32\textwidth]{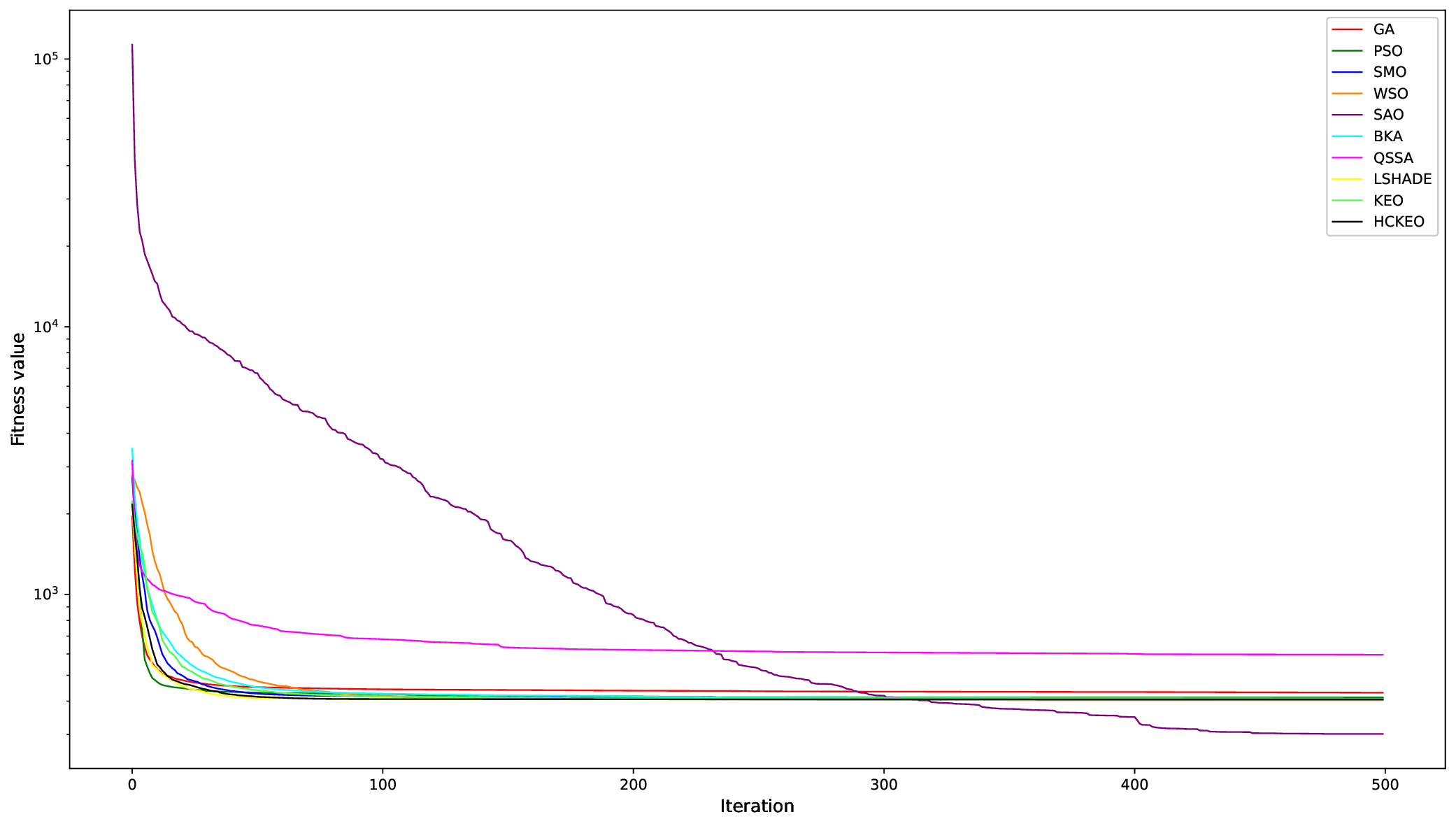}}
\subfloat[F3(D=10)]{\includegraphics[width=0.32\textwidth]{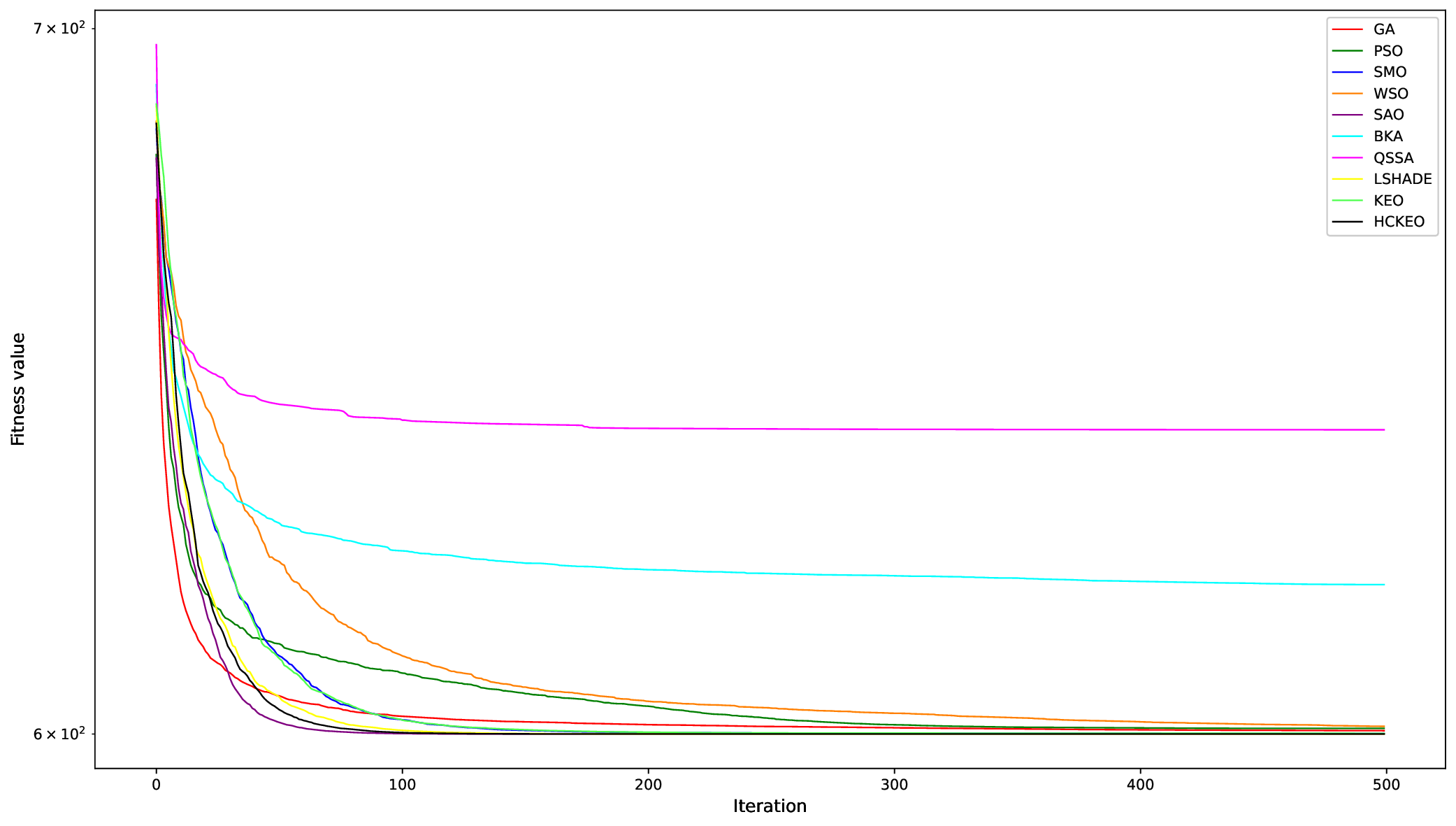}}\\
\subfloat[F4(D=10)]{\includegraphics[width=0.32\textwidth]{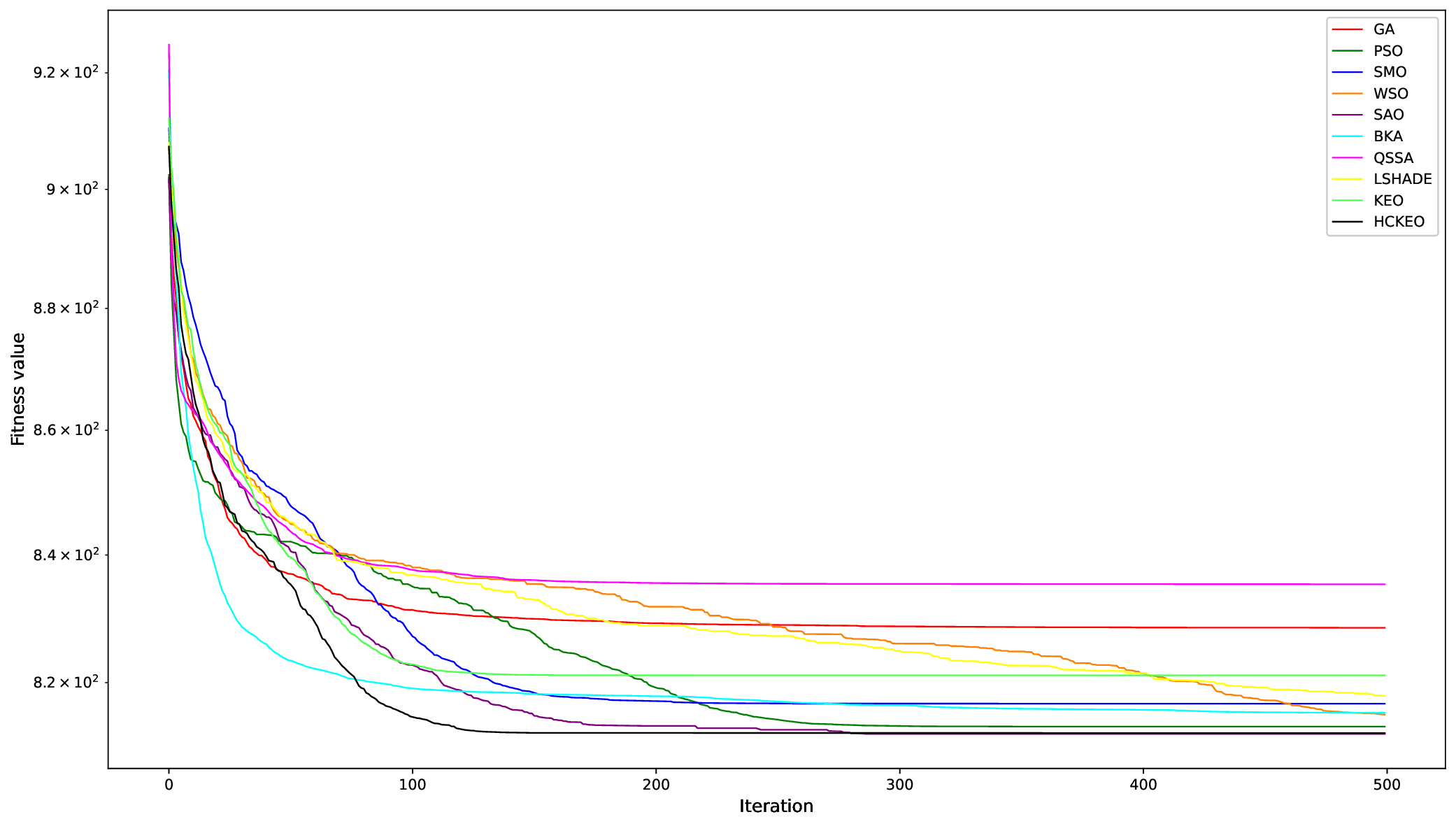}}
\subfloat[F5(D=10)]{\includegraphics[width=0.32\textwidth]{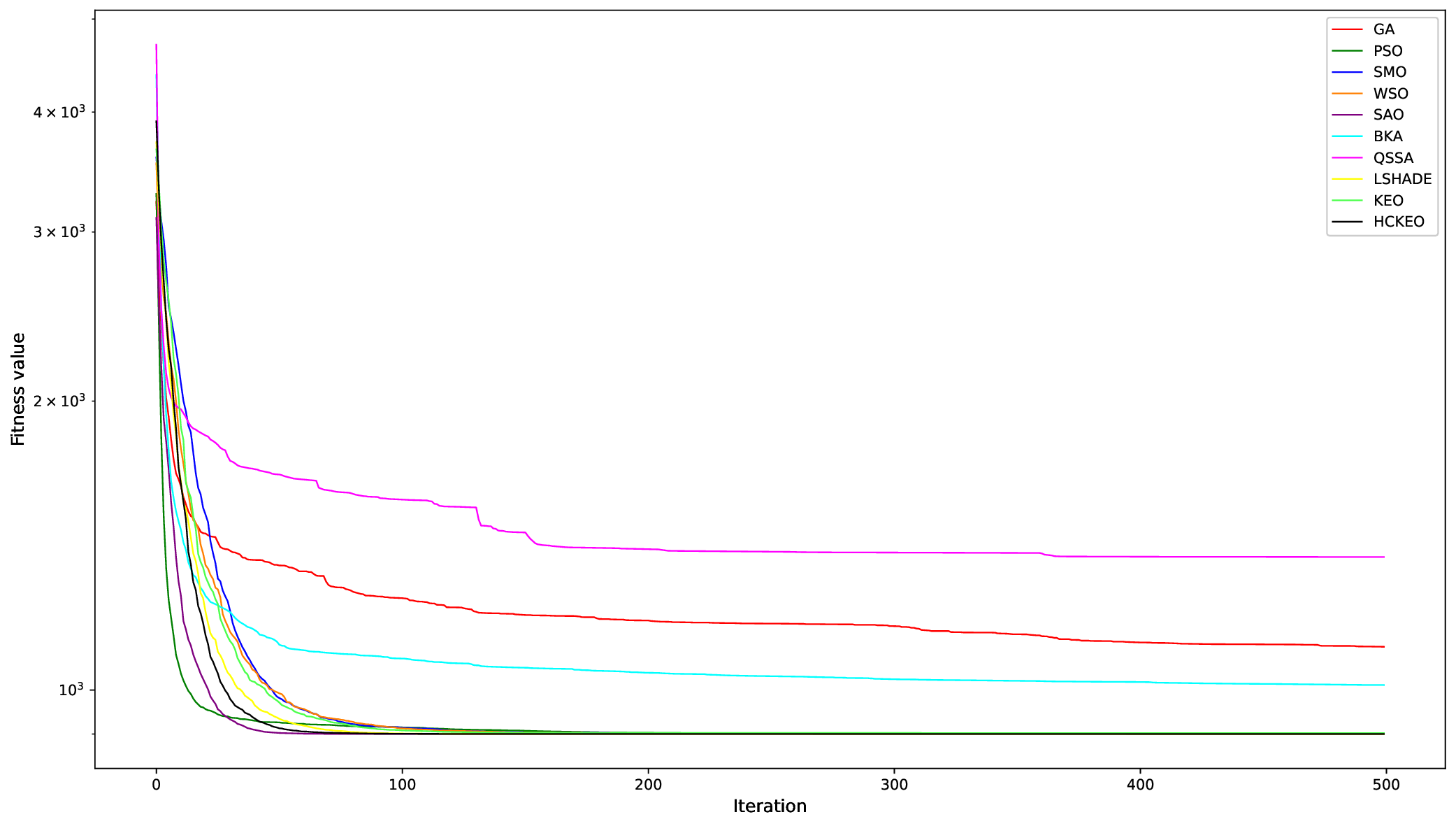}}
\subfloat[F6(D=10)]{\includegraphics[width=0.32\textwidth]{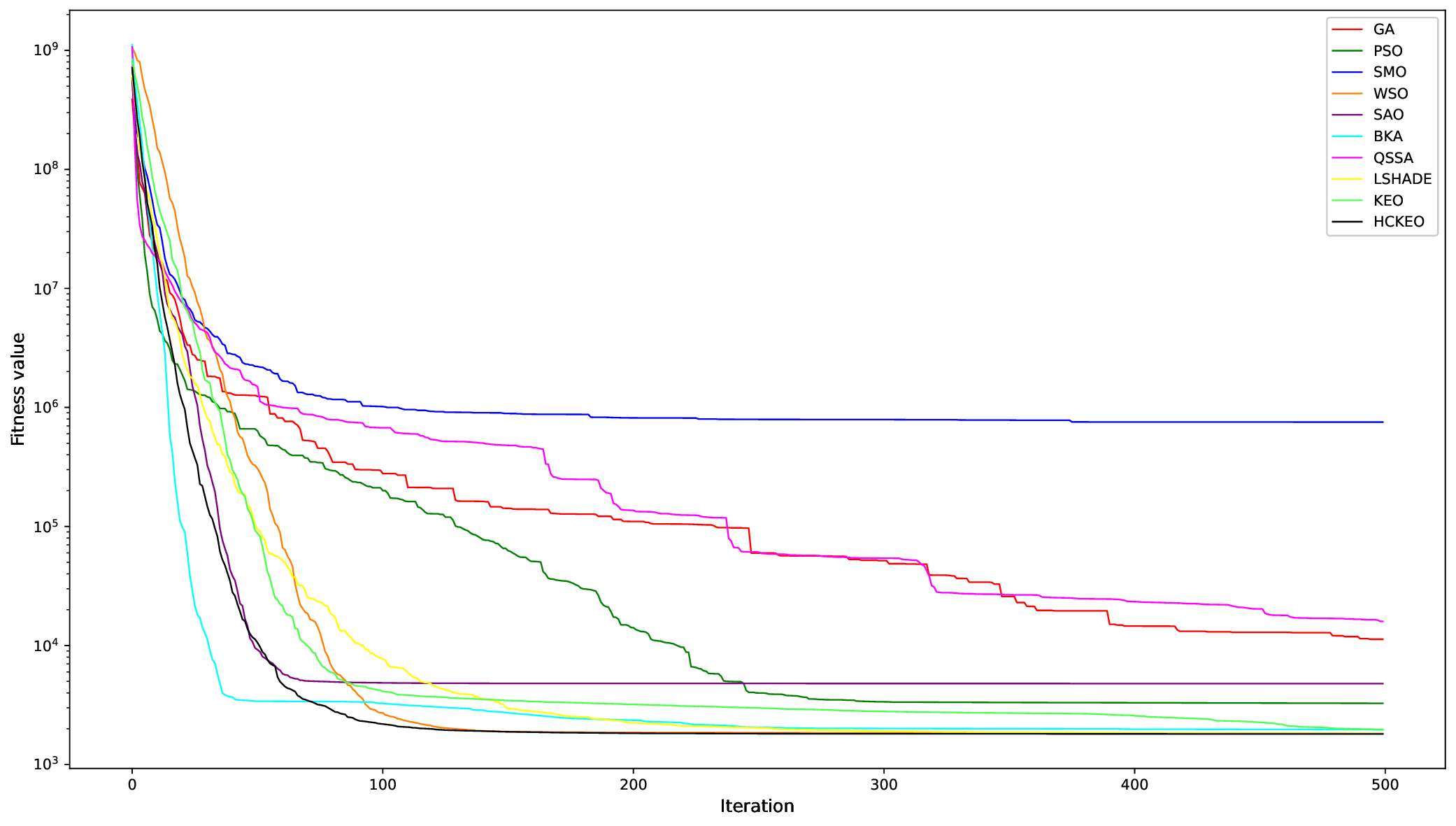}}\\
\subfloat[F7(D=10)]{\includegraphics[width=0.32\textwidth]{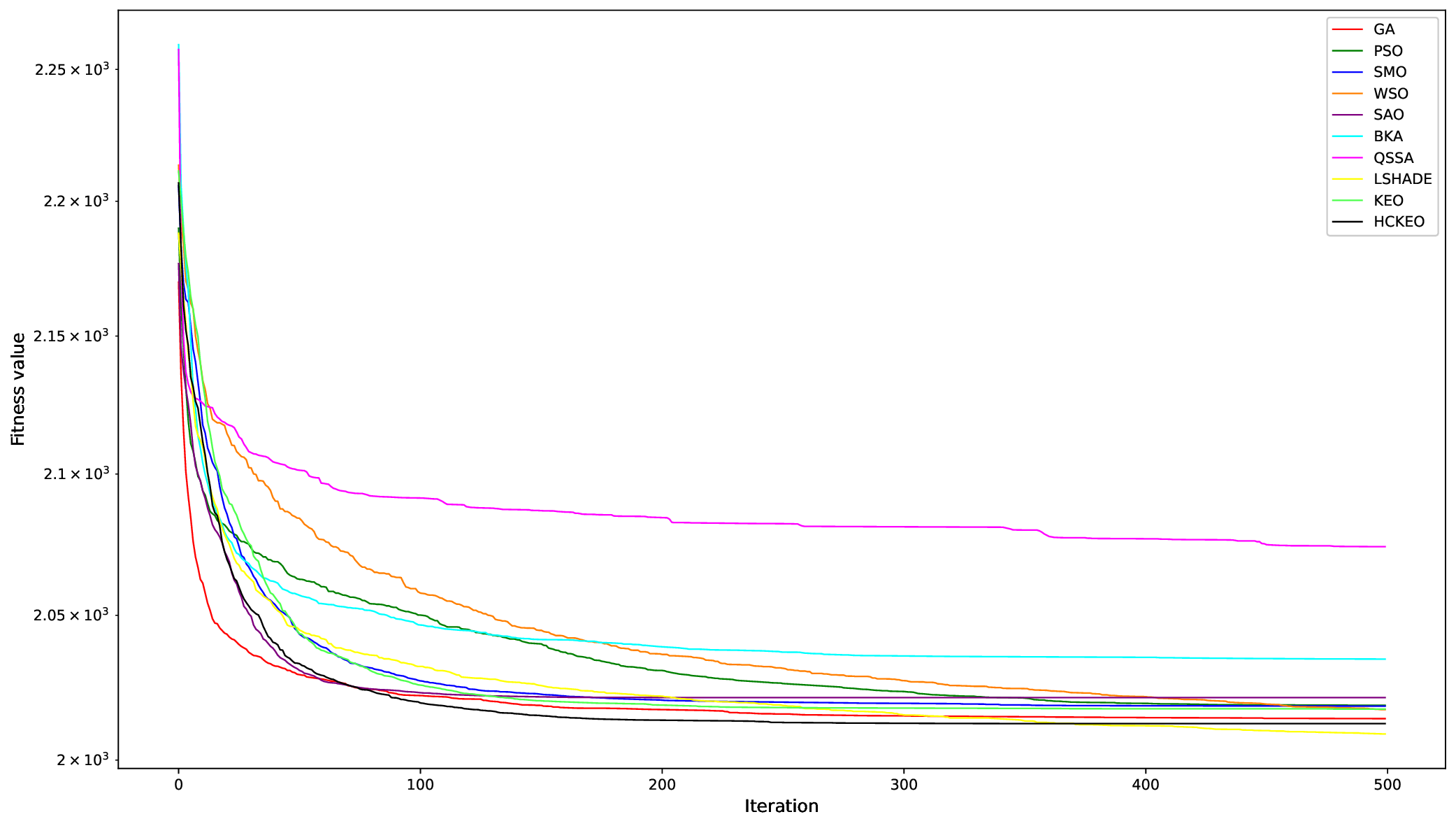}}
\subfloat[F8(D=10)]{\includegraphics[width=0.32\textwidth]{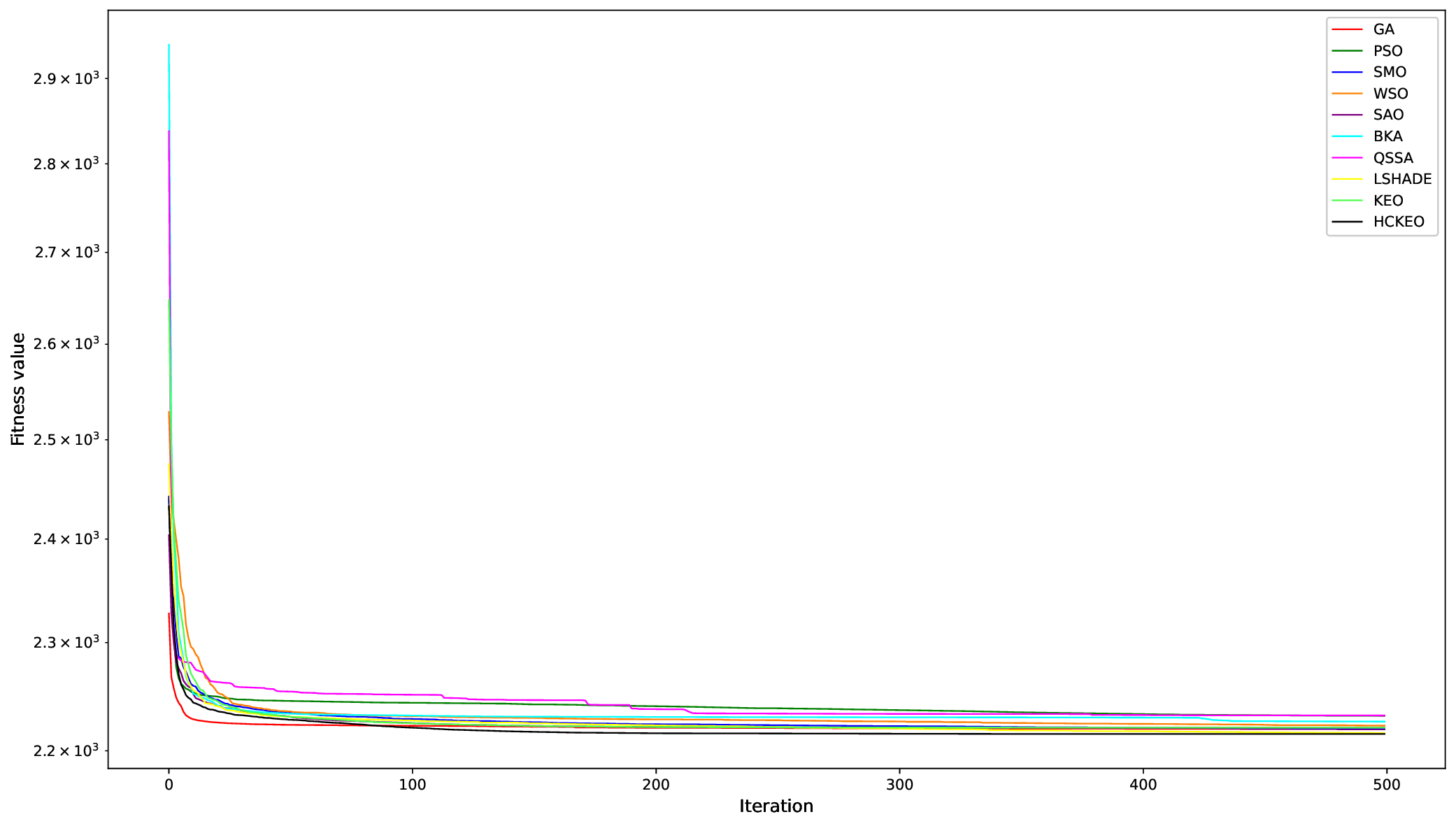}}
\subfloat[F9(D=10)]{\includegraphics[width=0.32\textwidth]{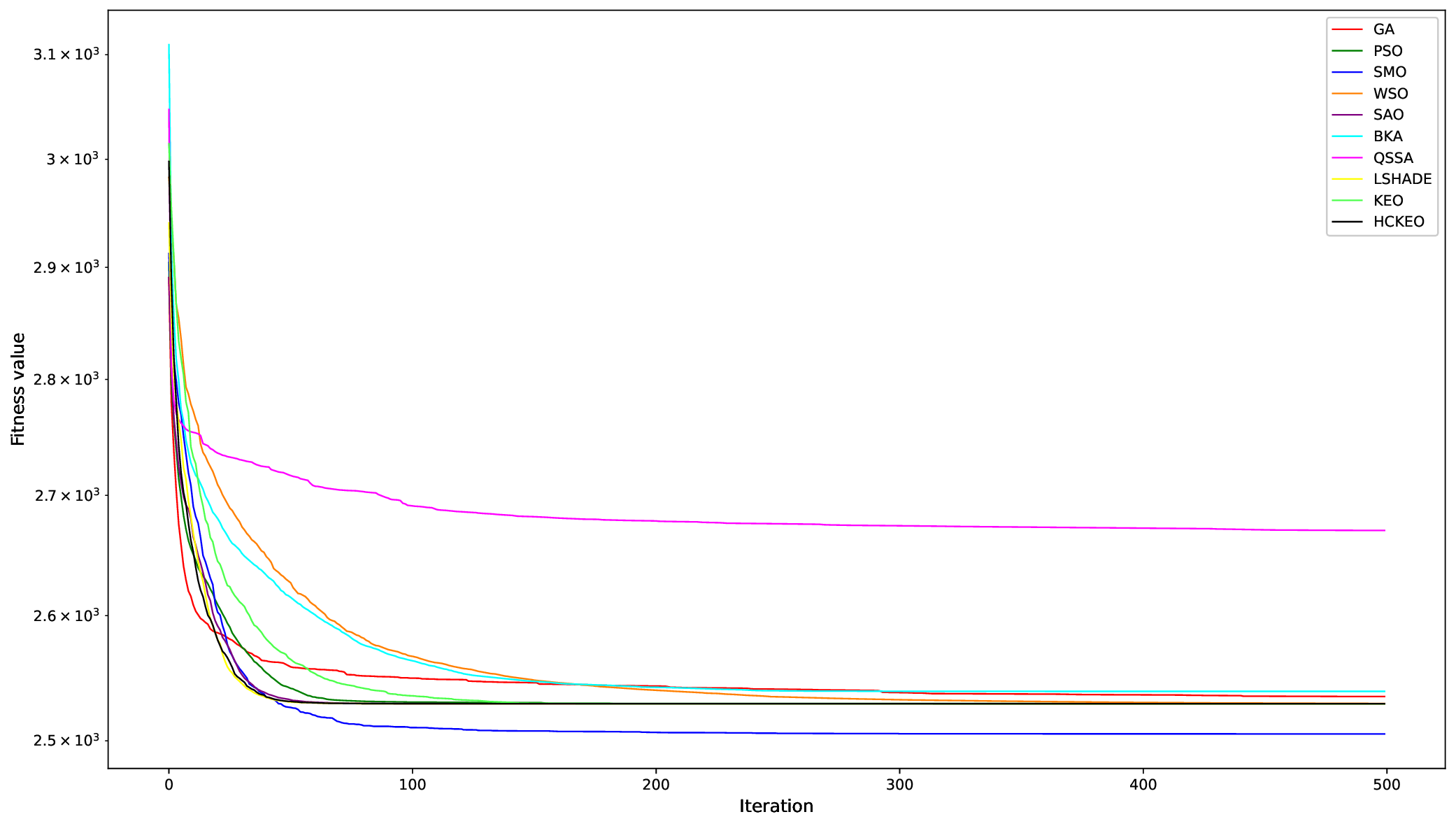}}\\
\subfloat[F10(D=10)]{\includegraphics[width=0.32\textwidth]{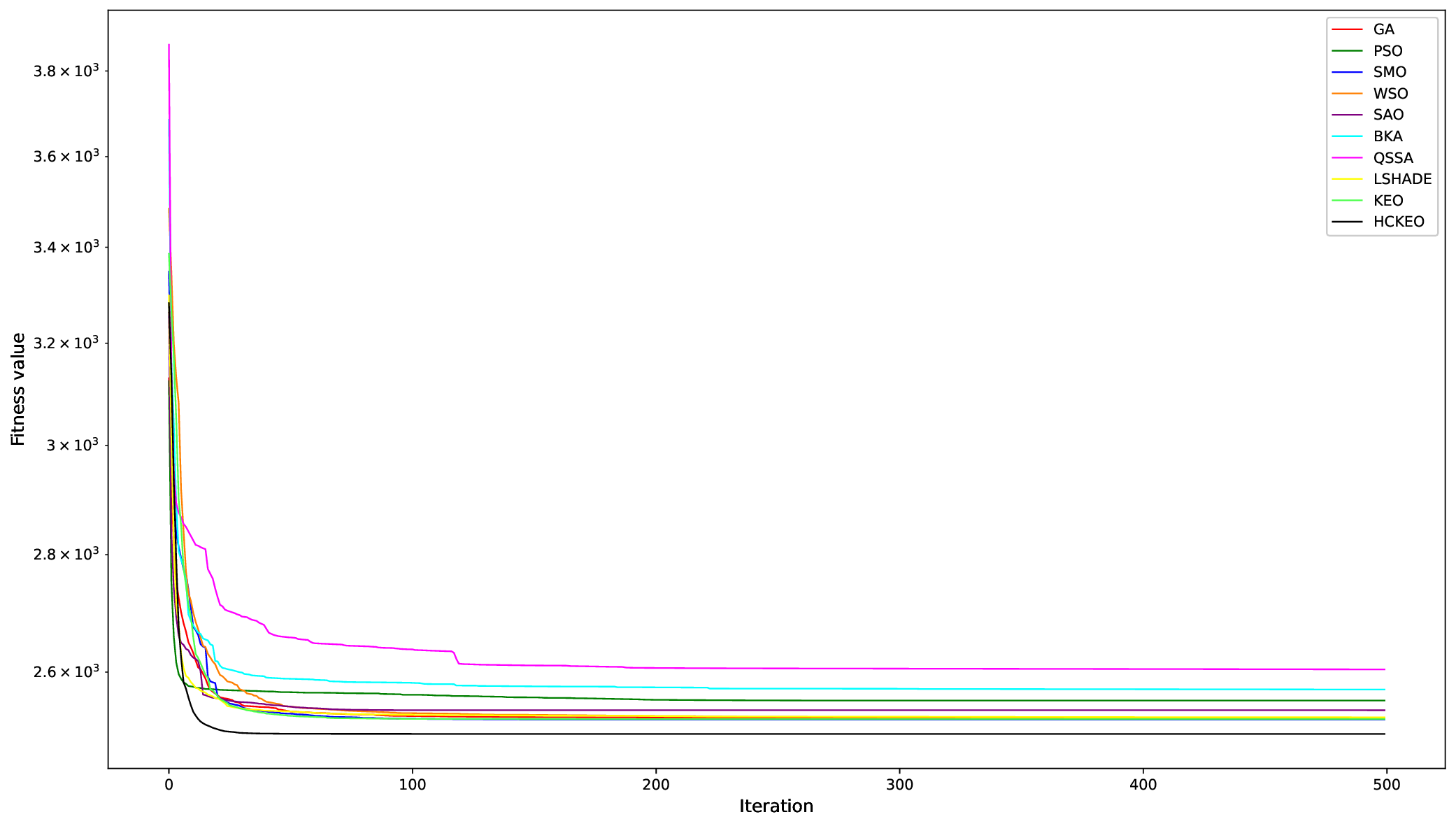}}
\subfloat[F11(D=10)]{\includegraphics[width=0.32\textwidth]{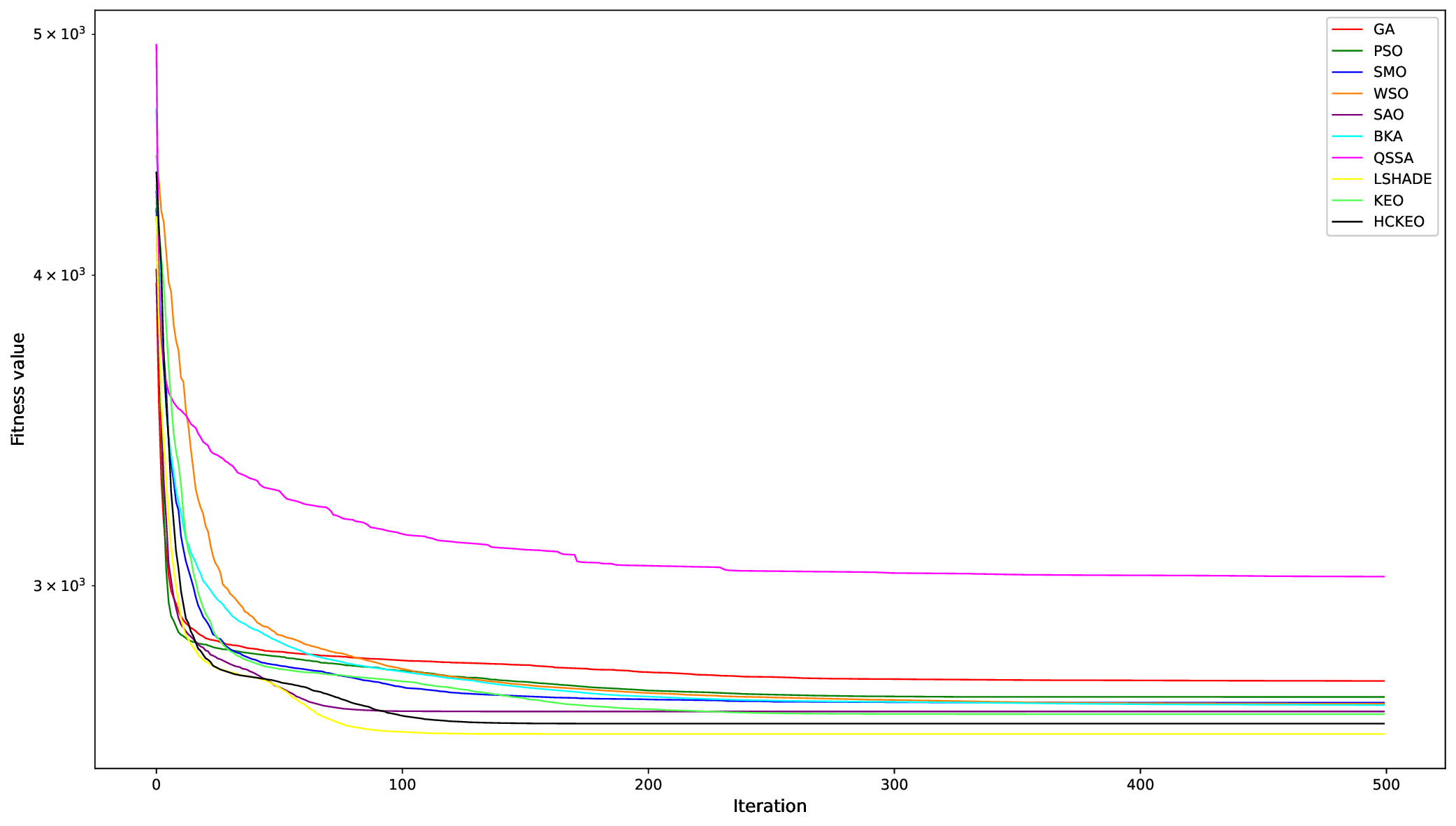}}
\subfloat[F12(D=10)]{\includegraphics[width=0.32\textwidth]{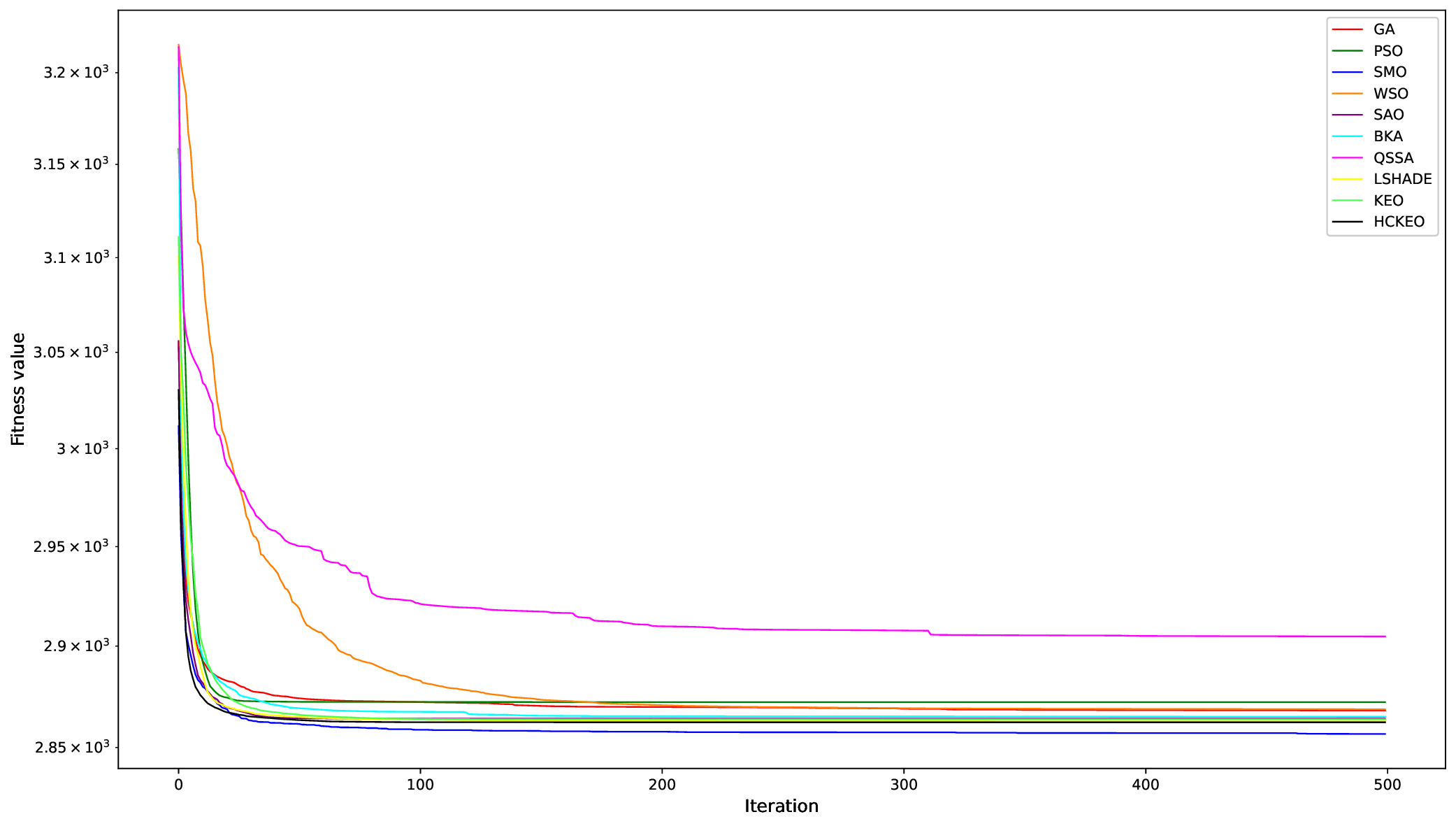}}\\
\caption{CEC2022 test functions convergence curve (Dim=10).}
\label{fig:curve10}
\end{figure}

\begin{figure}[!ht]
\centering
\subfloat[F1(D=20)]{\includegraphics[width=0.32\textwidth]{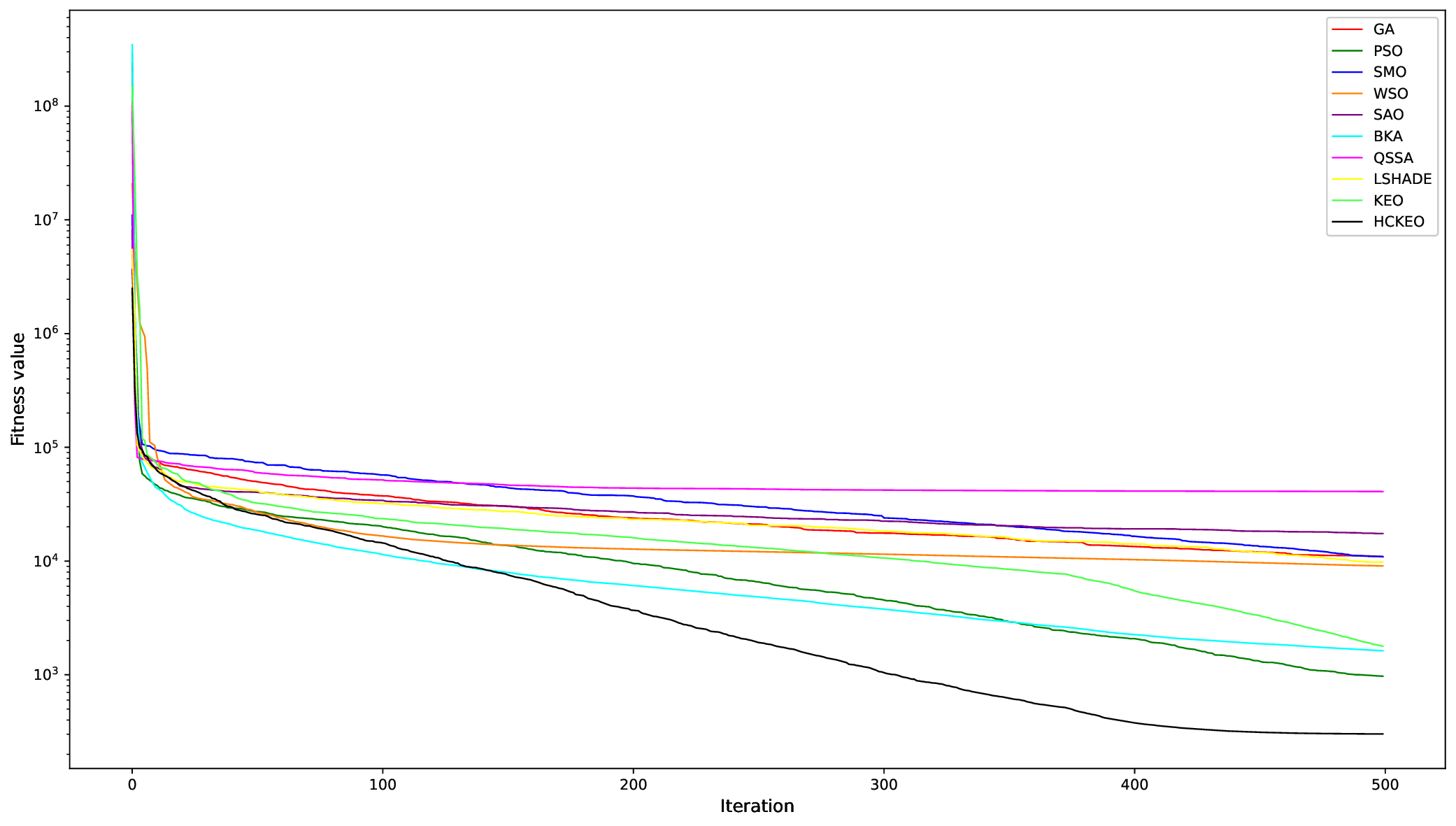}}
\subfloat[F2(D=20)]{\includegraphics[width=0.32\textwidth]{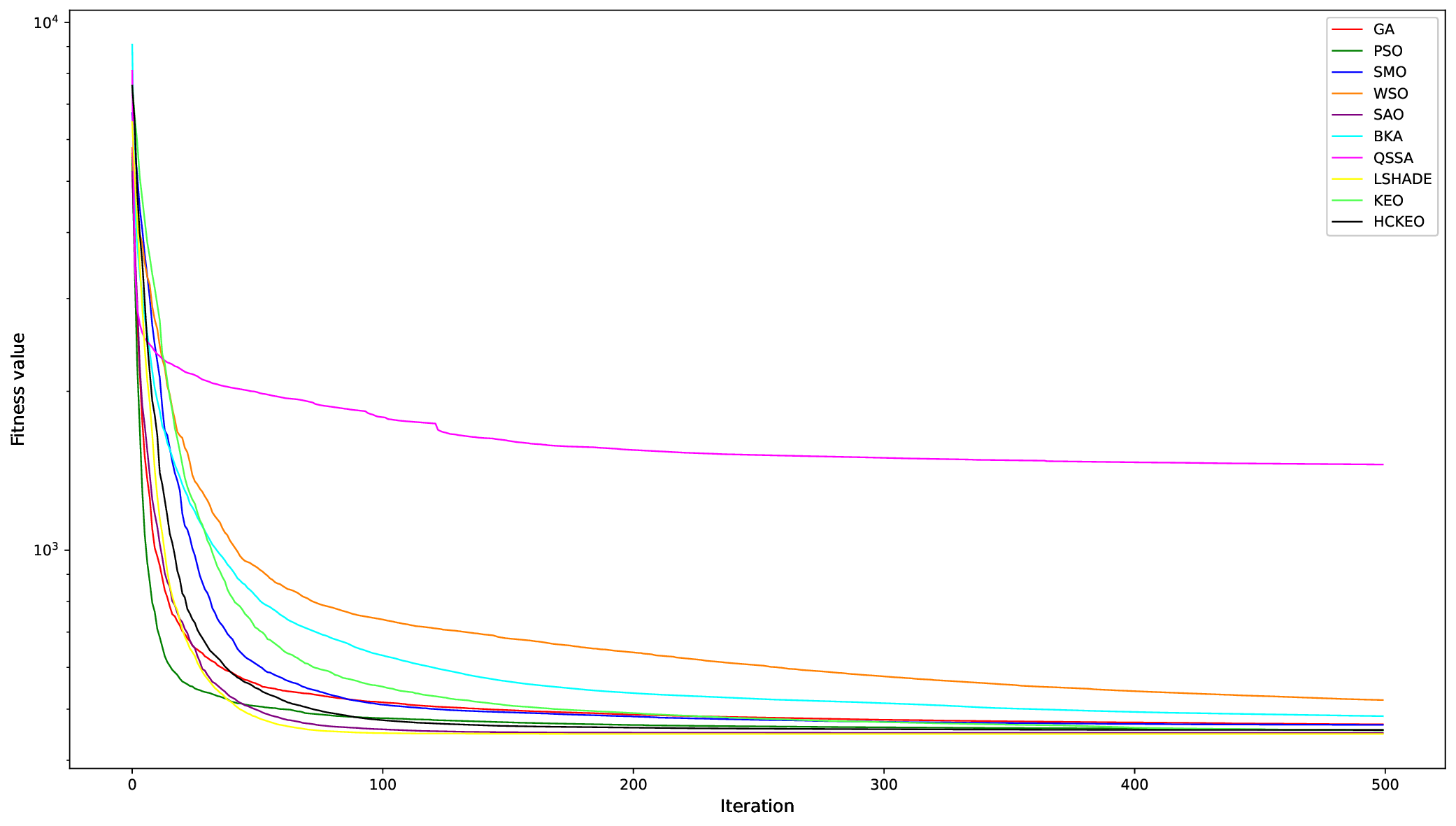}}
\subfloat[F3(D=20)]{\includegraphics[width=0.32\textwidth]{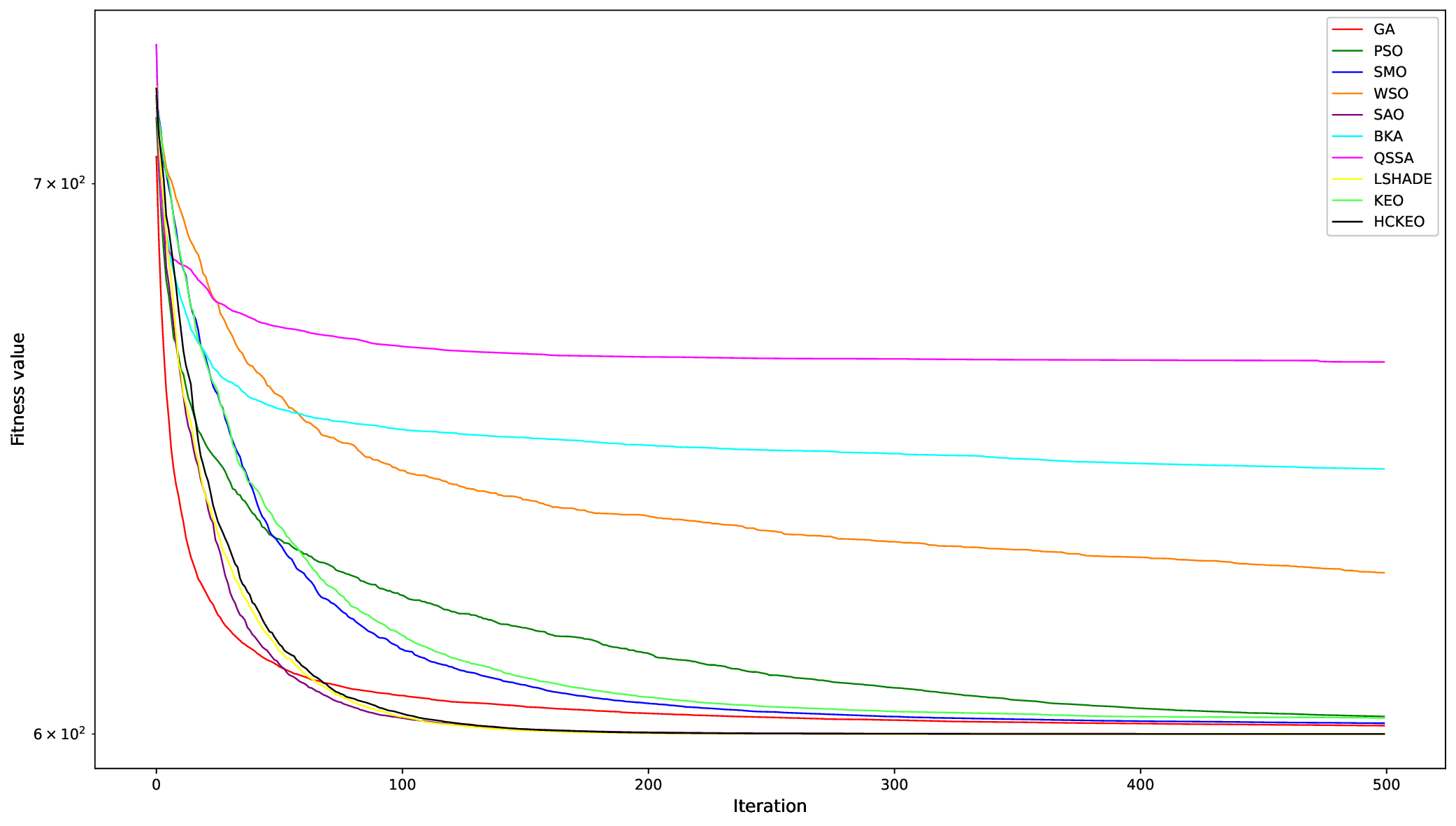}}\\
\subfloat[F3(D=20)]{\includegraphics[width=0.32\textwidth]{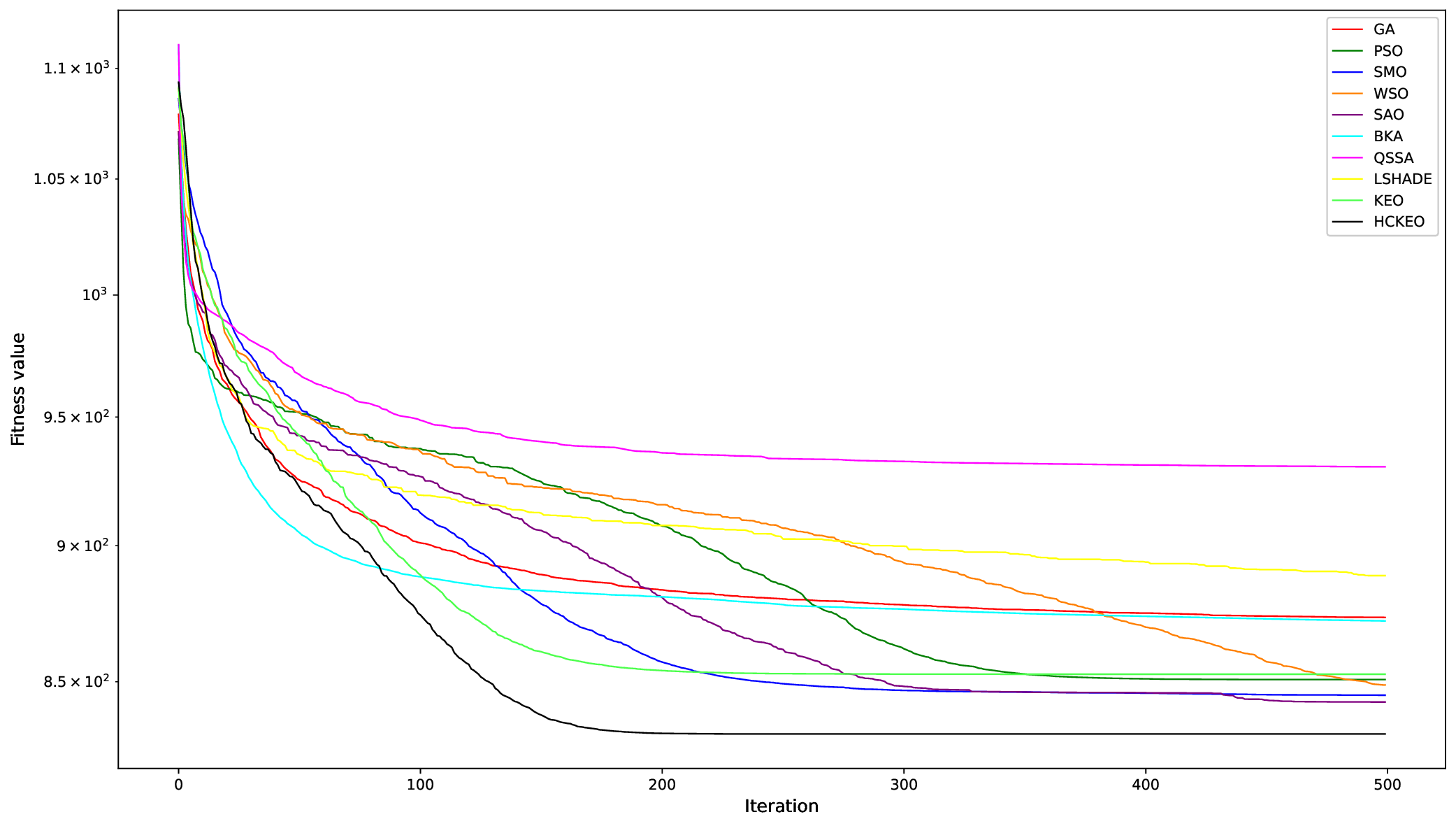}}
\subfloat[F4(D=20)]{\includegraphics[width=0.32\textwidth]{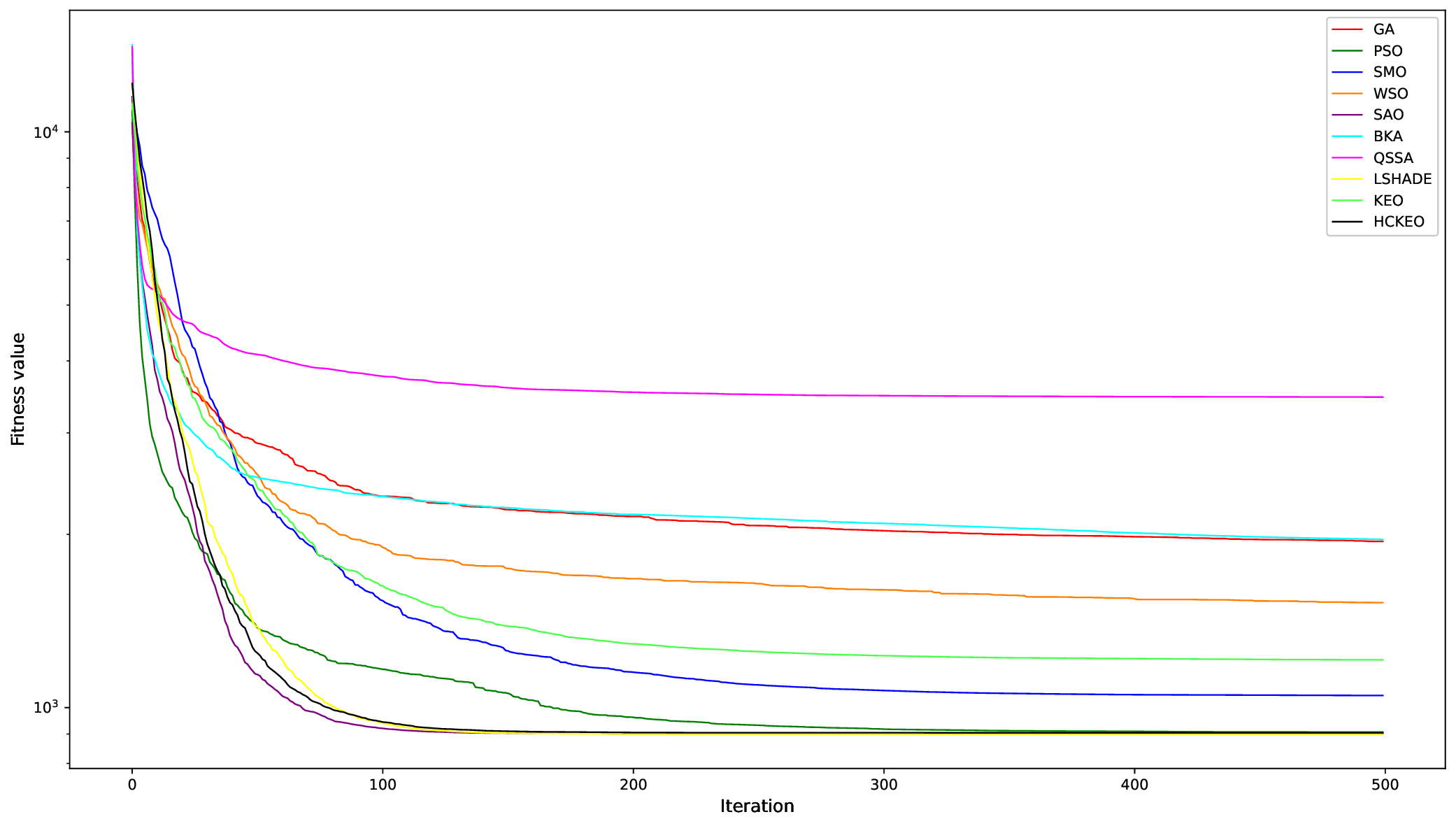}}
\subfloat[F6(D=20)]{\includegraphics[width=0.32\textwidth]{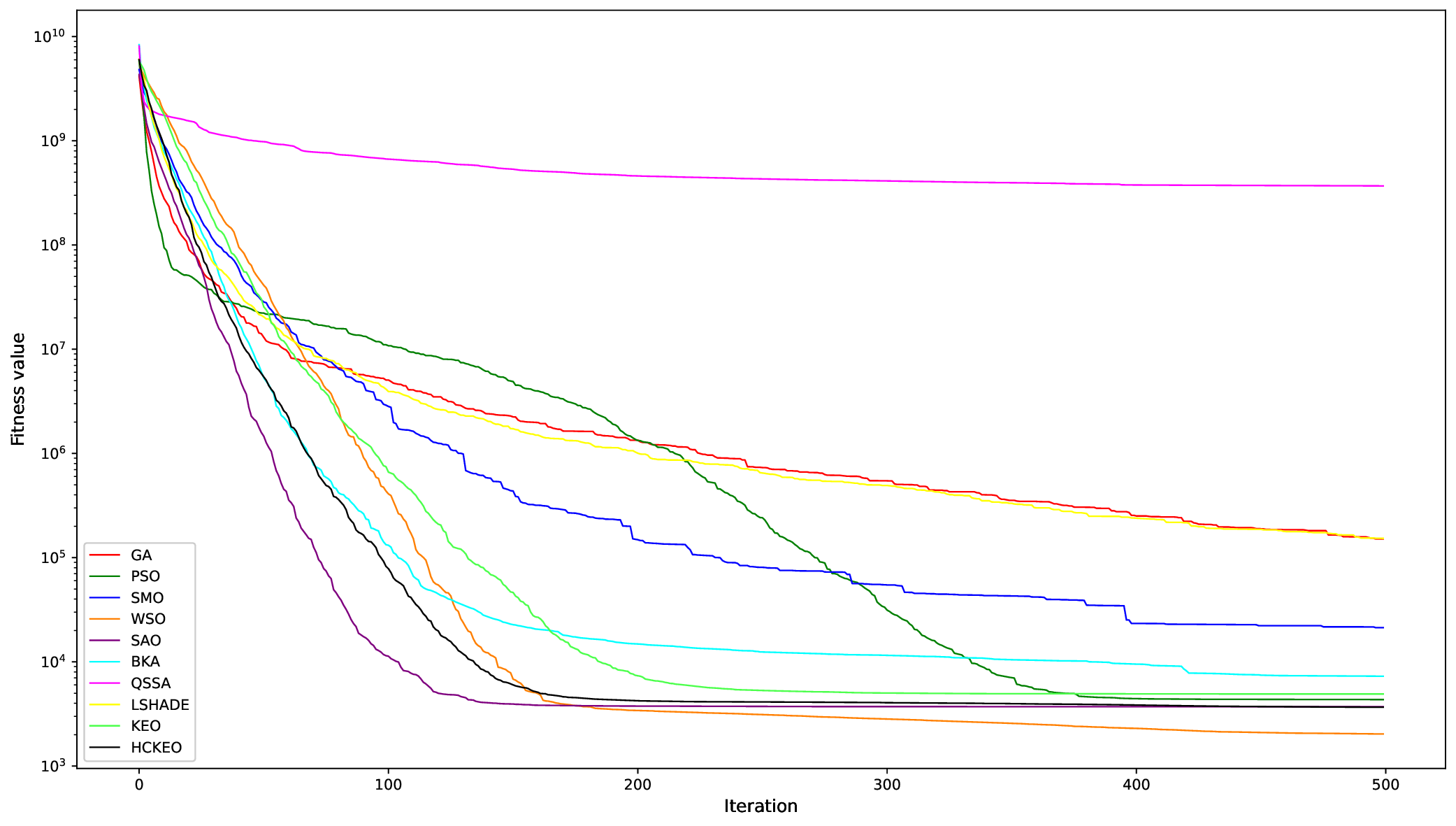}}\\
\subfloat[F7(D=20)]{\includegraphics[width=0.32\textwidth]{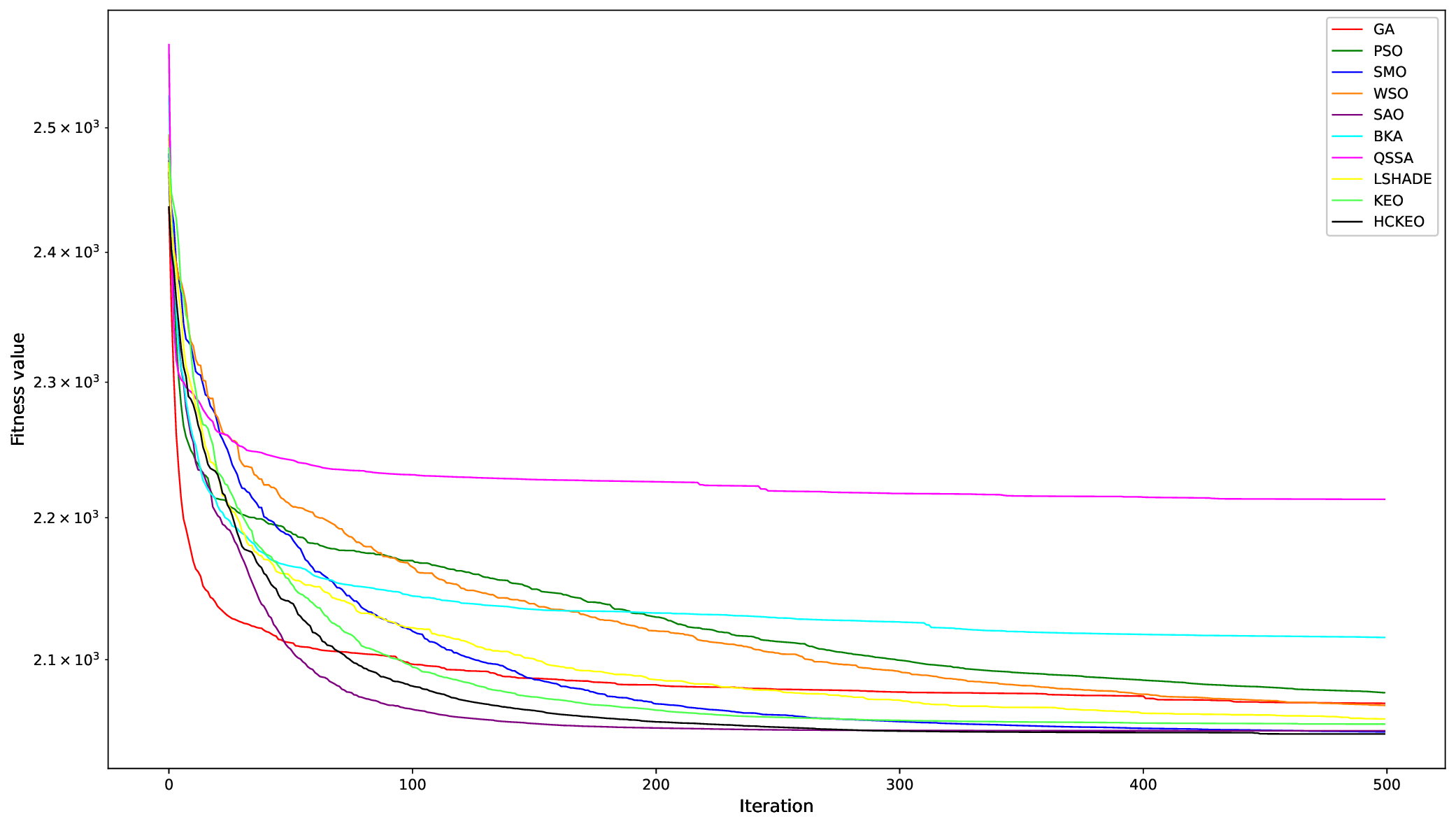}}
\subfloat[F8(D=20)]{\includegraphics[width=0.32\textwidth]{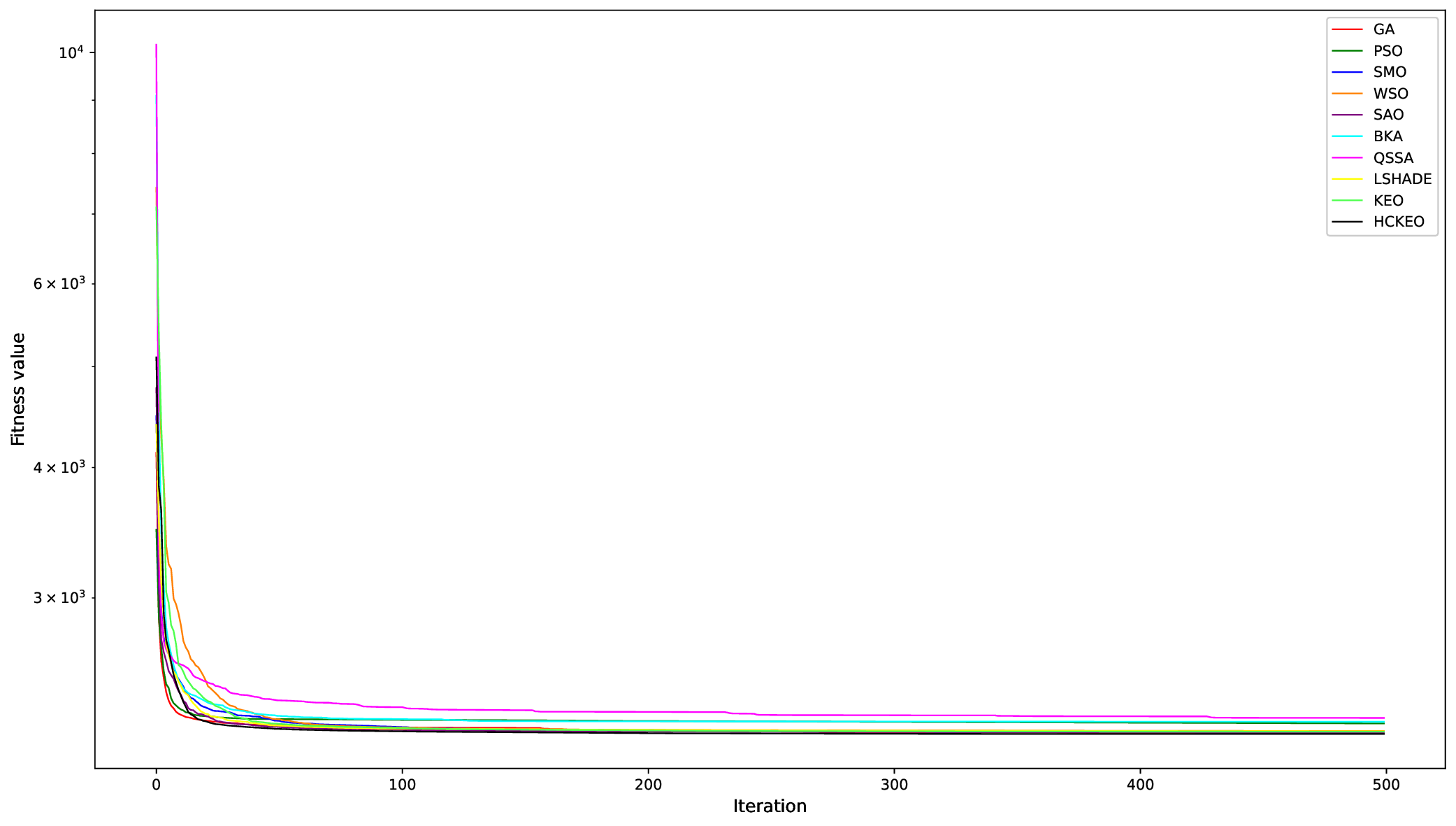}}
\subfloat[F9(D=20)]{\includegraphics[width=0.32\textwidth]{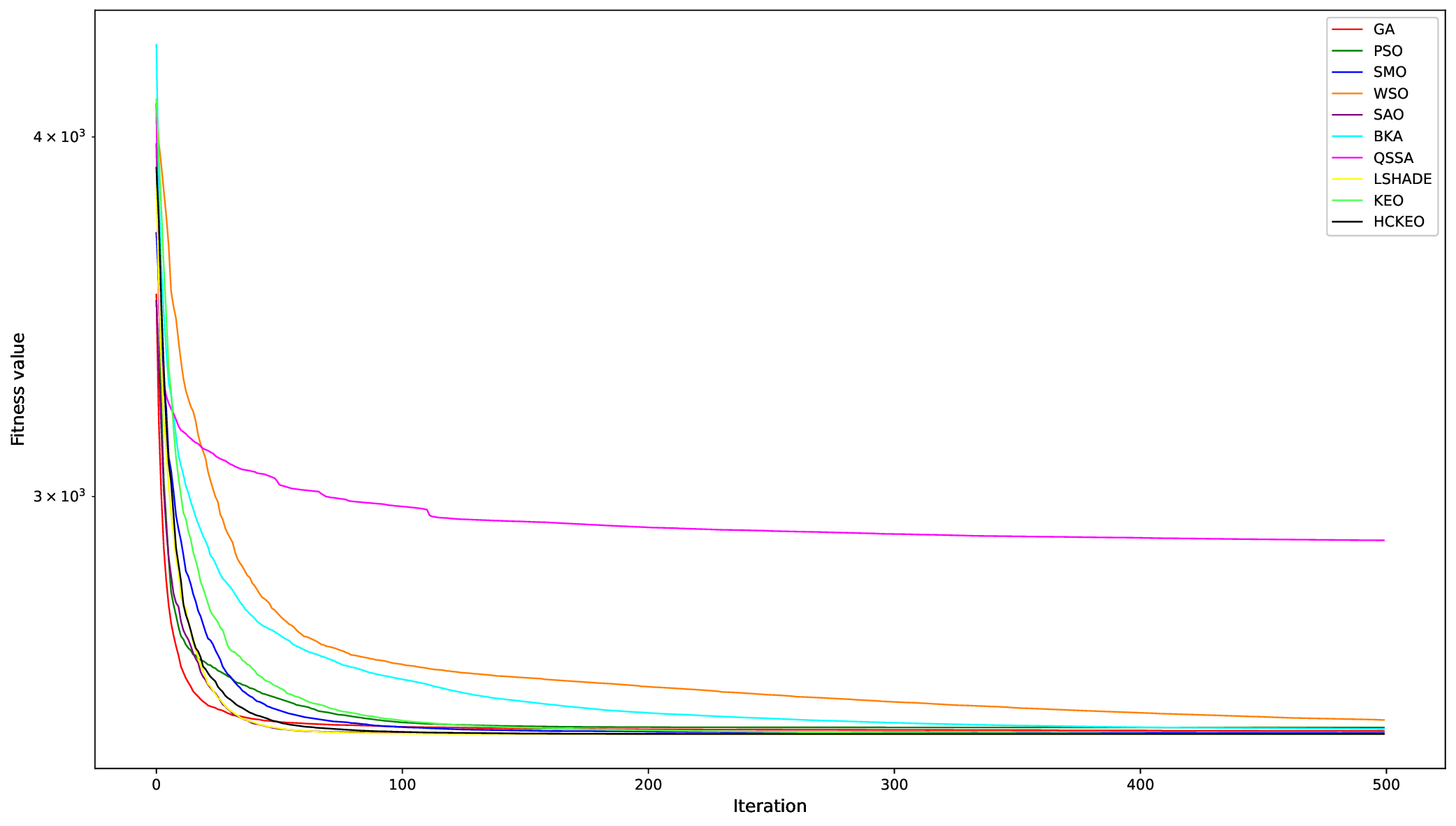}}\\
\subfloat[F10(D=20)]{\includegraphics[width=0.32\textwidth]{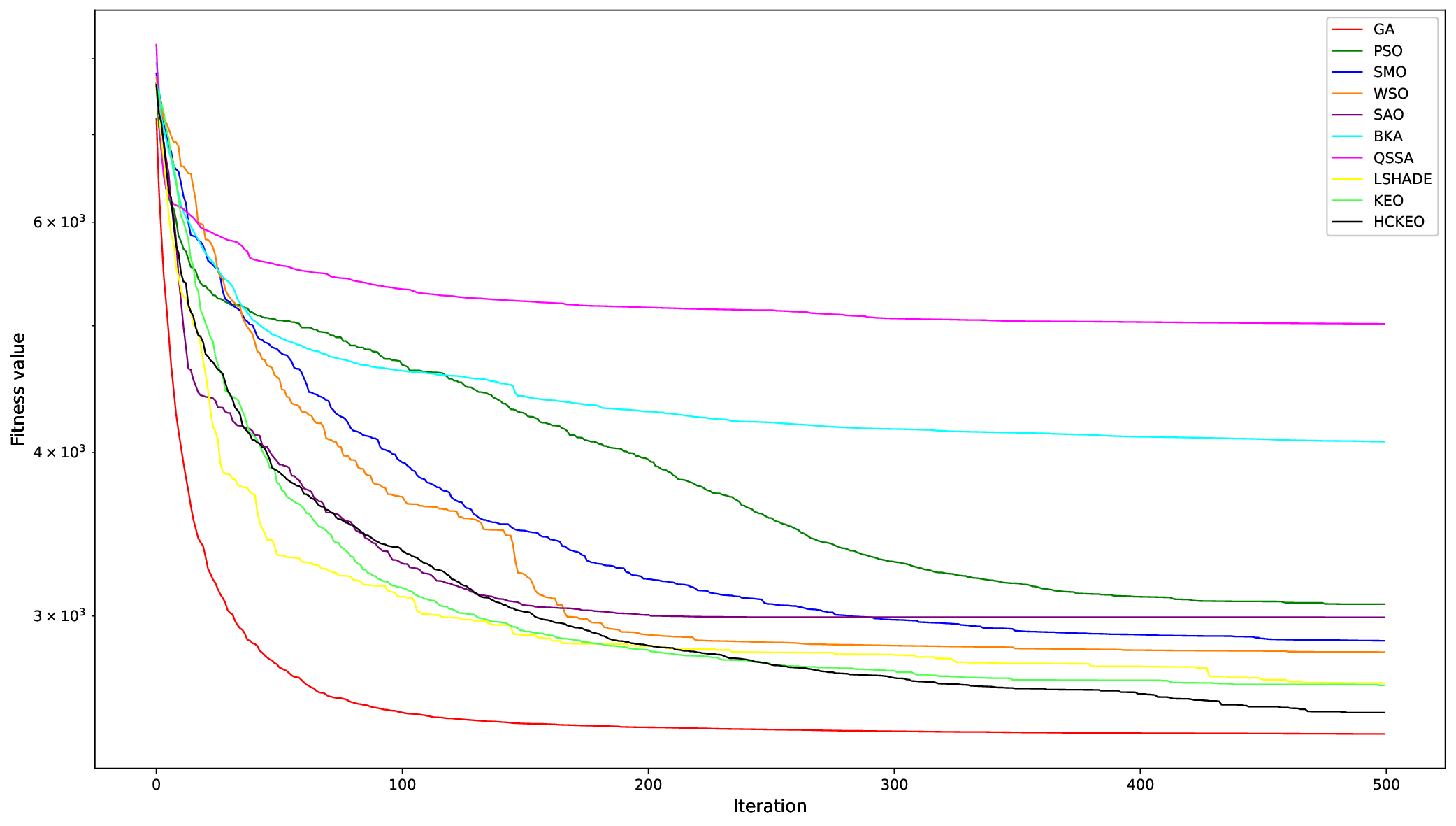}}
\subfloat[F11(D=20)]{\includegraphics[width=0.32\textwidth]{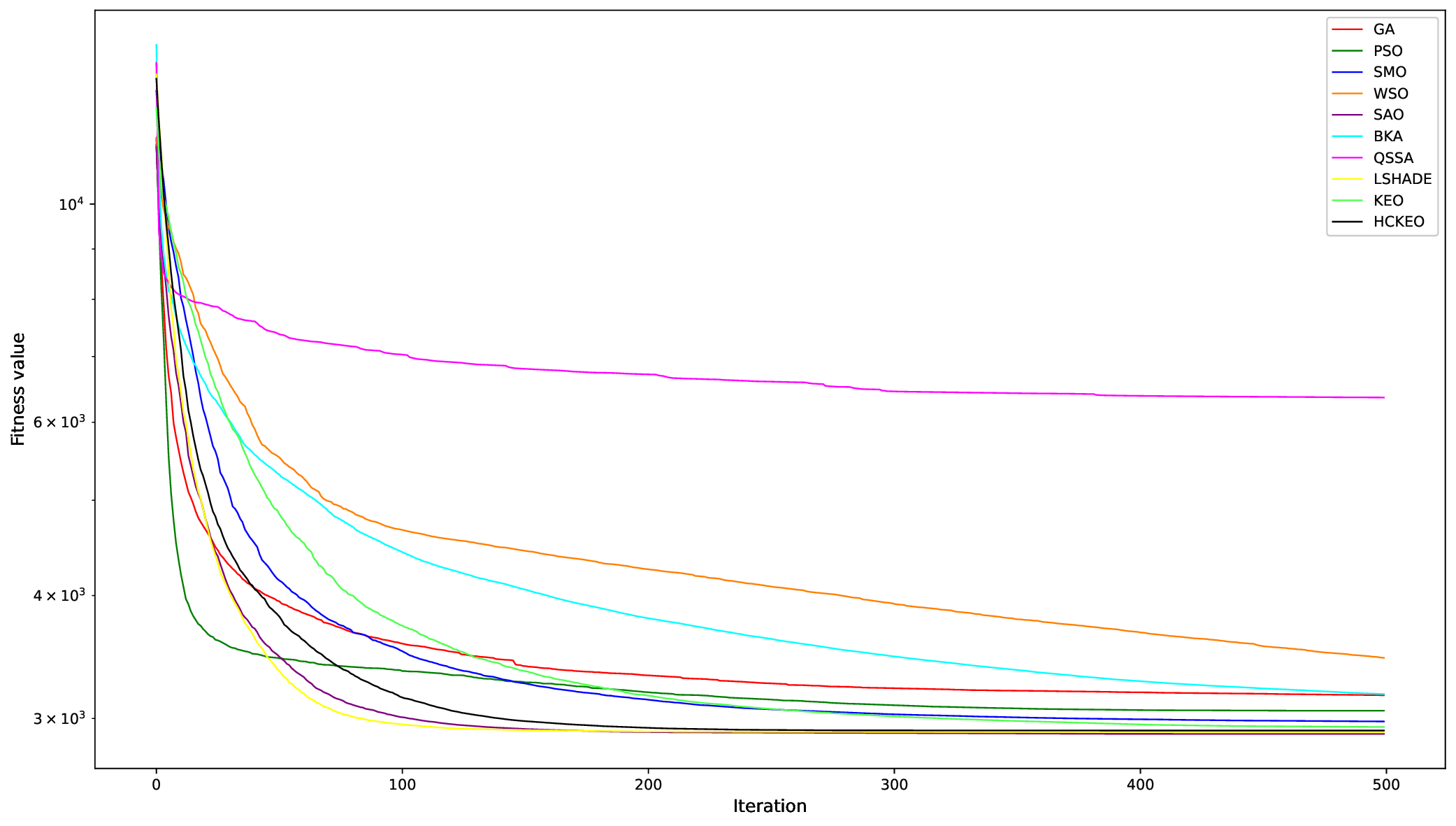}}
\subfloat[F12(D=20)]{\includegraphics[width=0.32\textwidth]{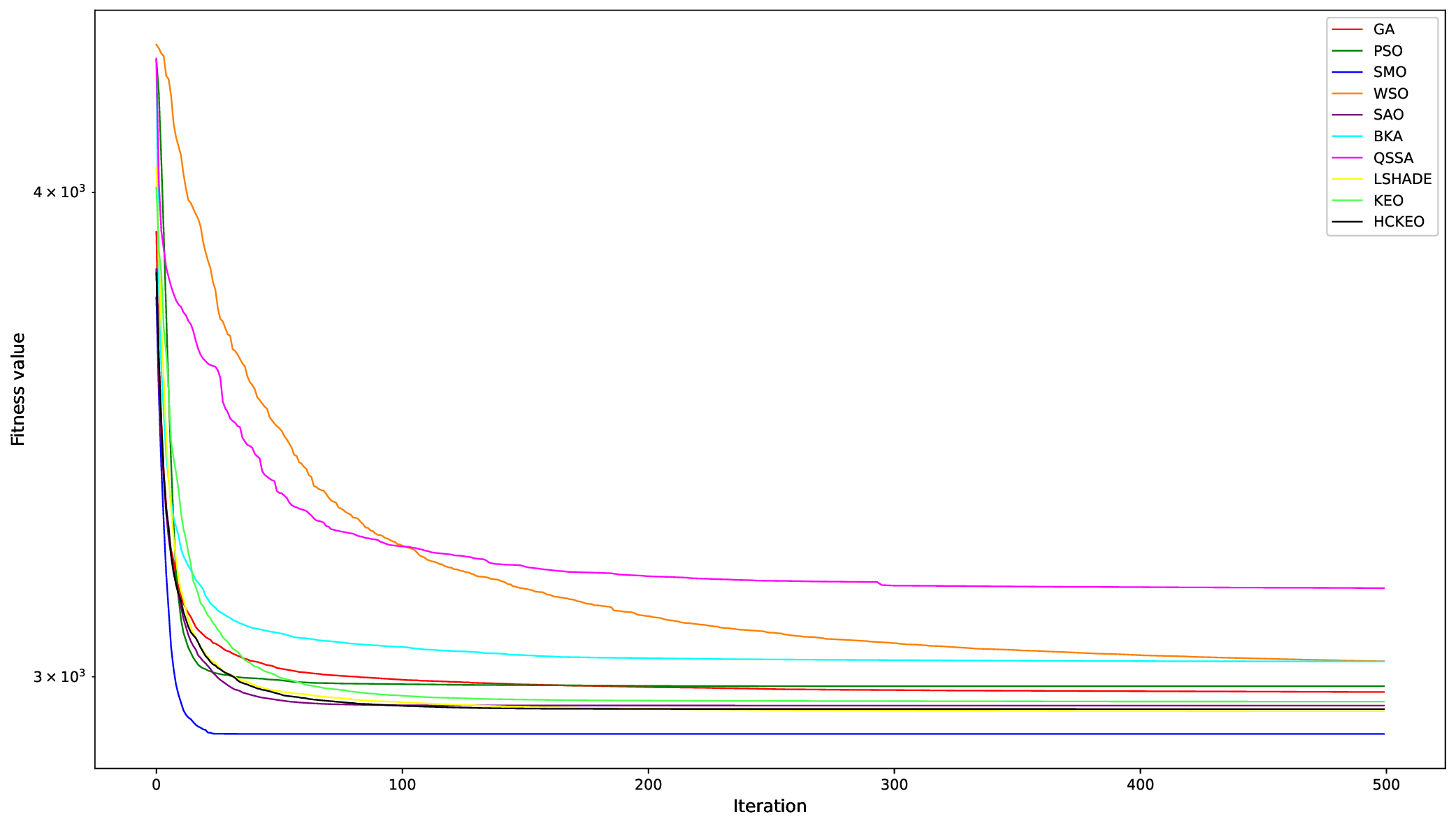}}\\
\caption{CEC2022 test functions convergence curve (Dim=20).}
\label{fig:curve20}
\end{figure}

Tables~\ref{tab:cec2022_10dim} and~\ref{tab:cec2022_20dim} present the comparative results of the proposed HCKEO algorithm and seven other algorithms on the CEC2022 benchmark functions with 10 and 20 dimensions, respectively.
The metrics include the mean, best, worst, and standard deviation (Std) of the final fitness values obtained after 30 independent runs.
In Table~\ref{tab:cec2022_20dim}, the overall performance trend remains consistent. HCKEO exhibits excellent robustness and global search capability even in higher-dimensional spaces.
In most function experiments (e.g., F1, F3, F4, F6, and F7), HCKEO's mean and optimal fitness values are significantly better than comparison algorithms.
Specifically, on the F1 function, HCKEO achieves a mean value of $301.74$, remarkably outperforming the other methods whose means exceed $900$.
In addition, the small standard deviations across almost all functions confirm the stability of the proposed algorithm under increased problem dimensionality.
These results verify that the hybrid strategy combining chaotic initialization and kernel updating significantly enhances the convergence precision and exploration–exploitation balance of the original KEO.
Across both 10D and 20D test suites, HCKEO demonstrates the best overall optimization performance, outperforming seven baseline algorithms in terms of accuracy, convergence stability, and scalability.
Compared with classical methods (GA, PSO), HCKEO converges faster and reaches better global minima.
Compared with recent metaheuristics (WSO, SAO, BKA), it achieves more reliable convergence with smaller fluctuations.
Even when compared with the advanced LSHADE algorithm, HCKEO maintains competitive or superior results on most functions, highlighting its potential as a robust optimizer for complex and high-dimensional search spaces.
Figure \ref{fig:curve10} and Figure \ref{fig:curve20} present the convergence curves from partial tests, indicating that the HCKEO algorithm exhibits superior convergence speed in most of the tests.

\subsubsection{Statistical test}

This work uses Wilcoxon rank sum test and Friedman test to further verify the statistical significance of the experiment and explain whether there is statistical significance between HCKEO and the comparative algorithm.

\paragraph{Wilcoxon rank sum test}

The Wilcoxon rank-sum test \cite{Phong2022} is conducted to perform pairwise comparisons between HCKEO and each competing algorithm.
This nonparametric test does not assume a specific distribution of the data and is widely adopted in optimization studies for evaluating algorithmic performance.

In this work, the significance level is set to $ \alpha = 0.05 $.
For each pairwise comparison, the null hypothesis $H_0$ assumes that there is no significant difference in the performance of the two algorithms statistically.
If the P-value is less than 0.05, it indicates that the performance difference between the two algorithms is statistically significant.

The resulting (p)-values obtained from the comparisons on the CEC2022 benchmark functions with dimensions (D=10) and (D=20) are reported in Tables \ref{tab:pvalue_cec2022_10}-\ref{tab:pvalue_cec2022_20}.

\begin{table}[!ht]
\centering
\caption{Wilcoxon rank-sum test $p$-values for pairwise comparisons on the CEC2022 benchmark functions (10D).}
\label{tab:pvalue_cec2022_10}
\resizebox{\textwidth}{!}{
\begin{tabular}{lccccccccc}
\hline
\textbf{Function} & \textbf{GA} & \textbf{PSO} & \textbf{SMO} & \textbf{WSO} & \textbf{SAO} & \textbf{BKA} & \textbf{QSSA} & \textbf{LSHADE} & \textbf{KEO} \\ \hline
F1  & 1.6040e-11 & 1.6040e-11 & 1.6040e-11 & 1.6040e-11 & 1.7746e-11 & 1.6040e-11 & 1.6040e-11 & 1.6040e-11 & 1.6040e-11 \\
F2  & 4.9726e-03 & \textbf{8.2007e-01} & \textbf{4.0824e-01} & 2.0701e-07 & \textbf{2.1134e-01} & \textbf{5.4619e-01} & 3.0623e-10 & \textbf{5.5021e-02} & 3.9010e-02 \\
F3  & 1.4498e-11 & 2.1925e-10 & 1.6342e-10 & 1.9762e-11 & \textbf{8.0625e-02} & 1.4498e-11 & 1.4498e-11 & 4.8728e-09 & 9.0261e-11 \\
F4  & 9.7037e-10 & \textbf{3.4019e-01} & 2.8333e-04 & \textbf{2.5181e-01} & \textbf{4.5911e-01} & 5.9565e-05 & 3.4547e-10 & 1.0972e-08 & 4.8951e-05 \\
F5  & 1.6040e-11 & 6.2937e-09 & 3.6637e-09 & 1.1042e-09 & \textbf{6.4763e-02} & 1.6040e-11 & 1.6040e-11 & 7.6872e-12 & 2.5431e-09 \\
F6  & 3.0199e-11 & 3.3384e-11 & 6.6955e-11 & 3.0939e-06 & 1.6132e-10 & 3.0199e-11 & 3.0196e-11 & 8.3520e-08 & 1.4110e-09 \\
F7  & \textbf{4.4642e-01} & \textbf{1.0233e-01} & 2.8789e-06 & 6.9125e-04 & 1.8731e-07 & 3.0199e-11 & 3.0199e-11 & 1.0035e-03 & 1.4412e-02 \\
F8  & 7.6973e-04 & 5.2640e-04 & 1.8681e-05 & 4.0839e-05 & 2.2360e-02 & 5.4617e-09 & 3.0199e-11 & 7.6171e-03 & 1.2362e-03 \\
F9  & 1.2118e-12 & 4.1926e-02 & 3.3592e-11 & 1.2118e-12 & \textbf{1} & 1.2118e-12 & 1.2117e-12 & \textbf{1} & 5.5419e-03 \\
F10 & 3.9881e-04 & 1.8368e-02 & 2.3884e-04 & 4.7138e-04 & 1.4423e-03 & 7.3803e-10 & 1.3288e-10 & 3.0317e-02 & 3.5201e-07 \\
F11 & 3.4737e-10 & 5.7470e-09 & 1.0931e-09 & 2.8125e-08 & 2.0471e-01 & 9.8351e-09 & 1.2590e-11 & \textbf{8.3657e-01} & 1.2822e-08 \\
F12 & 1.7162e-10 & 3.7407e-09 & 5.7964e-09 & 2.5222e-10 & 4.5259e-07 & 1.4621e-06 & 2.9082e-11 & 3.3897e-03 & 5.3323e-06 \\ \hline
\end{tabular}}
\end{table}

\begin{table}[!ht]
\centering
\caption{Wilcoxon rank-sum test $p$-values for pairwise comparisons on the CEC2022 benchmark functions (20D).}
\label{tab:pvalue_cec2022_20}
\resizebox{\textwidth}{!}{
\begin{tabular}{lccccccccc}
\hline
\textbf{Function} & \textbf{GA} & \textbf{PSO} & \textbf{SMO} & \textbf{WSO} & \textbf{SAO} & \textbf{BKA} & \textbf{QSSA} & \textbf{LSHADE} & \textbf{KEO} \\ \hline
F1  & 3.0199e-11 & 3.0199e-11 & 3.0199e-11 & 3.0198e-11 & 3.0199e-11 & 3.0199e-11 & 3.0199e-11 & 3.0199e-11 & 3.0199e-11 \\
F2  & 4.1127e-07 & \textbf{2.0095e-01} & \textbf{1.7612e-01} & 8.1014e-10 & \textbf{8.7710e-02} & 7.6950e-08 & 3.0198e-11 & 5.0912e-06 & 8.1200e-04 \\
F3  & 3.0199e-11 & 6.6955e-11 & 3.0199e-11 & 3.0199e-11 & 8.1465e-05 & 3.0199e-11 & 3.0198e-11 & 3.0199e-11 & 3.0199e-11 \\
F4  & 4.1784e-09 & 1.1929e-06 & 4.2159e-04 & 7.6947e-04 & \textbf{7.6182e-01} & 3.4703e-10 & 3.0161e-11 & 3.0161e-11  & 1.3359e-05 \\
F5  & 3.0199e-11 & \textbf{1.2597e-01} & 2.4386e-09 & 5.4941e-11 & 2.6784e-06 & 3.0199e-11 & 3.0198e-11 & 3.0066e-11 & 1.9568e-10 \\
F6  & 3.0199e-11 & \textbf{8.0727e-01} & 1.1710e-02 & 4.9818e-04 & 7.1719e-01 & \textbf{5.1877e-02} & 4.1997e-10  & 3.0199e-11 & 8.3026e-01 \\
F7  & 1.9963e-05 & 1.3249e-04 & 1.2212e-02 & 1.1747e-04 & \textbf{6.7869e-02} & 3.8202e-10 & 3.0198e-11 & 2.4913e-06 & 1.4932e-04 \\
F8  & 6.0658e-11 & 3.1589e-10 & 3.1588e-10 & 3.0199e-11 & \textbf{1.1882e-01} & 3.0199e-11 & 3.0198e-11 & 3.0199e-11 & 6.5486e-04 \\
F9  & 3.0199e-11 & 3.0199e-11 & \textbf{1.8576e-01} & 3.0199e-11 & 4.5043e-11 & 3.0199e-11 & 3.0198e-11 & 3.0123e-11 & 3.0199e-11 \\
F10 & 1.5014e-02 & 3.9881e-04 & 9.7917e-05 & 4.3531e-05 & 5.0842e-03 & 8.4848e-09 & 2.2273e-09 & 1.2597e-01 & 2.4994e-03 \\
F11 & 4.2175e-04 & 3.3681e-05 & 1.9526e-02 & 2.8314e-08 & 3.2651e-02 & 1.6980e-08 & 3.0198e-11 & 5.6073e-05 & 7.2884e-03 \\
F12 & 8.9934e-11 & 3.1589e-10 & 3.0199e-11 & 3.0199e-11 & 3.3679e-04 & 6.0658e-11 & 3.0198e-11 & 4.8413e-02 & 3.2555e-07 \\ \hline
\end{tabular}}
\end{table}

As shown in Tables~\ref{tab:pvalue_cec2022_10} and~\ref{tab:pvalue_cec2022_20}, most of the p-values between HCKEO and other algorithms are smaller than 0.05, indicating significant differences in optimization performance.
Only a few functions, such as F2, F4, F6, and F8, exhibit p-values greater than 0.05 when compared with PSO, SAO, or LSHADE, suggesting comparable performance in these specific cases.
Overall, the statistical evidence confirms that HCKEO outperforms the comparison algorithms in most benchmark functions across both 10- and 20-dimensional test sets, validating the robustness and superiority of the proposed improvements.

\paragraph{Freideman’s test}

To further evaluate the performance of HCKEO, the non-parametric Friedman average rank test is adopted to compare the results obtained by HCKEO and other algorithms on the CEC2022 benchmark suite.
This statistical test does not require assumptions about the underlying data distribution, which makes it well suited for comparing optimization algorithms over multiple benchmark functions.

\begin{table}[!ht]
\centering
\caption{Friedman ranking results on the CEC2022 benchmark functions.}
\label{tab:cec2022_friedman}
\begin{tabular}{lcc}
\toprule
\textbf{Algorithm} & \textbf{10D (Rank)} & \textbf{20D (Rank)} \\
\midrule
GA     & 7.0000 (8)  & 6.0000 (6) \\
PSO    & 6.4167 (7)  & 6.0000 (6) \\
SMO    & 4.7500 (3)  & 4.6667 (5) \\
WSO    & 5.5000 (6)  & 6.5833 (8) \\
SAO    & 4.7500 (3)  & 3.5000 (2) \\
BKA    & 7.8833 (9)  & 7.6667 (9) \\
QSSA   & 9.9167 (10) & 10.0000 (10) \\
LSHADE & 3.0000 (2)  & 3.9167 (3) \\
KEO    & 4.7500 (3)  & 4.5833 (4) \\
HCKEO  & \textbf{1.8333 (1)} & \textbf{2.0833 (1)} \\
\bottomrule
\end{tabular}
\end{table}

According to the Friedman test results in Table~\ref{tab:cec2022_friedman}, HCKEO achieves the lowest mean rank among all compared algorithms in both 10- and 20-dimensional problems, ranking first overall.
In contrast, conventional algorithms such as GA, PSO, and BKA obtain comparatively larger ranking values, suggesting weaker global search capability and slower convergence behavior.
This demonstrates that HCKEO consistently performs better than other algorithms across the majority of test functions.
The high rank values of GA, PSO, and BKA indicate that their global search ability is weaker compared to other algorithms.
Meanwhile, advanced optimizers like LSHADE and SAO achieve moderate ranks but still fall behind HCKEO in overall performance.
Therefore, the Friedman analysis provides further statistical evidence supporting that HCKEO maintains superior optimization accuracy and stability across different dimensional settings.

\subsubsection{Analysis of population diversity}

\begin{figure}[!ht]
\centering
\subfloat[F1]{\includegraphics[width=1.8in,height=1.8in]{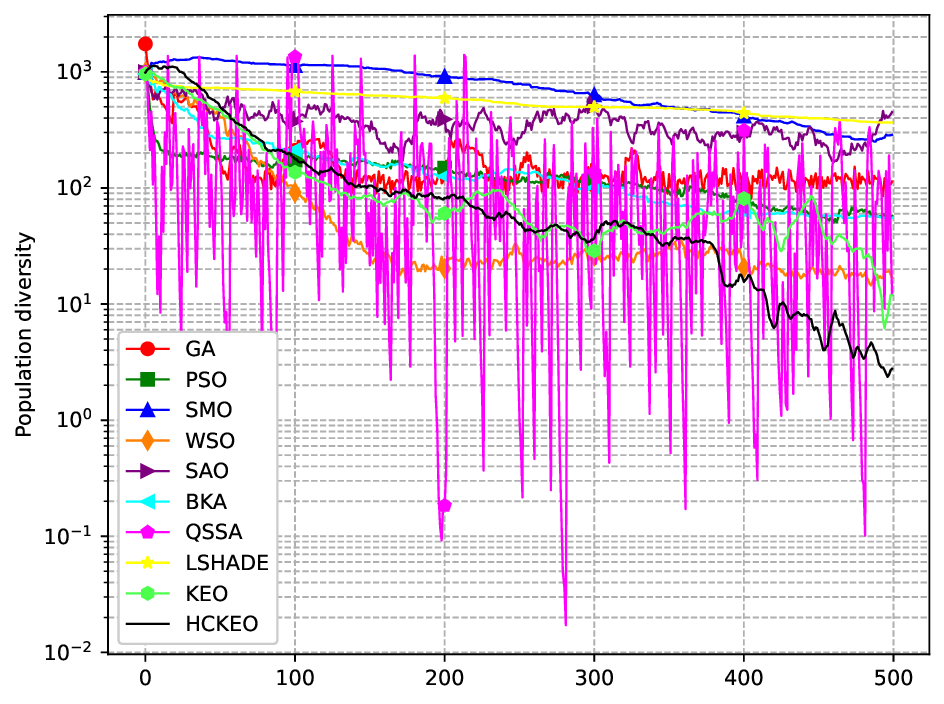}}
\subfloat[F2]{\includegraphics[width=1.8in,height=1.8in]{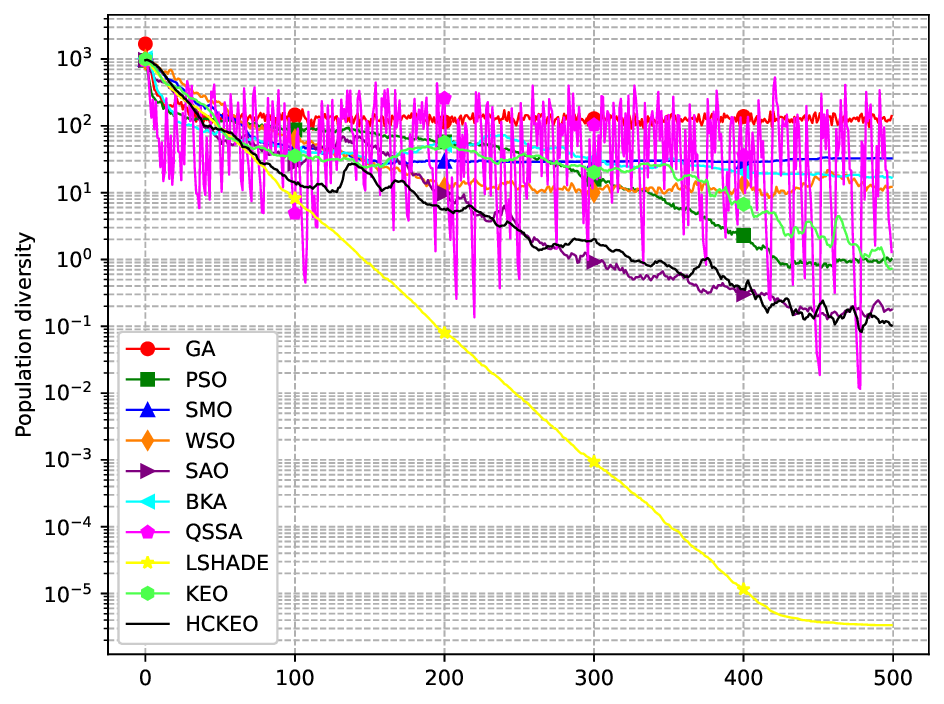}}
\subfloat[F3]{\includegraphics[width=1.8in,height=1.8in]{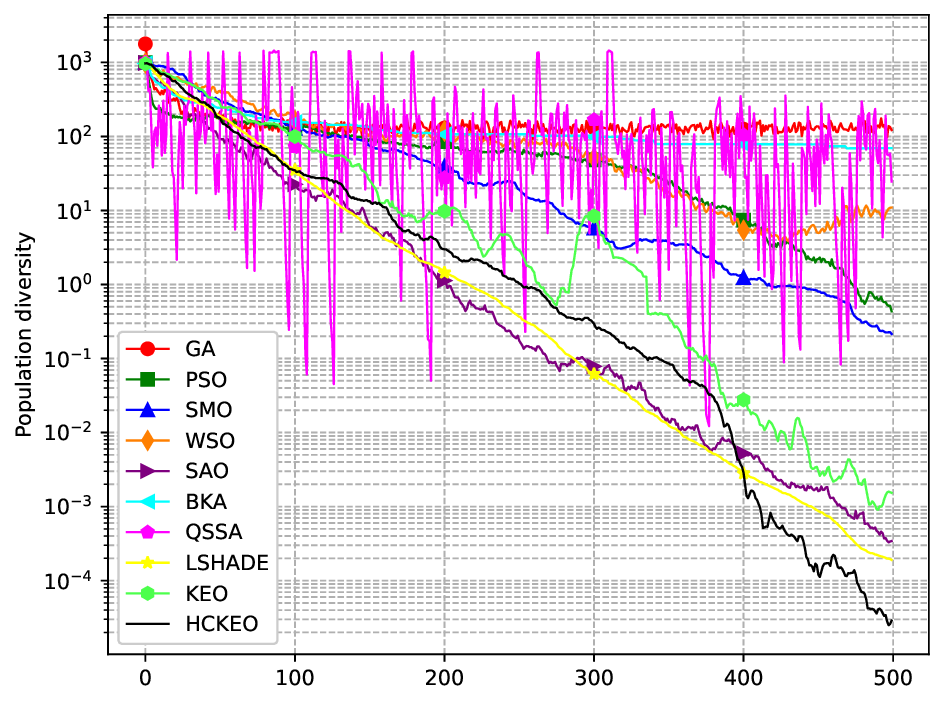}}\\
\subfloat[F4]{\includegraphics[width=1.8in,height=1.8in]{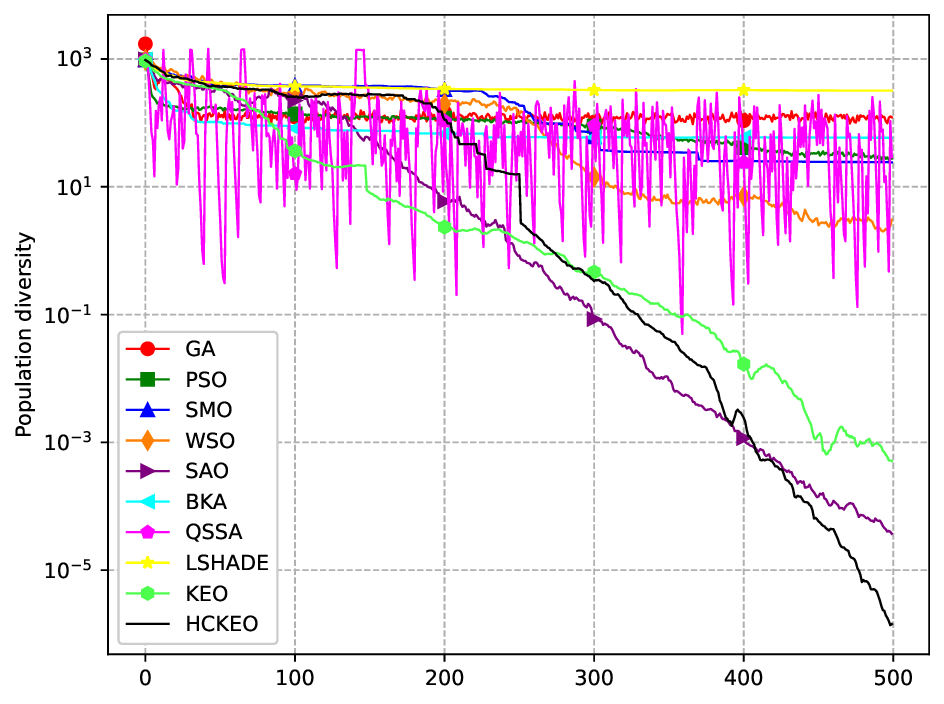}}
\subfloat[F5]{\includegraphics[width=1.8in,height=1.8in]{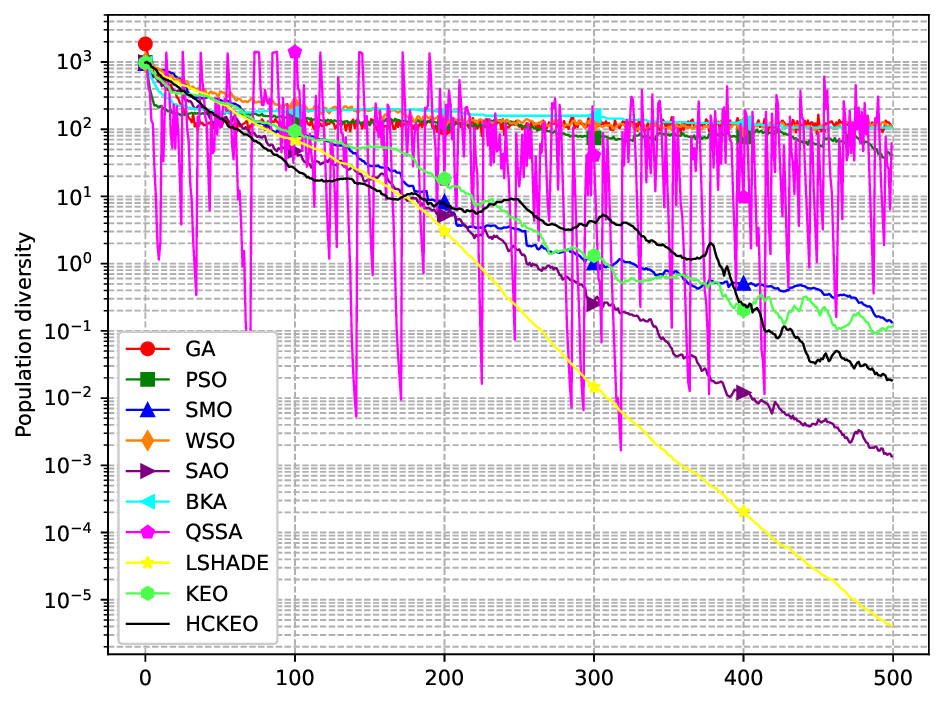}}
\subfloat[F6]{\includegraphics[width=1.8in,height=1.8in]{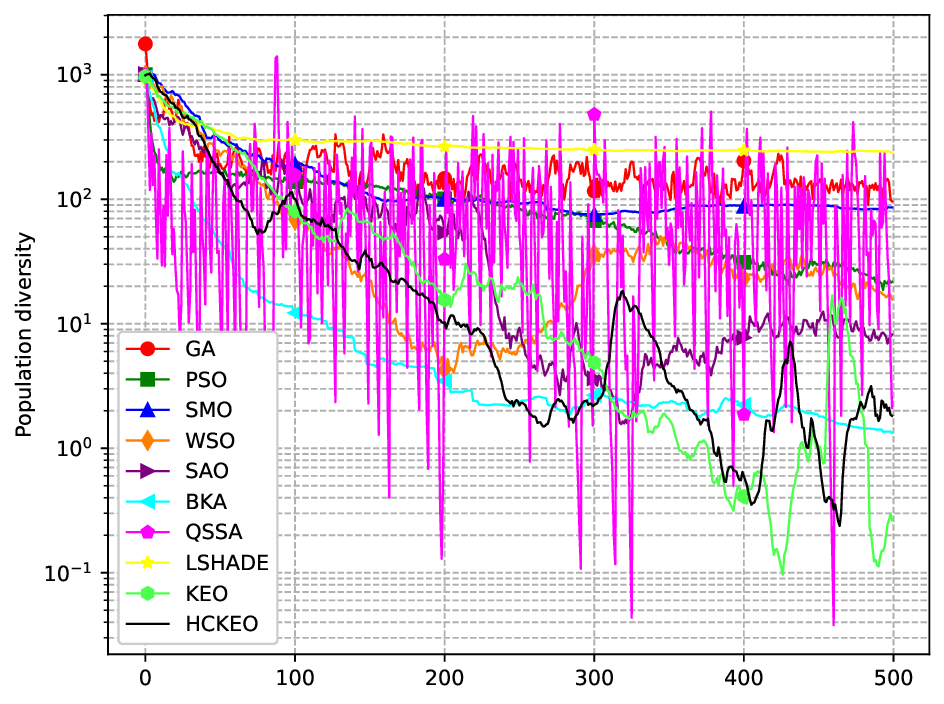}}\\
\subfloat[F7]{\includegraphics[width=1.8in,height=1.8in]{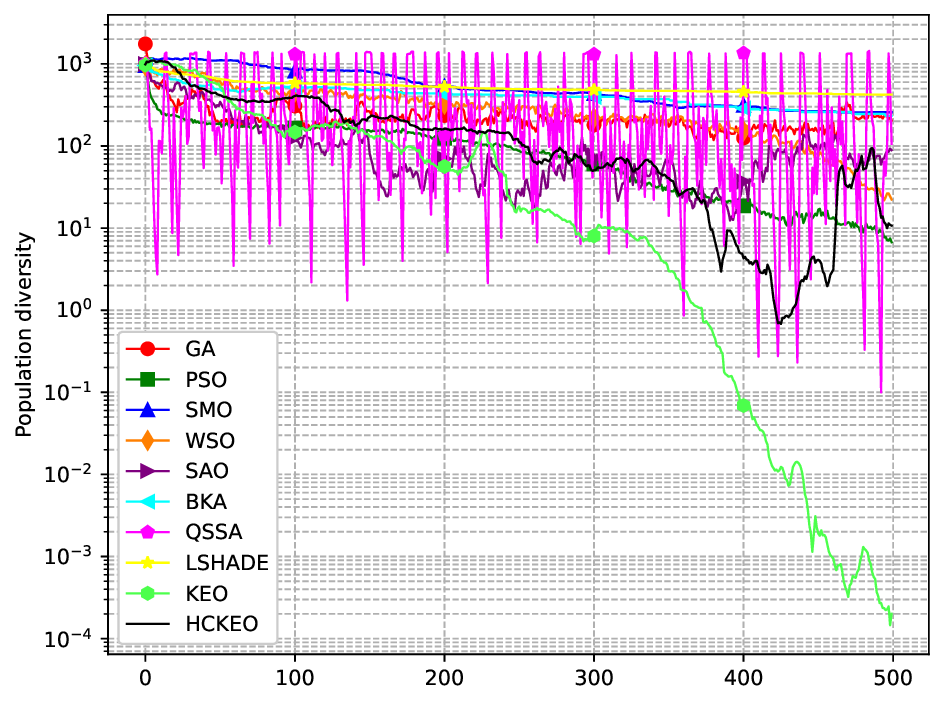}}
\subfloat[F8]{\includegraphics[width=1.8in,height=1.8in]{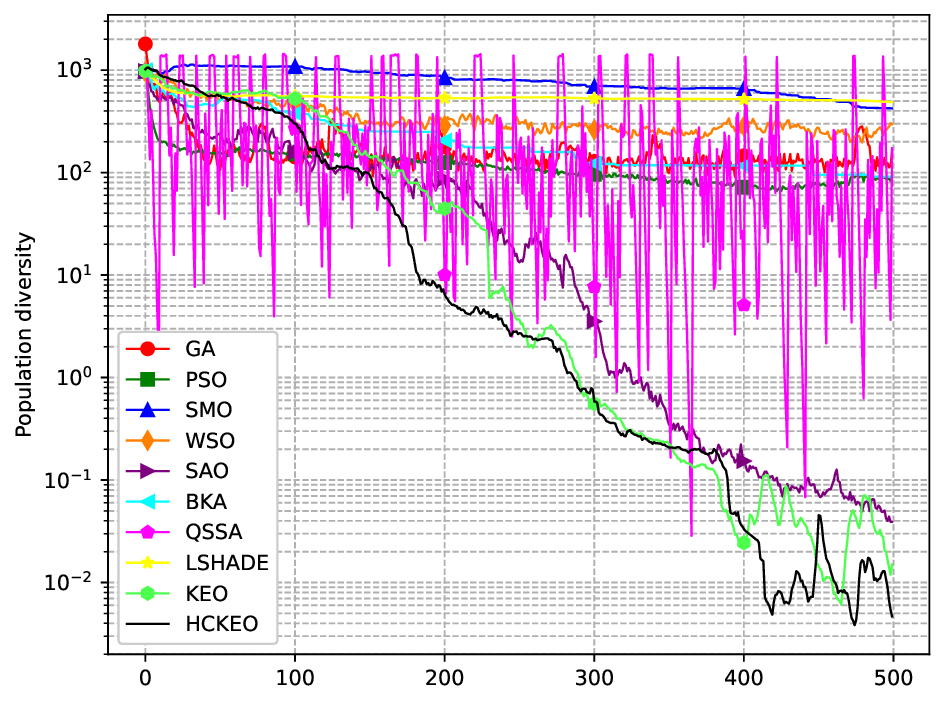}}
\subfloat[F9]{\includegraphics[width=1.8in,height=1.8in]{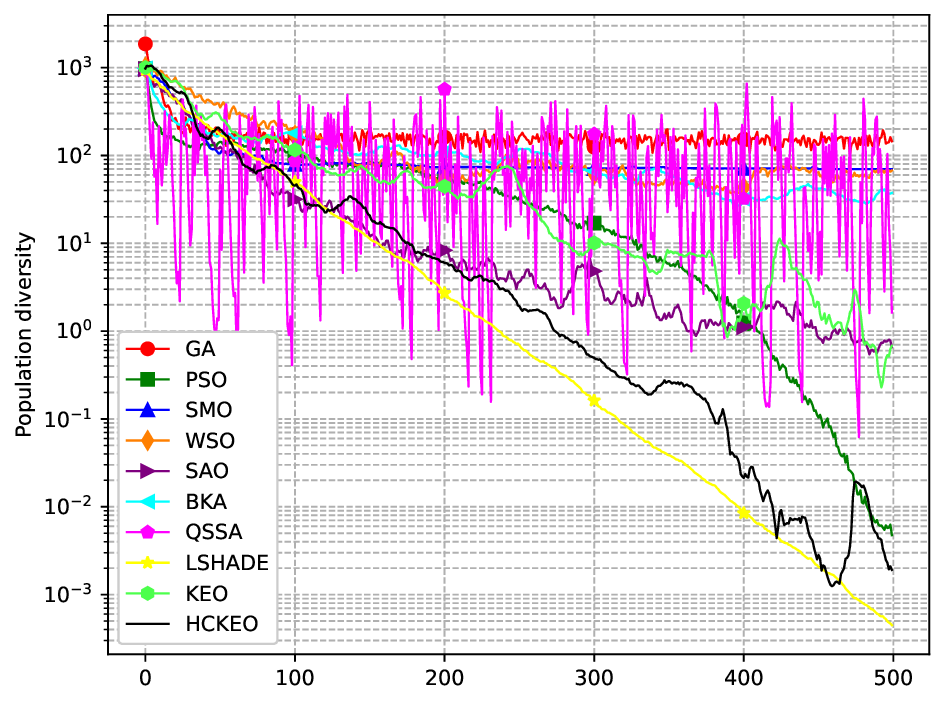}}\\
\subfloat[F10]{\includegraphics[width=1.8in,height=1.8in]{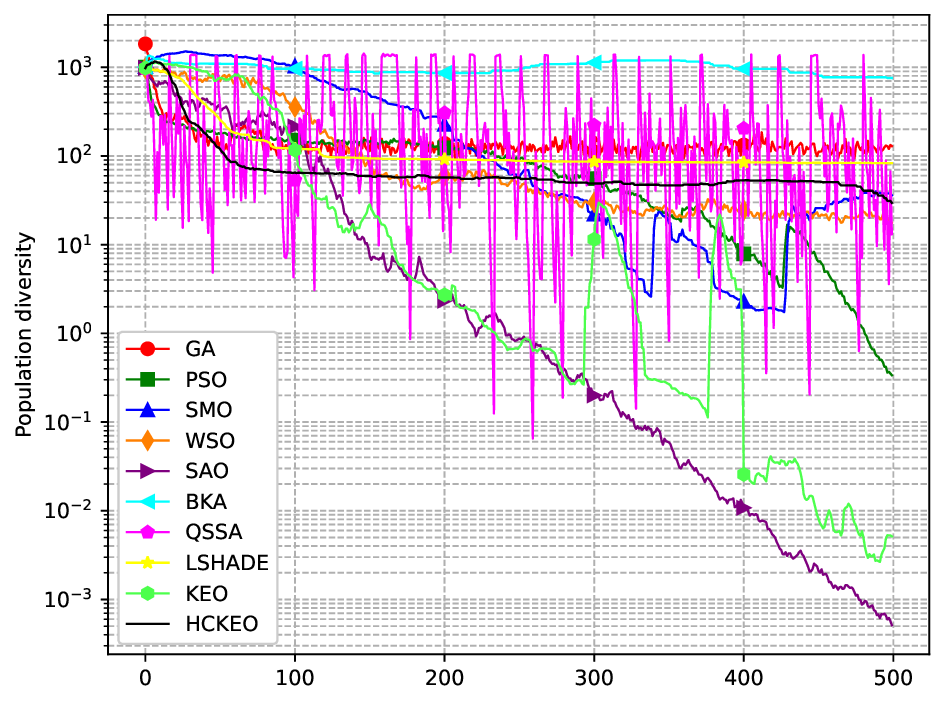}}
\subfloat[F11]{\includegraphics[width=1.8in,height=1.8in]{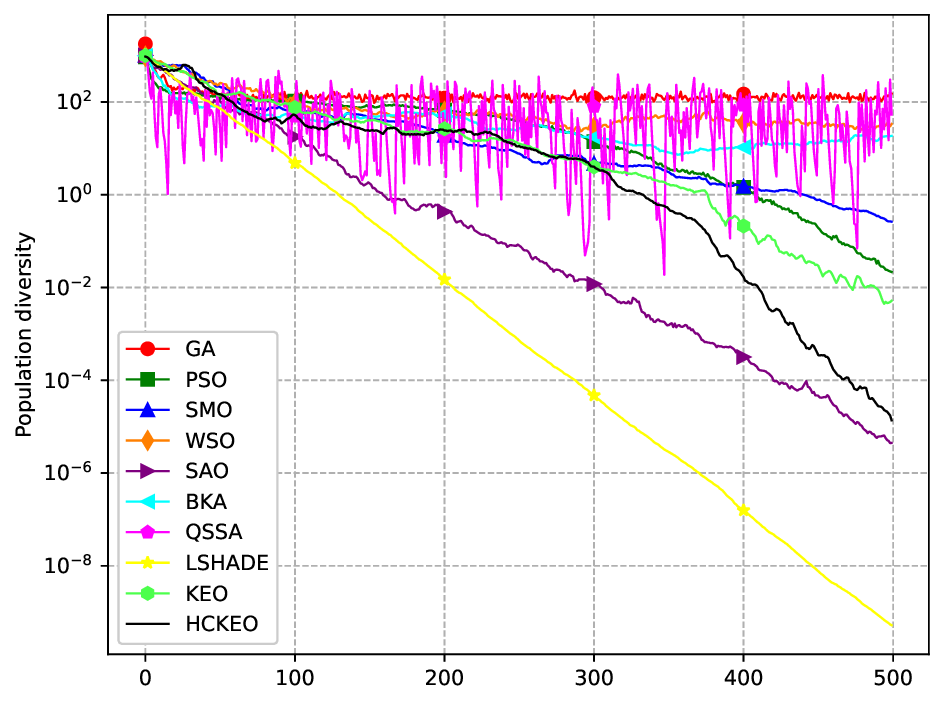}}
\subfloat[F12]{\includegraphics[width=1.8in,height=1.8in]{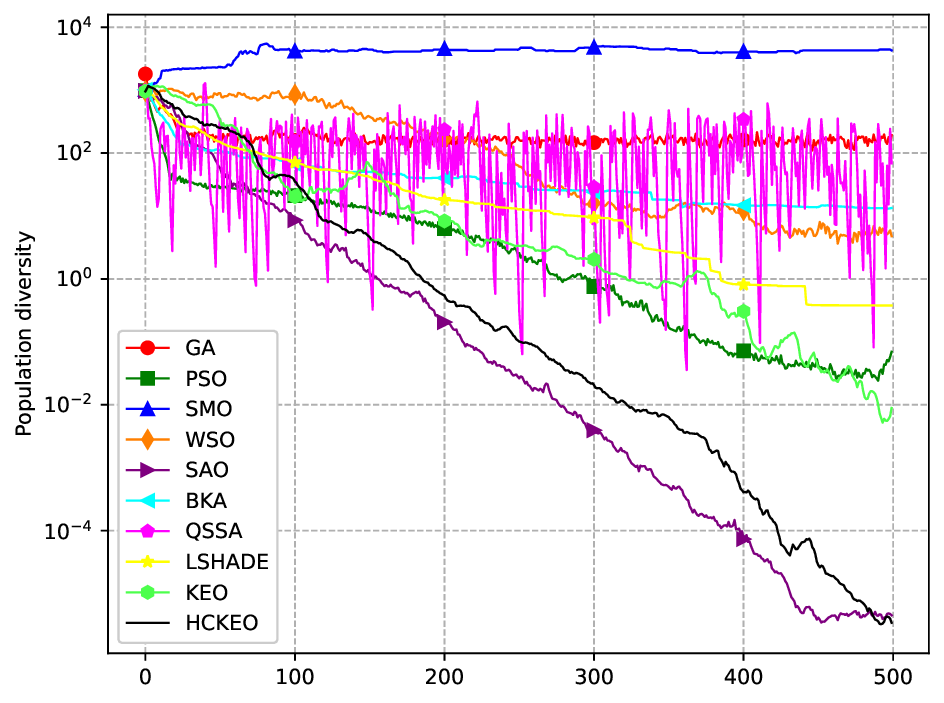}}\\
\caption{Evolution of population diversity on the CEC2022 benchmark functions (20D).}
\label{fig:div22_20}
\end{figure}

Figure \ref{fig:div22_20} presents the population diversity curves for the eight metaheuristic algorithms on the CEC2022 benchmark functions in 20 dimensions, calculated using Eq. (\ref{caldiv}).
All algorithms begin with similar high diversity levels (around $10^4$), reflecting the initial random population distribution.
As iterations progress, diversity generally decreases due to convergence toward optimal solutions, but the rate and pattern vary across algorithms and functions, highlighting differences in exploration-exploitation balance.

The proposed HCKEO algorithm consistently demonstrates improved diversity maintenance compared to its baseline KEO and other competitors in most functions.
For unimodal and basic functions (F1--F3), HCKEO exhibits a more gradual decline in diversity, remaining higher than KEO, LSHADE, and BKA in the mid-to-late stages (e.g., after 200 iterations), which supports better exploration and avoids premature convergence.
In the expanded function F4 and hybrid functions F5--F7, HCKEO shows fluctuations in diversity, particularly in F6 and F7, indicating adaptive mechanisms that reinvigorate exploration when needed, outperforming PSO and SAO, which drop rapidly to low levels ($  <10^{-4}  $).
For composition functions F8--F12, HCKEO maintains relatively stable and higher diversity longer than GA, WSO, and KEO, especially in F9, F10, and F12, where it ends at intermediate levels ($10^{-2}$ to $10^{0}$), facilitating escape from local optima.
In contrast, algorithms like LSHADE often converge too aggressively, reaching very low diversity early (e.g., $10^{-6}$ or below in F11), while HCKEO balances convergence with sustained diversity, contributing to its superior overall performance.

\subsection{Experiments on Real-World Engineering Optimization Problems}

To verify the ability of HCKEO to solve practical problems, this paper selected 10 real-world constrained engineering optimization problems.
These benchmark problems are widely used in the literature to assess the performance of MHAs under complex constraint conditions.
The selected problems include pressure vessel design, tension/compression spring design, speed reducer design, welded beam design, heat exchanger network design, three-bar truss design, step-cone pulley design, robot gripper design, planetary gear train design, and reactor network design \cite{Abhishek2020}.

Compared with standard benchmark functions, engineering optimization problems are typically characterized by nonlinear objective functions together with multiple inequality or equality constraints.
Therefore, an effective constraint-handling mechanism is required to guide the search process toward feasible regions.

In this work, a simple yet widely adopted constraint-handling approach, known as the death penalty method \cite{Coello2002}, is employed. Under this strategy, infeasible solutions are assigned a sufficiently large objective value when the optimization objective is minimization, which discourages their selection during the evolutionary process and promotes convergence toward feasible solutions.

\subsubsection{Pressure vessel design problem}

The pressure vessel design problem is a well-known constrained optimization benchmark in engineering design.
The objective is to determine an optimal configuration that minimizes the fabrication cost of a cylindrical pressure vessel while satisfying structural and volume constraints.

\begin{figure}[!ht]
	\begin{center}
		\includegraphics[scale=0.25]{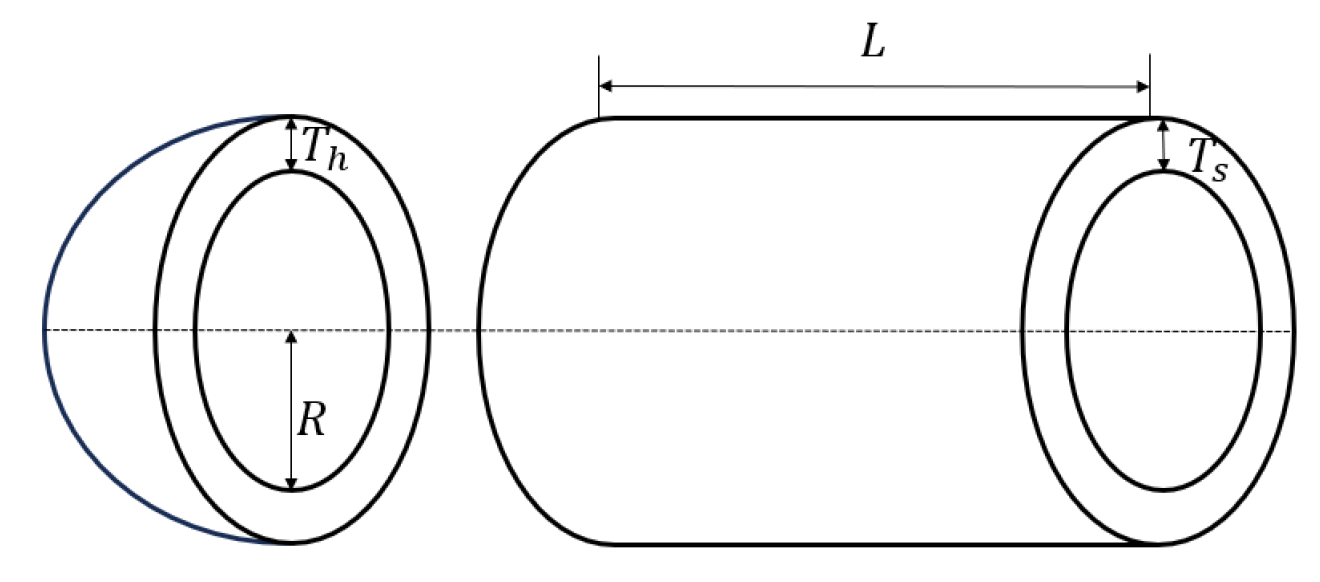}
		\caption{Pressure Vessel Schematic}
		\label{PDVFig}
	\end{center}
\end{figure}

The design involves four decision variables, namely the shell thickness ($T_s$), the head thickness ($T_h$), the inner radius ($R$), and the length of the cylindrical section ($L$). The schematic illustration of the pressure vessel is shown in Fig.~\ref{PDVFig}.

This problem's mathematical model is as follows:
\begin{itemize}
\item Consider the variables:
\[ \mathbf{X} = [x_1, x_2, x_3, x_4] = [T_s, T_h, R, L]. \]

\item The optimization problem is formulated as:
\begin{align*}
\text{Minimize:} \quad & f(\mathbf{x}) = 0.6224x_1x_3x_4 + 1.778x_2x_3^2 + 3.1661x_1^2x_4 + 19.84x_1^2x_3. \\
\text{Subject to:} \quad & g_1(\mathbf{x}) = -x_1 + 0.0193x_3 \leq 0, \\
& g_2(\mathbf{x}) = -x_2 + 0.00954x_3 \leq 0, \\
& g_3(\mathbf{x}) = -\pi x_3^2x_4 - \frac{4}{3}\pi x_3^3 + 1,\!296,\!000 \leq 0, \\
& g_4(\mathbf{x}) = x_4 - 240 \leq 0. \\
\text{With bounds:} \quad & 0 \leq x_1, x_2 \leq 100 \quad \text{and} \quad 10 \leq x_3, x_4 \leq 200.
\end{align*}
\end{itemize}

\begin{table}[!ht]
    \centering
    \renewcommand{\arraystretch}{1.2}
    \caption{Optimization results for the pressure vessel design problem.}
    \label{tab:pvd_results}
    \resizebox{\textwidth}{!}{
    \begin{tabular}{lccccc}
        \toprule
        \multirow{2}{*}{Algorithm} & \multicolumn{4}{c}{Optimal values for variables}  & \multirow{2}{*}{Optimal value} \\
        \cmidrule(lr){2-5}
         & $T_s$ & $T_h$ & $R$ & $L$ & \\
        \midrule	
        GA    & 1.263285628  & 0.670982883  & 63.88312924  & 16.02998820  & 7777.870937  \\
        PSO   & 4.189655760  & 20.40674138  & 70.52963736  & 195.9236791  & 251982.6757  \\
        SMO   & 0.852086386  & 0.431481982  & 44.14003262  & 152.9108038  & 6061.652003  \\
        WSO   & 0.778168642  & 0.384649163  & 40.31961877  & 199.9999993  & 5885.332775  \\
        SAO   & 0.894747549  & 0.442274180  & 46.35997664  & 130.1281218  & 6115.948600  \\
        BKA   & 0.778172825  & 0.384673488  & 40.31962084  & 200.0000000  & 5885.433773  \\
        QSSA   & 0.955915621  & 0.472506305  & 49.52895598  & 102.1268451  & 6263.865463  \\
        LSAHDE & 0.778170365 & 0.384650698  & 40.31970279  & 199.9990555  & 5885.343249  \\
        KEO  & 0.938591962  & 0.456198333  & 46.195713815  & 134.801356900  & 6552.296264  \\
        \textbf{HCKEO} & \textbf{0.778168641} & \textbf{0.384649163} & \textbf{40.31961872} & \textbf{200.000000} & \textbf{5885.332774} \\
        \bottomrule
    \end{tabular}}
\end{table}

\subsubsection{Tension/compression spring design problem}

\begin{figure}[!ht]
	\begin{center}
		\includegraphics[scale=0.3]{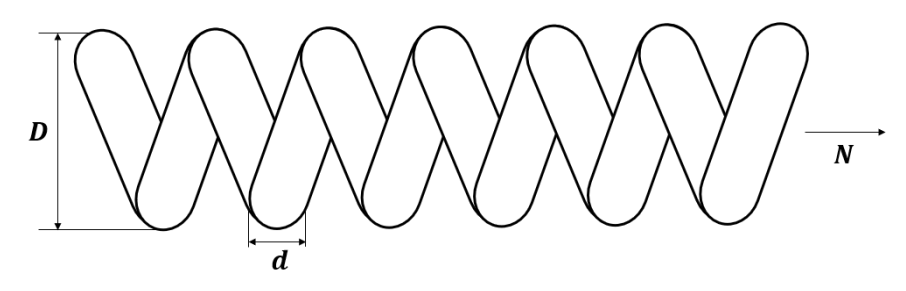}
		\caption{Pressure Vessel Schematic}
		\label{SPRINGFig}
	\end{center}
\end{figure}

The classic constraint optimization problem in this experiment is tension/compression spring design problem, whose optimization objective is to obtain the lightest spring structure while meeting different mechanical requirements.
A schematic representation of the spring structure is shown in Fig.~\ref{SPRINGFig}.

The design process involves three decision variables: the wire diameter ($d$), the mean coil diameter ($D$), and the number of active coils ($N$).
These parameters directly affect the structural strength, deformation characteristics, and vibration behavior of the spring.
During optimization, several engineering constraints must be satisfied, including limitations on shear stress, deflection, and surge frequency.
These constraints ensure that the resulting design maintains acceptable mechanical performance under operating conditions.
The mathematical model is as follows:
\begin{itemize}
\item Consider the variables:
\[ \mathbf{X} = [x_1, x_2, x_3] = [d, D, N]. \]

\item The optimization problem is formulated as:
\begin{align*}
\text{Minimize:} \quad & f(\mathbf{x}) = (x_3 + 2)x_2x_1^2. \\
\text{Subject to:} \quad & g_1(\mathbf{x}) = 1 - \frac{x_2^3x_3}{71785x_1^4} \leq 0, \\
& g_2(\mathbf{x}) = \frac{4x_2^2 - x_1x_2}{12566(x_2x_1^3)} + \frac{1}{5108x_1^2} - 1 \leq 0, \\
& g_3(\mathbf{x}) = 1 - \frac{140.45x_1}{x_2^2x_3} \leq 0, \\
& g_4(\mathbf{x}) = \frac{x_1 + x_2}{1.5} - 1 \leq 0. \\
\text{With bounds:} \quad & 0.05 \leq x_1 \leq 2 \quad 0.25 \leq x_2 \leq 1.3 \quad \text{and} \quad 2 \leq x_3 \leq 15.
\end{align*}
\end{itemize}

\begin{table}[!ht]
\centering
\caption{Optimization results for the tension/compression spring design problem.}
\label{tab:spring_results}
\begin{tabular}{lcccc}
\toprule
\multirow{2}{*}{Algorithm} & \multicolumn{3}{c}{Optimal values for variables} & \multirow{2}{*}{Optimal value} \\
\cmidrule(lr){2-4}
 & $d$ & $D$ & $N$ & \\
\midrule
GA & 0.070246839 & 0.917461209 & 2.264891525 & 0.019308532 \\
PSO & 0.067758313 & 0.750177860 & 4.444386828 & 0.022195810 \\
SMO & 0.053532910 & 0.402718561 & 9.027283172 & 0.012726585 \\
WSO & 0.052002309 & 0.364300097 & 10.85792520 & 0.012666932 \\
SAO & 0.062441594 & 0.674810594 & 3.551272119 & 0.014605699 \\
BKA & 0.055281082 & 0.449342903 & 7.389349163 & 0.012893370 \\
QSSA & 0.050000000 & 0.317399111 & 14.03125773 & 0.01272076 \\
LSAHDE & 0.051688218 & 0.356697481 & 11.29015351 & 0.01266523 \\
KEO & 0.050000000 & 0.317423896 & 14.0279719 & 0.012719153 \\
\textbf{HCKEO} & \textbf{0.051691219} & \textbf{0.35676965} & \textbf{11.28592329} & \textbf{0.012665233} \\
\bottomrule
\end{tabular}
\end{table}

\subsubsection{Speed reducer design problem}

The speed reducer design problem is a widely studied benchmark in mechanical engineering optimization.
The objective is to determine a configuration of the reducer components that minimizes the total weight of the transmission system while satisfying a set of mechanical and geometric constraints.

\begin{figure}[!ht]
	\begin{center}
		\includegraphics[scale=0.4]{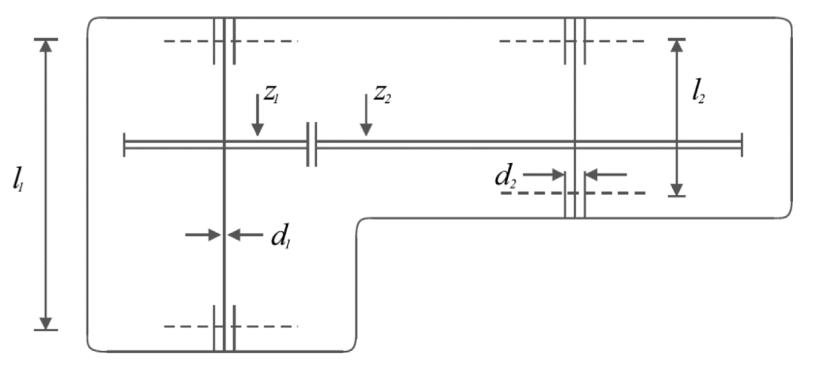}
		\caption{Pressure Vessel Schematic}
		\label{SPEEDig}
	\end{center}
\end{figure}

A schematic illustration of the speed reducer structure is presented in Fig.~\ref{SPEEDig}.
The design involves several structural parameters related to the gear train and shaft configuration, which jointly influence the strength, stability, and manufacturability of the reducer.
The corresponding mathematical formulation of the speed reducer design problem can be expressed as follows.:
\begin{itemize}
\item Consider the variables:
\[ \mathbf{X} = [x_1, x_2, x_3, x_4, x_5, x_6, x_7] = [b, m, p, l_1, l_2, d_1, d_2]. \]

\item The optimization problem is formulated as:
\begin{align*}
\text{Minimize:} \quad & f(\mathbf{x}) = 0.7854x_1x_2^2 (3.3333x_3^2 + 14.9334x_3 - 43.0934) \\
& \qquad - 1.508x_1 (x_6^2 + x_7^2) + 7.4777 (x_6^3 + x_7^3) \\
& \qquad + 0.7854(x_4x_6^2 + x_5x_7^2). \\
\text{Subject to:} \quad & \begin{aligned}
    g_1(\mathbf{x}) &= \frac{27}{x_1x_2^2x_3} - 1 \leq 0, \\
    g_2(\mathbf{x}) &= \frac{397.5}{x_1x_2^2x_3^2} - 1 \leq 0, \\
    g_3(\mathbf{x}) &= \frac{1.93x_4^3}{x_2x_3x_6^4} - 1 \leq 0, \\
    g_4(\mathbf{x}) &= \frac{1.93x_5^3}{x_2x_3x_7^4} - 1 \leq 0, \\
    g_5(\mathbf{x}) &= \frac{1}{110x_6^3} \sqrt{\left( \frac{745x_4}{x_2x_3} \right)^2 + 16.9 \times 10^6} - 1 \leq 0, \\
    g_6(\mathbf{x}) &= \frac{1}{85x_7^3} \sqrt{\left( \frac{745x_5}{x_2x_3} \right)^2 + 157.5 \times 10^6} - 1 \leq 0, \\
    g_7(\mathbf{x}) &= \frac{x_2x_3}{40} - 1 \leq 0, \\
    g_8(\mathbf{x}) &= \frac{5x_2}{x_1} - 1 \leq 0, \\
    g_9(\mathbf{x}) &= \frac{x_1}{12x_2} - 1 \leq 0, \\
     g_{10}(\mathbf{x}) &= \frac{1.5x_6 + 1.9}{x_4} - 1 \leq 0, \\
    g_{11}(\mathbf{x}) &= \frac{1.1x_7 + 1.9}{x_5} - 1 \leq 0.
\end{aligned} \\
\text{With bounds:} \quad & 2.6 \leq x_1 \leq 3.6, \quad 0.7 \leq x_2 \leq 0.8, \quad 17 \leq x_3 \leq 28, \\
& 7.3 \leq x_4 \leq 8.3, \quad 7.3 \leq x_5 \leq 8.3, \quad 2.9 \leq x_6 \leq 3.9, \\
& \text{and } 5 \leq x_7 \leq 5.5.
\end{align*}
\end{itemize}

\begin{table}[!ht]
\centering
\caption{Optimization results for the speed reducer design problem.}
\label{tab:speed_results}
\resizebox{\textwidth}{!}{
\begin{tabular}{lcccccccr}
\toprule
\multirow{2}{*}{Algorithm} & \multicolumn{7}{c}{Optimum variables} & \multirow{2}{*}{Optimum cost} \\
\cmidrule(lr){2-8}
 & $b$ & $m$ & $p$ & $l_1$ & $l_2$ & $d_1$ & $d_2$ & \\
\midrule
GA      & 3.500468106 & 0.700051444 & 17.00230823 & 7.300195341 & 7.716471632 & 3.350620984 & 5.287121979 & 2995.713913 \\
PSO     & 3.395686519 & 0.711580678 & 27.06970830 & 7.818530060 & 7.984860605 & 3.328936322 & 5.236036583 & 3.36291E+98 \\
SMO     & 3.60000000 & 0.70000000 & 15.3772941 & 7.69503073 & 8.30000000 & 3.503266761 & 5.33716 & 2847.385749 \\
WSO     & 3.500000006 & 0.700000001 & 17.0000005 & 7.30000000 & 7.715356853 & 3.350214948 & 5.286654545 & 2994.471999 \\
SAO     & 3.500011324 & 0.700000000 & 17.0000000 & 7.30000000 & 7.715327163 & 3.350222934 & 5.286654641 & 2994.477895 \\
BKA     & 3.500000124 & 0.700000000 & 17.0000000 & 7.95549509 & 7.953650139 & 3.372574947 & 5.288299759 & 3012.343349 \\
QSSA     & 3.500000138 & 0.70000000 & 17.0000000 & 7.560290190 & 8.171001745 & 3.350709256 & 5.286813266 & 3006.996505 \\
LSAHDE  & 3.500000000 & 0.70000000 & 17.0000000 & 7.30000000 & 7.715319911 & 3.350214666 & 5.286654464 & 2994.471066 \\
KEO    & 3.500000003 & 0.70000000 & 17.0000000 & 7.30000000 & 7.808771427 & 3.350214673 & 5.286686233 & 2996.542621 \\
\textbf{HCKEO}   & \textbf{3.50000000} & \textbf{0.70000000} & \textbf{17.0000000} & \textbf{7.30000000} & \textbf{7.71531991147825} & \textbf{3.35021466609645} & \textbf{5.28665446498022} & \textbf{2994.47106614682} \\
\bottomrule
\end{tabular}}
\end{table}

\subsubsection{Welded beam design problem}

\begin{figure}[!ht]
	\begin{center}
		\includegraphics[scale=0.5]{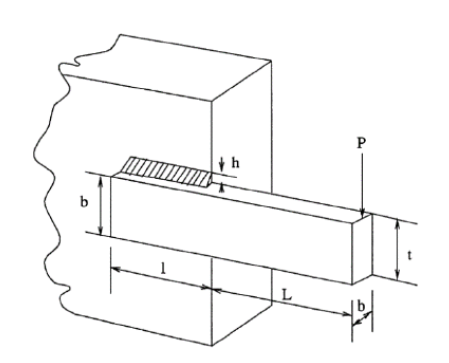}
		\caption{Welded Beam Design}
		\label{fig:Welded Beam Design}
	\end{center}
\end{figure}

Welded beam design focuses on minimizing costs by optimizing the design of a welded beam.
The design process encompasses four variables: height(h), length(l), thickness(t) and breadth(b).
The schematic representation of the welded beam is illustrated in Fig. \ref{fig:Welded Beam Design}.
The mathematical formulation of the optimization problem associated with welded beam design is presented as follows:
\begin{itemize}
\item Consider the variables:
\[ \mathbf{X} = [x_1, x_2, x_3, x_4] = [h, l, t, b]. \]

\item The optimization problem is formulated as:
\begin{align*}
\text{Minimize:} \quad & f(\mathbf{x}) = 0.04811 x_3 x_4 (x_2 + 14) + 1.10471 x_1^2 x_2 \\
\text{Subject to:} \quad & \begin{aligned}
g_1(\mathbf{x}) &= x_1 - x_4 \leq 0 \\
g_2(\mathbf{x}) &= \delta(\mathbf{x}) - \delta_{\text{max}} \leq 0 \\
g_3(\mathbf{x}) &= P - P_c(\mathbf{x}) \leq 0 \\
g_4(\mathbf{x}) &= \tau(\mathbf{x}) - \tau_{\text{max}} \leq 0 \\
g_5(\mathbf{x}) &= \sigma(\mathbf{x}) - \sigma_{\text{max}} \leq 0
\end{aligned} \\
& \begin{aligned}
    g_5(\mathbf{x}) &= \frac{1}{110x_6^3} \sqrt{\left( \frac{745x_4}{x_2x_3} \right)^2 + 16.9 \times 10^6} - 1 \leq 0, \\
    g_6(\mathbf{x}) &= \frac{1}{85x_7^3} \sqrt{\left( \frac{745x_5}{x_2x_3} \right)^2 + 157.5 \times 10^6} - 1 \leq 0,
\end{aligned} \\
\text{With bounds:} \quad & 0.125 \leq x_1 \leq 2, \quad 0.1 \leq x_2, x_3 \leq 10, \quad 0.1 \leq x_4 \leq 2.
\end{align*}
\end{itemize}

\begin{table}[!ht]
    \centering
    \renewcommand{\arraystretch}{1.2}
    \caption{Optimization results for the welded beam design problem.}
    \label{tab:wbd_results}
    \resizebox{\textwidth}{!}{
    \begin{tabular}{lccccc}
        \toprule
        \multirow{2}{*}{Algorithm} & \multicolumn{4}{c}{Optimal values for variables}  & \multirow{2}{*}{Optimal value} \\
        \cmidrule(lr){2-5}
         & $T_s$ & $T_h$ & $R$ & $L$ & \\
        \midrule
        GA    & 0.310085068845545 & 2.59294439090531 & 6.44114431410583 & 0.405346825609125 & 2.35966918896065  \\
        PSO   & 0.437019933412312&2.28566439992982&6.44674625531403&0.742669104399893  & 4.23349660745509  \\
        SMO   & 0.348377307815625 & 3.96309536657826 & 8.58108179456618 & 0.578412600289239  & 4.82074999933099  \\
        WSO   & 0.205729649397326&3.01039888580124&9.03662333452380&0.205729677172381  & 1.66218893444694  \\
        SAO   & 0.205729639784955&3.01039750927435&9.03662391085411&0.205729639788443  & 1.66218887323887  \\
        BKA   & 0.205084248243580&3.02128289060828&9.03626598683563&0.205849313274611  & 1.66361176407930  \\
        QSSA   & 0.125000000000000&5.23514231996968&9.03662197950160&0.205729727786598  & 1.81078208357137  \\
        LSAHDE & 0.205729639780189&3.01039750950179&9.03662391035895&0.205729639786136  & 1.66218887316290  \\
        KEO  & 0.205727608117569&3.01030933005766&9.03700827772803&0.205727722212788  & 1.66222461471345  \\
        \textbf{HCKEO} & \textbf{0.205729639786079} & \textbf{0.205729639786079} & \textbf{9.03662391035763} & \textbf{0.205729639786080} & \textbf{1.66218887315602} \\
        \bottomrule
    \end{tabular}}
\end{table}

\subsubsection{Heat Exchanger Network Design (Case 2)}

The Heat Exchanger Network (HEN) design problem is a representative large-scale nonlinear programming (NLP) problem widely studied in chemical engineering optimization.
The objective is to minimize the total investment and operating costs of the heat exchanger network while satisfying energy balance constraints and heat transfer requirements.
Due to the presence of multiple equality constraints, nonlinear logarithmic terms, and decision variables with significantly different scales, this problem poses considerable challenges for metaheuristic optimization algorithms.

In this study, the considered HEN design model consists of 11 continuous decision variables and 9 nonlinear equality constraints.
The mathematical formulation involves logarithmic functions related to temperature differences and heat transfer areas, making the feasible region highly nonlinear and sensitive to constraint violations.
This problem's mathematical model is as follows:
\begin{itemize}
\item The decision variable vector is defined as
\[
\mathbf{X} = [x_1, x_2, x_3, x_4, x_5, x_6, x_7, x_8, x_9, x_{10}, x_{11}],
\]
\item the optimization problem is formulated as follows:
\begin{align*}
\text{Minimize:} \quad
& f(\mathbf{x}) = \left( \frac{x_1}{120x_4} \right)^{0.6}
+ \left( \frac{x_2}{80x_5} \right)^{0.6}
+ \left( \frac{x_3}{40x_6} \right)^{0.6}, \\
\text{Subject to:} \quad
& \begin{aligned}
h_1(\mathbf{x}) &= x_1 - 10^4(x_7 - 100) = 0, \\
h_2(\mathbf{x}) &= x_2 - 10^4(x_8 - x_7) = 0, \\
h_3(\mathbf{x}) &= x_3 - 10^4(500 - x_8) = 0, \\
h_4(\mathbf{x}) &= x_1 - 10^4(300 - x_9) = 0, \\
h_5(\mathbf{x}) &= x_2 - 10^4(400 - x_{10}) = 0, \\
h_6(\mathbf{x}) &= x_3 - 10^4(600 - x_{11}) = 0, \\
h_7(\mathbf{x}) &= x_4 \ln(x_9 - 100) - x_4 \ln(300 - x_7) - x_9 - x_7 + 400 = 0, \\
h_8(\mathbf{x}) &= x_5 \ln(x_{10} - x_7) - x_5 \ln(400 - x_8) - x_{10} + x_7 - x_8 + 400 = 0, \\
h_9(\mathbf{x}) &= x_6 \ln(x_{11} - x_8) - x_6 \ln(100) - x_{11} + x_8 + 100 = 0,
\end{aligned} \\
\text{With bounds:} \quad
& \begin{aligned}
&0 \leq x_1 \leq 81.9 \times 10^4, \quad
10^4 \leq x_2 \leq 113.1 \times 10^4, \\
&10^4 \leq x_3 \leq 205 \times 10^4, \quad
0 \leq x_4, x_5, x_6 \leq 5.074 \times 10^{-2}, \\
&100 \leq x_7 \leq 200, \quad
100 \leq x_8, x_9, x_{10} \leq 300, \quad
100 \leq x_{11} \leq 400.
\end{aligned}
\end{align*}
\end{itemize}

\begin{figure}[!ht]
\centering
\begin{tikzpicture}[
    scale=1.0,
    line cap=round,
    line join=round,
    >=Stealth,
    hx/.style={draw, rectangle, minimum width=2.6cm, minimum height=1.2cm, line width=1pt},
    stream/.style={line width=1pt},
    label/.style={font=\small}
]

\node[hx] (HX1) at (0,0)   {HX$_1$};
\node[hx] (HX2) at (4,0)   {HX$_2$};
\node[hx] (HX3) at (8,0)   {HX$_3$};

\draw[stream, ->] (-1,0.6) -- (HX1.west |- 0,0.6);
\draw[stream, ->] (HX1.east |- 0,0.6) -- (HX2.west |- 0,0.6);
\draw[stream, ->] (HX2.east |- 0,0.6) -- (HX3.west |- 0,0.6);
\draw[stream, ->] (HX3.east |- 0,0.6) -- (9.5,0.6);

\node[label] at (-1.2,0.6) {Hot in};
\node[label] at (9.8,0.6) {Hot out};

\draw[stream, <-] (-1,-0.6) -- (HX1.west |- 0,-0.6);
\draw[stream, <-] (HX1.east |- 0,-0.6) -- (HX2.west |- 0,-0.6);
\draw[stream, <-] (HX2.east |- 0,-0.6) -- (HX3.west |- 0,-0.6);
\draw[stream, <-] (HX3.east |- 0,-0.6) -- (9.5,-0.6);

\node[label] at (-1.3,-0.6) {Cold out};
\node[label] at (9.8,-0.6) {Cold in};

\node[label] at (0,-1.2) {$x_1,\; x_4$};
\node[label] at (4,-1.2) {$x_2,\; x_5$};
\node[label] at (8,-1.2) {$x_3,\; x_6$};

\node[label] at (1.3,0.95) {$x_7$};
\node[label] at (5.3,0.95) {$x_8$};
\node[label] at (1.3,-1.0) {$x_9$};
\node[label] at (5.3,-1.0) {$x_{10}$};
\node[label] at (8.9,-1.0) {$x_{11}$};

\end{tikzpicture}

\caption{Schematic diagram of the Heat Exchanger Network (HEN) design problem for Case~2.
The network consists of three heat exchangers connected in series, with a hot stream flowing from left to right and a cold stream flowing in the opposite direction.
The decision variables $x_1$--$x_3$ represent the heat duties, $x_4$--$x_6$ denote heat transfer area-related parameters, and $x_7$--$x_{11}$ correspond to intermediate temperature variables.}
\label{fig:HEND2}
\end{figure}
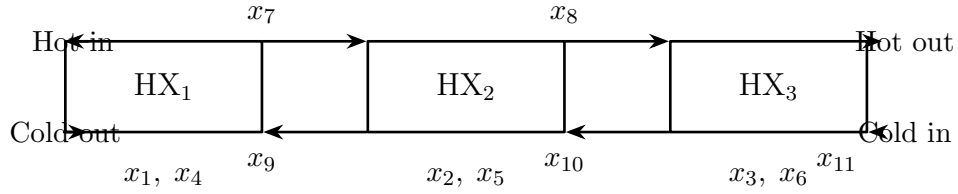

A schematic illustration of the Heat Exchanger Network configuration for Case 2 is shown in Fig.~\ref{fig:HEND2}.

Table~\ref{tab:heat_exchanger_results} summarizes the results of HCKEO and other comparative algorithms in this experiment, which lists the optimal design variables together with their corresponding objective function values.

\begin{table}[!ht]
    \centering
    \renewcommand{\arraystretch}{1.2}
    \caption{Optimization results for the heat exchanger network design (case 2) problem.}
    \label{tab:heat_exchanger_results}
    \resizebox{\textwidth}{!}{
    \begin{tabular}{lccccc}
        \toprule
        \multirow{2}{*}{Algorithm} & \multicolumn{4}{c}{Optimal values for variables} & \multirow{2}{*}{Optimal value} \\
        \cmidrule(lr){2-5}
         & $T_s$ & $T_h$ & $R$ & $L$ & \\
        \midrule
        GA
        & 809461.060933405 & 1130854.62691484 & 2049565.26951451 & 0.0507391201276824
        & $5.29866533597473\times10^{15}$ \\

        PSO
        & 746016.607355141 & 1009450.12983996 & 2050000.00000000 & 0.0266502188375417
        & $8.90122617984458\times10^{18}$ \\

        SMO
        & 818493.135665050 & 1144403.23753470 & 2012052.99661648 & 0.0300720209858634
        & 3.49188577751531e+16 \\

        WSO
        & 818999.998409066 & 1130999.99899977 & 2049999.99966074 & 0.0477233533823672
        & 10244.2763837370 \\

        SAO
        & 819000.000000000 & 1131000.00000000 & 2050000.00000000 & 0.0507400000000000
        & 7049.03695431410 \\

        BKA
        & 818999.964829814 & 1130998.66989143 & 2049999.99999964 & 0.0507127174194137
        & 62322786.4893148 \\

        QSSA
        & 818290.158391332 & 1130999.99891003 & 2047742.80404034 & 0.0506764712460471
        & 2.42460500010439e+19 \\

        LSAHDE
        & 819000.000000000 & 1131000.00000000 & 2050000.00000000 & 0.0507400000000000
        & 7049.03695431410 \\

        KEO
        & 819000.000000000 & 1131000.00000000 & 2050000.00000000 & 0.0507400000000000
        & 7049.03695431410 \\

        \textbf{HCKEO}
        & \textbf{819000.000000000} & \textbf{1131000.00000000}
        & \textbf{2050000.00000000} & \textbf{0.0507400000000000}
        & \textbf{7049.03695431410} \\
        \bottomrule
    \end{tabular}}
\end{table}

\subsubsection{Three-bar truss design problem}

This case study represents a classic structural engineering task characterized by a linear objective function coupled with highly non-linear stress constraints.
Although the search space is limited to two dimensions, the "accidented" nature of the constrained feasible region makes it a challenging benchmark for testing the precision of the HCKEO algorithm.
It evaluates whether the algorithm can push the solution toward the constraint boundary without violating the safety limits ($g_1, g_2, g_3$).
This problem’s mathematical model is as follows:

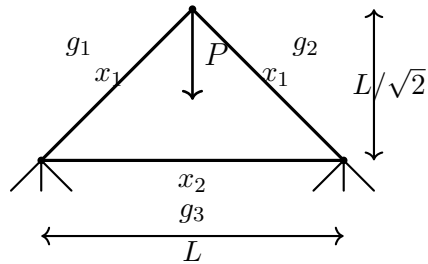
\begin{figure}[!ht]
\centering
\begin{tikzpicture}[scale=2.0, line cap=round, line join=round]

\tikzset{
    truss/.style={line width=1.2pt},
    support/.style={line width=0.8pt},
    load/.style={->, line width=1.0pt},
    dim/.style={<->, line width=0.8pt},
    label/.style={font=\small},
}

\coordinate (A) at (-1,0);   
\coordinate (B) at (1,0);    
\coordinate (C) at (0,1);    

\draw[truss] (A) -- (C);    
\draw[truss] (C) -- (B);    
\draw[truss] (A) -- (B);    

\draw[support] (A) -- ++(-0.2,-0.2);
\draw[support] (A) -- ++(0.0,-0.2);
\draw[support] (A) -- ++(0.2,-0.2);

\draw[support] (B) -- ++(-0.2,-0.2);
\draw[support] (B) -- ++(0.0,-0.2);
\draw[support] (B) -- ++(0.2,-0.2);

\draw[load] (C) -- ++(0,-0.6) node[midway,right] {$P$};

\node[label] at (-0.55,0.55) {$x_1$};
\node[label] at (0.55,0.55) {$x_1$};
\node[label] at (0,-0.15) {$x_2$};

\node[label] at (-0.75,0.75) {$g_1$};
\node[label] at (0.75,0.75) {$g_2$};
\node[label] at (0,-0.35) {$g_3$};

\draw[dim] (-1,-0.5) -- (1,-0.5);
\node[label] at (0,-0.6) {$L$};

\draw[dim] (1.2,0) -- (1.2,1);
\node[label] at (1.3,0.5) {$L/\sqrt{2}$};

\fill (A) circle (0.025);
\fill (B) circle (0.025);
\fill (C) circle (0.025);

\end{tikzpicture}
\caption{Schematic diagram of the three-bar truss design problem.
The structure consists of two symmetric diagonal members and one horizontal member.
The design variables $x_1$ and $x_2$ represent the cross-sectional areas of the diagonal and horizontal bars, respectively.
A vertical load $P$ is applied at the top joint, and stress constraints $g_1$, $g_2$, and $g_3$ are imposed on the corresponding members.}
\label{fig:three_bar_truss}
\end{figure}

\begin{itemize}
\item Consider the variables:
\[ \mathbf{X} = [x_1, x_2]. \]

\item The optimization problem is formulated as:
\begin{align*}
\text{Minimize:} \quad
& f(\mathbf{x}) = L(x_2 + 2\sqrt{2}x_1), \\
\text{Subject to:} \quad
& \begin{aligned}
g_1(\mathbf{x}) &= \frac{x_2}{2x_2x_1 + \sqrt{2}x_1^2} P - \sigma \leq 0 \\
g_2(\mathbf{x}) &= \frac{x_2 + \sqrt{2}x_1}{2x_2x_1 + \sqrt{2}x_1^2} P - \sigma \leq 0 \\
g_3(\mathbf{x}) &= \frac{1}{x_1 + \sqrt{2}x_2} P - \sigma \leq 0,
\end{aligned} \\
\text{Parameters:} \quad
& \begin{aligned}
L = 100, \quad P = 2, \quad \sigma = 2.
\end{aligned} \\
\text{With bounds:} \quad
& \begin{aligned}
0 \leq x_1, x_2 \leq 1.
\end{aligned}
\end{align*}
\end{itemize}

\begin{table}[!ht]
\centering
\caption{Optimization results for the three-bar truss design problem.}
\label{tab:three_bar_comparison}
\begin{tabular}{lccc}
\toprule
\multirow{2}{*}{Algorithm} & \multicolumn{2}{c}{Optimal values for variables} & \multirow{2}{*}{Optimal value} \\
\cmidrule(lr){2-3}
 & $x_1$ (A1) & $x_2$ (A2) &  \\
\midrule
GA      & 0.79099573 & 0.40234393 & 263.96176997 \\
PSO     & 0.79563224 & 0.50790109 & 275.82889047 \\
SMO     & 0.78867514 & 0.40824827 & 263.89584338 \\
WSO     & 0.78867504 & 0.40824855 & 263.89584338 \\
SAO     & 0.78648329 & 0.41448333 & 263.89940072 \\
BKA     & 0.78811921 & 0.40982297 & 263.89607093 \\
QSSA    & 0.78493853 & 0.41892094 & 263.90623799 \\
LSAHDE  & 0.78867513 & 0.40824829 & 263.89584338 \\
KEO     & 0.78866712 & 0.40827095 & 263.89584342 \\
\textbf{HCKEO} & \textbf{0.78867513} & \textbf{0.40824829} & \textbf{263.89584338} \\
\bottomrule
\end{tabular}
\end{table}

\subsubsection{Step-cone pulley problem}

The step-cone pulley problem is a representative mechanical design benchmark featuring a highly constrained nonlinear optimization landscape.
The main difficulty arises from the requirement that a single belt must maintain identical effective lengths across four distinct transmission stages, which is enforced through a set of equality constraints.
From an optimization standpoint, these constraints restrict the feasible solutions to a very narrow manifold embedded in a five-dimensional decision space.
In addition, the presence of transcendental functions, including exponential and inverse trigonometric terms, in the stress and power constraints substantially increases the numerical complexity and the risk of premature convergence to local optima.

\begin{figure}[!ht]
\centering
\begin{tikzpicture}[
    scale=1.0,
    line cap=round,
    line join=round,
    >=Stealth,
    pulley/.style={draw, circle, line width=1pt},
    belt/.style={line width=1pt},
    label/.style={font=\small}
]

\coordinate (L) at (0,0);
\coordinate (R) at (6,0);

\draw[pulley] (L) circle (0.5);
\draw[pulley] (L) circle (0.7);
\draw[pulley] (L) circle (0.9);
\draw[pulley] (L) circle (1.1);

\draw[pulley] (R) circle (0.6);
\draw[pulley] (R) circle (0.8);
\draw[pulley] (R) circle (1.0);
\draw[pulley] (R) circle (1.2);

\draw[belt] ($(L)+(0,1.1)$) -- ($(R)+(0,1.2)$);
\draw[belt] ($(L)+(0,-1.1)$) -- ($(R)+(0,-1.2)$);

\node[label] at (-0.9,1.15) {$d_1$};
\node[label] at (-1.0,0.95) {$d_2$};
\node[label] at (-1.1,0.75) {$d_3$};
\node[label] at (-1.2,0.55) {$d_4$};

\draw[<->, line width=0.8pt] (3,1.35) -- (3,1.65);
\node[label] at (3,1.8) {$w$};

\draw[<->, line width=0.8pt] (0,-1.6) -- (6,-1.6);
\node[label] at (3,-1.8) {$a$};

\draw[->, line width=0.8pt] (0.37,0.15) arc (20:320:0.4);
\draw[->, line width=0.8pt] (6.37  ,0.15) arc (20:320:0.4);

\node[label] at (0,1.9) {$N$};
\node[label] at (6,1.9) {$N_i$};

\end{tikzpicture}
\caption{Schematic diagram of the step-cone pulley design problem.
A single belt with width $w$ transmits power between two stepped pulleys with four diameter levels $d_1$--$d_4$.
The center distance $a$ is fixed, and belt-length consistency constraints enforce identical belt lengths across all transmission stages.}
\label{fig:step_cone_pulley}
\end{figure}
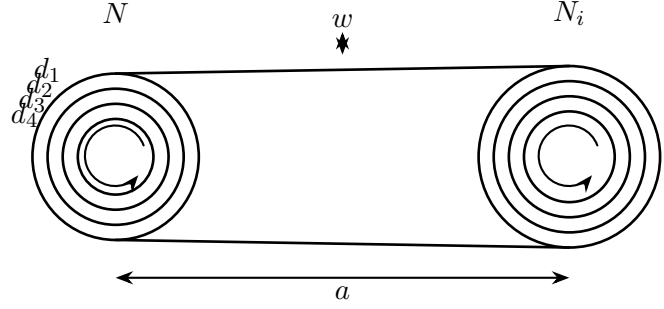

Fig.~\ref{fig:step_cone_pulley} presents a schematic illustration of the step-cone pulley system.
A single belt with width $w$ transmits power between two stepped pulleys with four diameter levels $d_1$--$d_4$ under different transmission ratios.
The center distance between the pulleys is fixed, and belt-length consistency is enforced across all transmission stages.
The problem’s mathematical model is as follows:
\begin{itemize}
\item The decision variable vector is defined as
\[
\mathbf{X} = [d_1, d_2, d_3, d_4, w].
\]

\item The optimization problem is formulated as follows:
\begin{align*}
\text{Minimize:} \quad
& f(\mathbf{x}) = \rho w \frac{\pi}{4} \sum_{i=1}^{4} d_i^2 \left[ 1 + \left( \frac{N_i}{N} \right)^2 \right], \\
\text{Subject to:} \quad
& \begin{aligned}
h_1(\mathbf{x}) &= C_1 - C_2 = 0, \\
h_2(\mathbf{x}) &= C_1 - C_3 = 0, \\
h_3(\mathbf{x}) &= C_1 - C_4 = 0, \\
g_{1,2,3,4}(\mathbf{x}) &= R_i - 2 \leq 0, \\
g_{5,6,7,8}(\mathbf{x}) &= 0.75 \times 745.7 - P_i \leq 0,
\end{aligned} \\
\text{Parameters:} \quad
& \begin{aligned}
C_i &= \frac{\pi d_i}{2} \left( 1 + \frac{N_i}{N} \right)
      + \frac{(N_i/N - 1)^2 d_i^2}{4a} + 2a, \\
R_i &= \exp \left( \mu \left( \pi - 2\arcsin\left( \frac{(N_i/N - 1)d_i}{2a} \right) \right) \right), \\
P_i &= s t w (1 - R_i) \frac{\pi d_i N}{60000}, \\
\rho &= 7200, \quad N = 250, \quad N_i = \{100,150,200,250\}, \\
a &= 3, \quad \mu = 0.35, \quad s = 1.75, \quad t = 8,
\end{aligned} \\
\text{With bounds:} \quad
& 1 \leq d_i, w \leq 100.
\end{align*}
\end{itemize}

The optimization results of HCKEO and the seven comparison algorithms are summarized in Table~\ref{tab:pulley_comparison}.

\begin{table}[!ht]
\centering
\caption{Optimization results for the step-cone pulley problem.}
\label{tab:pulley_comparison}
\begin{tabular}{lcccccc}
\toprule
\multirow{2}{*}{Algorithm} & \multicolumn{5}{c}{Optimal values for variables} & \multirow{2}{*}{Optimal value} \\
\cmidrule(lr){2-6}
 & $d_1$ & $d_2$ & $d_3$ & $d_4$ & $w$ & \\
\midrule
GA      & 39.3750 & 54.1773 & 72.2306 & 86.6130 & 88.0906 & $6.60 \times 10^{91}$ \\
PSO     & 45.6943 & 60.0000 & 90.0000 & 90.0000 & 90.0000 & $9.41 \times 10^{97}$ \\
SMO     & 40.7969 & 56.1402 & 74.8475 & 89.7363 & 89.9784 & 18.1448 \\
WSO     & 40.0594 & 55.1246 & 73.4935 & 88.1145 & 87.6804 & $8.80 \times 10^{73}$ \\
SAO     & 40.9086 & 56.2941 & 75.0525 & 89.9819 & 84.4974 & 17.1329 \\
BKA     & 40.7671 & 56.0992 & 74.7928 & 89.6708 & 88.9354 & 17.9082 \\
QSSA    & 40.5094 & 55.7443 & 74.3196 & 89.1040 & 88.9900 & $1.11 \times 10^{78}$ \\
LSAHDE  & 39.6825 & 54.6056 & 72.8015 & 87.2855 & 88.4048 & $4.63 \times 10^{78}$ \\
KEO     & 40.4146 & 55.6138 & 74.1456 & 88.8956 & 85.5331 & 16.9265 \\
\textbf{HCKEO} & \textbf{40.0141} & \textbf{55.0621} & \textbf{73.4102} & \textbf{88.0147} & \textbf{86.3916} & \textbf{16.7591} \\
\bottomrule
\end{tabular}
\end{table}

\subsubsection{Robot gripper problem}

The robot gripper problem is a classical benchmark in mechanical linkage synthesis, widely used to evaluate the performance of optimization algorithms under strong geometric and nonlinear constraints.
The primary objective is to minimize the fluctuation of the gripping force by reducing the difference between its maximum and minimum values over the allowable operating displacement range of the gripper.
This problem is particularly challenging due to the following characteristics:
\begin{enumerate}
\item \textit{Transcendental nonlinearity}: The mathematical model contains nested inverse trigonometric functions ($\cos^{-1}$, $\tan^{-1}$) combined with sinusoidal terms, which imposes high requirements on numerical stability during optimization.
\item \textit{Geometric feasibility}: The design variables must strictly satisfy triangle inequality conditions to ensure the validity of angular computations; otherwise, the objective function becomes undefined.
\item \textit{Strongly coupled constraints}: The displacement constraints and force requirements are highly interdependent, resulting in a narrow and fragmented feasible region within the search space.
\end{enumerate}

The problem’s mathematical model is as follows:
\begin{itemize}

\item The decision variable vector is defined as
\[
\mathbf{X} = [a, b, c, e, f, l, \delta].
\]

\item The optimization problem is formulated as follows:
\begin{align*}
\text{Minimize:} \quad
& f(\mathbf{x}) = \max_{z} F_k(\mathbf{x}, z) - \min_{z} F_k(\mathbf{x}, z), \\
\text{Subject to:} \quad
& \begin{aligned}
g_1(\mathbf{x}) &= Y_{min} - y(\mathbf{x}, Z_{max}) \leq 0, \\
g_2(\mathbf{x}) &= -y(\mathbf{x}, Z_{max}) \leq 0, \\
g_3(\mathbf{x}) &= y(\mathbf{x}, 0) - Y_{max} \leq 0, \\
g_4(\mathbf{x}) &= Y_G - y(\mathbf{x}, 0) \leq 0, \\
g_5(\mathbf{x}) &= (a+b)^2 - l^2 - e^2 \leq 0, \\
g_6(\mathbf{x}) &= (a-e)^2 + (l-Z_{max})^2 - b^2 \leq 0, \\
g_7(\mathbf{x}) &= Z_{max} - l \leq 0.
\end{aligned} \\
\text{Parameters:} \quad
& \begin{aligned}
F_k &= \frac{P b \sin(\alpha + \beta + \delta)}{2 c \cos \alpha}, \quad y(\mathbf{x}, z) = 2(f + e + c \sin(\beta + \delta)), \\
\alpha &= \cos^{-1} \left( \frac{a^2 + g^2 - b^2}{2ag} \right) + \phi, \quad \beta = \cos^{-1} \left( \frac{b^2 + g^2 - a^2}{2bg} \right) - \phi, \\
g &= \sqrt{e^2 + (z-l)^2}, \quad \phi = \tan^{-1} \left( \frac{e}{l-z} \right), \\
Y_{min} &= 50, Y_{max} = 100, Y_G = 150, Z_{max} = 100, P = 100
\end{aligned} \\
\text{With bounds:} \quad
& \begin{aligned}
0 \le a, b, f \le 150, \\ 100 \le c \le 200, \\ 0 \le e \le 50 \\
1 \le \delta \le 3.14, \\ 100 \le l \le 300
\end{aligned}
\end{align*}
\end{itemize}

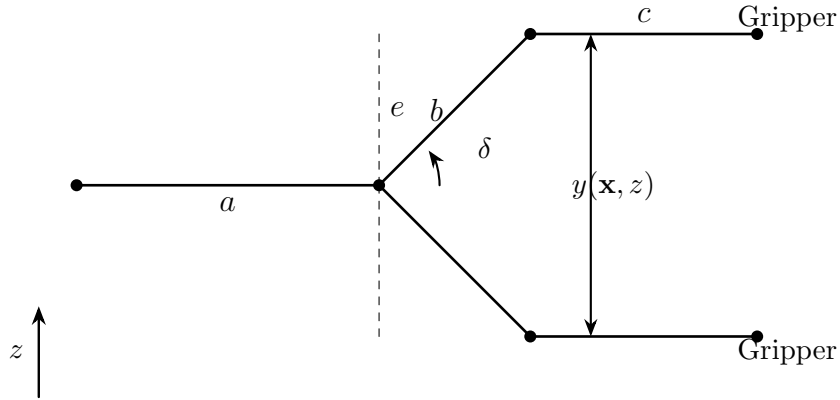
\begin{figure}[!ht]
\centering
\begin{tikzpicture}[
    scale=1.0,
    line cap=round,
    line join=round,
    >=Stealth,
    link/.style={line width=1pt},
    joint/.style={fill=black, circle, inner sep=1.6pt},
    label/.style={font=\small}
]

\coordinate (O) at (0,0);          
\coordinate (A) at (4,0);          
\coordinate (B) at (6,2);          
\coordinate (C) at (6,-2);         
\coordinate (D) at (9,2);          
\coordinate (E) at (9,-2);         

\draw[link] (O) -- (A) node[midway,below] {$a$};
\draw[link] (A) -- (B) node[midway,left] {$b$};
\draw[link] (A) -- (C);
\draw[link] (B) -- (D) node[midway,above] {$c$};
\draw[link] (C) -- (E);

\draw[dashed] (A) -- ++(0,2) node[midway,right] {$e$};
\draw[dashed] (A) -- ++(0,-2);

\draw[<->, line width=0.8pt] (6.8,-2) -- (6.8,2);
\node[label] at (7.1,0) {$y(\mathbf{x},z)$};

\draw[->, line width=0.9pt] (-0.5,-2.8) -- (-0.5,-1.6);
\node[label] at (-0.8,-2.2) {$z$};

\foreach \p in {O,A,B,C,D,E}
    \node[joint] at (\p) {};

\draw[->, line width=0.8pt] (4.8,0) arc (0:35:0.8);
\node[label] at (5.4,0.5) {$\delta$};

\node[label] at (9.4,2.2) {Gripper};
\node[label] at (9.4,-2.2) {Gripper};

\end{tikzpicture}

\caption{Schematic diagram of the robot gripper mechanism.
The design variables \(a,b,c,e,f,l,\delta\) define the geometric configuration of the linkage.
The vertical opening distance \(y(\mathbf{x},z)\) varies with the prismatic displacement \(z\), and the objective is to minimize the fluctuation of gripping force over the operating range.}
\label{fig:robot_gripper}
\end{figure}

Figure~\ref{fig:robot_gripper} illustrates a schematic diagram of the robot gripper mechanism.
The design variables \(a, b, c, e, f, l,\) and \(\delta\) define the geometric configuration of the linkage system.
The vertical opening distance \(y(\mathbf{x},z)\) varies with the prismatic displacement \(z\), and the optimization objective is to minimize the variation of the gripping force throughout the entire operating range.

Table~\ref{tab:gripper_results} summarizes the experimental results obtained by eight optimization algorithms on the robot gripper problem.

\begin{table}[!ht]
\centering
\caption{Optimization results for the robot gripper problem.}
\label{tab:gripper_results}
\resizebox{\textwidth}{!}{
\begin{tabular}{lcccccccc}
\toprule
\multirow{2}{*}{Algorithm} & \multicolumn{7}{c}{Optimal values for variables} & \multirow{2}{*}{Optimal value} \\
\cmidrule(lr){2-8}
 & $a$ & $b$ & $c$ & $e$ & $f$ & $l$ & $\delta$ & \\
\midrule
GA & 142.2758 & 119.1540 & 148.2965 & 15.5537 & 80.6864 & 169.1391 & 2.4368 & 5.7103 \\
PSO & 150.0000 & 129.1068 & 200.0000 & 46.8871 & 114.3722 & 217.2609 & 3.1400 & 67.3019 \\
SMO & 280.6236 & 178.6229 & 264.8124 & 86.2115 & 173.4440 & 268.9487 & 3.1202 & 2.0257 \\
WSO & 149.7431 & 146.4406 & 199.7222 & 2.5334 & 138.6806 & 128.3933 & 2.3743 & 2.9272 \\
SAO & 149.8850 & 149.6754 & 199.9815 & 0.0918 & 149.6032 & 100.9752 & 2.2965 & 2.5474 \\
BKA & 149.9934 & 149.8454 & 184.2192 & 0.0097 & 53.0803 & 102.9344 & 1.8082 & 2.8002 \\
QSSA & 149.9999 & 146.1486 & 176.4694 & 3.2526 & 26.7694 & 123.8566 & 1.7483 & 3.3371 \\
LSHADE & 148.2569 & 141.4447 & 200.0000 & 6.6377 & 150.0000 & 104.8494 & 2.3948 & 2.7165 \\
KEO & 149.9990 & 149.8651 & 200.0000 & 0.0167 & 150.0000 & 100.9442 & 2.2975 & 2.5440 \\
\textbf{HCKEO} & \textbf{149.9383} & \textbf{148.9958} & \textbf{200.0000} & \textbf{0.8335} & \textbf{150.0000} & \textbf{100.0000} & \textbf{2.3002} & \textbf{2.5338} \\
\bottomrule
\end{tabular}}
\end{table}

\subsubsection{Planetary Gear Train Design}

The planetary gear train design problem aims to minimize the maximum deviation between the achieved gear ratios and their prescribed target values in an automatic transmission system.
This benchmark problem is widely used in automotive engineering, as it determines the optimal combination of gear tooth numbers and modules to ensure accurate speed ratios and efficient torque transmission.
The design space is mixed-integer and discrete in nature, consisting of six integer variables representing gear tooth numbers, one discrete variable indicating the number of planet gears, and two discrete variables corresponding to gear modules.
A total of eleven constraints are imposed to satisfy geometric feasibility and assembly requirements.
The problem’s mathematical model is as follows:
\begin{itemize}

\item The decision variable vector is defined as
\[
\mathbf{x} = [p, N_6, N_5, N_4, N_3, N_2, N_1, m_3, m_1],
\]
where $N_i$ denotes the number of teeth of the corresponding gear, $p$ is the number of planet gears, and $m_1$, $m_3$ are the selected gear modules.

\item The optimization problem is formulated as
\begin{align*}
\text{Minimize:} \quad
& f(\mathbf{x}) = \max_{k=1,\dots,R} \left| i_k - i_{0k} \right|, \\
\text{Subject to:} \quad
& \begin{aligned}
g_1(\mathbf{x}) &= m_3(N_6 + 2.5) - D_{\max} \le 0, \\
g_2(\mathbf{x}) &= m_1(N_1 + N_2) + m_1(N_2 + 2) - D_{\max} \le 0, \\
g_3(\mathbf{x}) &= m_3(N_4 + N_5) + m_3(N_5 + 2) - D_{\max} \le 0, \\
g_4(\mathbf{x}) &= \left| m_1(N_1 + N_2) - m_3(N_6 - N_3) \right| - m_1 - m_3 \le 0, \\
g_5(\mathbf{x}) &= -(N_1 + N_2)\sin(\pi/p) + N_2 + 2 + \delta_{22} \le 0, \\
g_6(\mathbf{x}) &= -(N_6 - N_3)\sin(\pi/p) + N_3 + 2 + \delta_{33} \le 0, \\
g_7(\mathbf{x}) &= -(N_4 + N_5)\sin(\pi/p) + N_5 + 2 + \delta_{55} \le 0, \\
g_8(\mathbf{x}) &= (N_3 + N_5 + 2 + \delta_{35})^2
- (N_6 - N_3)^2 - (N_4 + N_5)^2 \\
&\quad + 2(N_6 - N_3)(N_4 + N_5)\cos\!\left(\frac{2\pi}{p} - \beta\right) \le 0, \\
g_9(\mathbf{x}) &= N_4 - N_6 + 2N_5 + 2\delta_{56} + 4 \le 0, \\
g_{10}(\mathbf{x}) &= 2N_3 - N_6 + N_4 + 2\delta_{34} + 4 \le 0, \\
h_1(\mathbf{x}) &= \frac{N_6 - N_4}{p} \in \mathbb{Z}.
\end{aligned}
\end{align*}

\item The corresponding gear ratios are defined as
\[
i_1 = \frac{N_6}{N_4}, \quad
i_2 = \frac{N_6(N_1N_3 + N_2N_4)}{N_1N_3(N_6 - N_4)}, \quad
i_R = -\frac{N_2N_6}{N_1N_3},
\]
with target values
\[
i_{01} = 3.11, \quad i_{02} = 1.84, \quad i_{0R} = -3.11.
\]

\item Variable bounds
\begin{align*}
& N_i \in \mathbb{Z}, \quad p \in \{3,4,5\}, \\
& 17 \le N_1 \le 96, \quad 14 \le N_2 \le 54, \quad 14 \le N_3 \le 51, \\
& 17 \le N_4 \le 46, \quad 14 \le N_5 \le 51, \quad 48 \le N_6 \le 124, \\
& m_1, m_3 \in \{1.75, 2.0, 2.25, 2.5, 2.75, 3.0\}.
\end{align*}

\end{itemize}

\begin{figure}[!ht]
\centering
\begin{tikzpicture}[
    scale=1.0,
    line cap=round,
    line join=round,
    >=Stealth,
    gear/.style={draw, circle, line width=1pt},
    carrier/.style={line width=1pt, dashed},
    label/.style={font=\small}
]

\draw[gear] (0,0) circle (0.8);
\node[label] at (0,0) {$N_1$};

\draw[gear] (0,0) circle (3.2);
\node[label] at (0,3.7) {Ring gear $N_6$};

\foreach \angle/\name in {45/$N_2$,135/$N_3$,225/$N_4$,315/$N_5$} {
    \coordinate (P) at (\angle:2.0);
    \draw[gear] (P) circle (0.6);
    \node[label] at (\angle:2.9) {\name};
}

\foreach \angle in {45,135,225,315} {
    \draw[carrier] (0,0) -- (\angle:2.0);
}
\node[label] at (-3.0,-3.0) {Carrier ($p=4$)};

\draw[->,line width=0.8pt] (0.55,0.1) arc (10:350:0.55);
\node[label] at (1.55,0.65) {$i_k$};

\end{tikzpicture}

\caption{Schematic diagram of the planetary gear train design problem.
The system consists of a sun gear ($N_1$), four planet gears ($N_2$--$N_5$)
mounted on a carrier, and an internal ring gear ($N_6$).
The carrier holds $p=4$ equally spaced planets, and the gear ratios $i_k$
are determined by the corresponding tooth numbers and module selections.}
\label{fig:planetary_gear}
\end{figure}
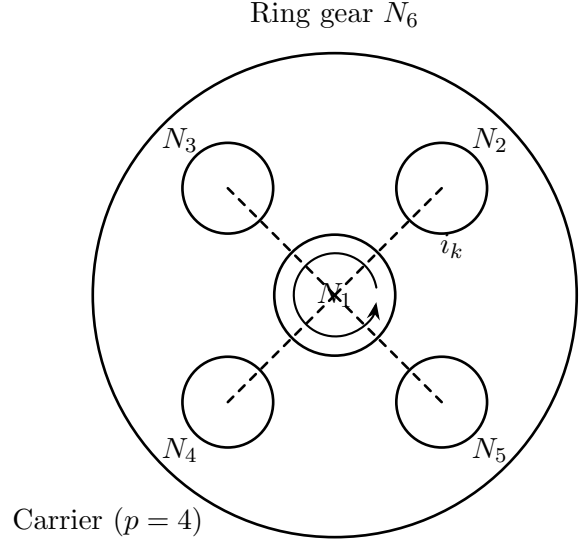

Figure~\ref{fig:planetary_gear} illustrates the schematic configuration of the planetary gear train, consisting of a central sun gear, multiple planet gears mounted on a carrier, and an internal ring gear.

Table~\ref{tab:planetary_results} reports the comparative performance of different optimization algorithms on the planetary gear train design problem.

\begin{table}[!ht]
\centering
\caption{Comparison of experimental results for the planetary gear train design problem.}
\label{tab:planetary_results}
\begin{tabular}{lcccccccccc}
\toprule
\multirow{2}{*}{Algorithm} & \multicolumn{9}{c}{Optimal values for variables} & \multirow{2}{*}{Optimal value} \\
\cmidrule(lr){2-10}
  & $p$ & $N_6$ & $N_5$ & $N_4$ & $N_3$ & $N_2$ & $N_1$ & $m_3$ & $m_1$ &   \\
\midrule
GA   & 30 & 21 & 21 & 22 & 17 & 78 & 2 & 1.75 & 4 & 0.5742857143 \\
PSO  & 20 & 29 & 15 & 17 & 15 & 54 & 3 & 3.00 & 3 & $1.13 \times 10^{104}$ \\
WSO  & 44 & 28 & 20 & 24 & 21 & 84 & 2 & 2.25 & 4 & 0.5328468900 \\
SAO  & 23 & 14 & 14 & 17 & 14 & 56 & 2 & 1.75 & 3 & 0.6752173913 \\
BKA  & 32 & 20 & 15 & 17 & 18 & 62 & 1.75 & 2.00 & 3 & 0.5370588235 \\
LSHADE  & 37 & 22 & 20 & 24 & 19 & 87 & 2.25 & 2.00 & 3 & 0.5262805663 \\
KEO  & 35 & 28 & 28 & 25 & 20 & 91 & 2.25 & 2.25 & 3 & 0.5300000000 \\
\textbf{HCKEO}
     & \textbf{41} & \textbf{28} & \textbf{24} & \textbf{25}
     & \textbf{22} & \textbf{91} & \textbf{2.00} & \textbf{2.00} & \textbf{3} & \textbf{0.5257687075} \\
\bottomrule
\end{tabular}
\end{table}

\subsubsection{Reactor network design}

The reactor network design problem is a benchmark optimization task in chemical engineering. The core objective is to design an optimal sequence consisting of two Continuously Stirred Tank Reactors (CSTRs). The primary goal of this optimization is to maximize the final concentration of the desired product through the manipulation of various process variables.
The problem’s mathematical model is formulated as follows:
\begin{itemize}

\item The decision variable vector is defined as
\[
\mathbf{x} = [x_1, x_2, x_3, x_4, x_5, x_6],
\]
where $x_1$--$x_4$ are continuous design variables bounded in $[0,1]$, and
$x_5$, $x_6$ are auxiliary positive variables associated with nonlinear constraints.

\item The optimization problem is stated as
\begin{align*}
\text{Minimize:} \quad
& f(\mathbf{x}) = x_4, \\
\text{Subject to:} \quad
& \begin{aligned}
h_1(\mathbf{x}) &= k_1 x_5 x_2 + x_1 - 1 = 0, \\
h_2(\mathbf{x}) &= k_3 x_5 x_3 + x_3 + x_1 - 1 = 0, \\
h_3(\mathbf{x}) &= k_2 x_6 x_2 - x_1 + x_2 = 0, \\
h_4(\mathbf{x}) &= k_4 x_6 x_4 + x_2 - x_1 + x_4 - x_3 = 0, \\
g_1(\mathbf{x}) &= \sqrt{x_5} + \sqrt{x_6} - 4 \le 0.
\end{aligned}
\end{align*}

\item The model parameters are fixed as
\[
k_1 = 0.09755988, \quad
k_2 = 0.99\,k_1, \quad
k_3 = 0.0391908, \quad
k_4 = 0.9\,k_3.
\]

\item The variable bounds are given by
\begin{align*}
& 0 \le x_1, x_2, x_3, x_4 \le 1, \\
& 10^{-5} \le x_5, x_6 \le 16.
\end{align*}

\end{itemize}

\begin{figure}[!ht]
\centering
\begin{tikzpicture}[
    scale=1.0,
    line cap=round,
    line join=round,
    >=Stealth,
    reactor/.style={draw, rectangle, rounded corners=3pt, minimum width=2.8cm, minimum height=2cm, line width=1pt},
    stream/.style={->, line width=0.9pt},
    label/.style={font=\small}
]

\node[reactor] (R1) at (0,0) {CSTR 1};
\node[reactor] (R2) at (5,0) {CSTR 2};

\draw[stream] (-3,0) -- (R1.west);
\node[label] at (-2.2,0.3) {Feed};

\draw[stream] (R1.east) -- (R2.west);

\draw[stream] (R2.east) -- (8,0);
\node[label] at (7.2,0.3) {Product};

\draw[dashed] (0,-1.5) -- (0,-2.4);
\node[label] at (0,-2.7) {$x_5$};

\draw[dashed] (5,-1.5) -- (5,-2.4);
\node[label] at (5,-2.7) {$x_6$};

\node[label] at (0,1.4) {$x_1,\, x_2$};
\node[label] at (5,1.4) {$x_3,\, x_4$};

\end{tikzpicture}

\caption{Schematic diagram of the reactor network design problem.
The system consists of two continuously stirred tank reactors (CSTRs) connected in series.
The residence times $x_5$ and $x_6$ are the primary decision variables, while the state variables $x_1$--$x_4$ represent the component concentrations.
The objective is to maximize the final product concentration at the outlet of the second reactor.}
\label{fig:reactor_network}
\end{figure}
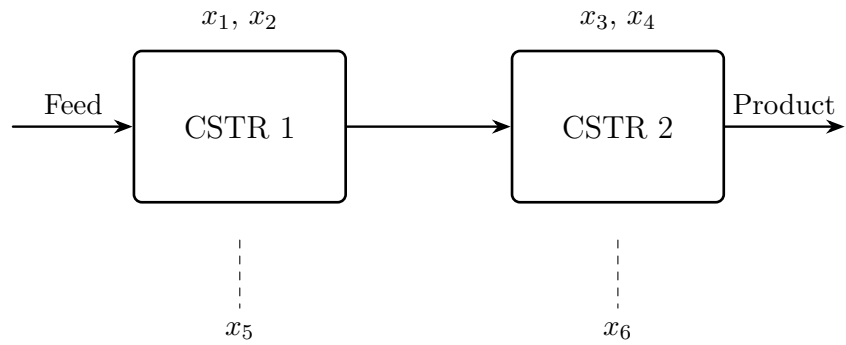

\begin{table}[!ht]
\centering
\caption{Optimization results for the reactor network design problem.}
\label{tab:reactor_results}
\begin{tabular}{lccccccc}
\toprule
\multirow{2}{*}{Algorithm} & \multicolumn{5}{c}{Optimal values for variables} & \multirow{2}{*}{Optimal value} \\
\cmidrule(lr){2-6}
  & $x_1$ & $x_2$ & $x_3$ & $x_5$ & $x_6$ &   \\
\midrule
GA      & 0.4940 & 0.4925 & 0.3576 & 10.5462 & 0.0234 & -0.3472 \\
PSO     & 0.6718 & 0.5864 & 0.2804 & 4.5071 & 4.0056 & 1127.0736 \\
SMO     & 0.7402 & 0.6094 & 0.3605 & 5.0385 & 2.8299 & 651.36910 \\
WSO     & 0.5937 & 0.5268 & 0.3102 & 7.9033 & 1.3149 & -0.3604 \\
SAO     & 0.9628 & 0.8295 & 0.0366 & 0.4600 & 1.6624 & -0.1604 \\
BKA     & 1.0000 & 1.0000 & 0.0000 & 0.0000 & 0.0000 & $1.88 \times 10^{-8}$ \\
QSSA     & 0.9993 & 0.9986 & 0.0006 & 0.0064 & 0.0072 & -0.0013 \\
LSHADE     & 0.7717 & 0.5628 & 0.1963 & 4.1576 & 3.8445 & -0.3569 \\
KEO     & 0.9600 & 0.8771 & 0.0391 & 0.4580 & 0.9766 & -0.1100 \\
\textbf{HCKEO} & \textbf{1.0000} & \textbf{0.3933} & \textbf{0.0000} & \textbf{0.0000} & \textbf{15.9742} & \textbf{-0.3881} \\
\bottomrule
\end{tabular}
\end{table}

\subsubsection{Experimental analysis and discussion}

The experimental results (Table \ref{tab:pvd_results}--\ref{tab:reactor_results}) on ten real-world engineering optimization problems demonstrate the overall robustness and effectiveness of the proposed HCKEO algorithm when compared with a variety of state-of-the-art metaheuristic methods. These benchmark problems cover a wide range of engineering characteristics, including structural weight minimization, system and network synthesis, and mechanical transmission design, thereby providing a comprehensive evaluation of algorithmic performance.

For structural and weight minimization problems, including the pressure vessel, tension/compression spring, speed reducer, and welded beam design cases, HCKEO consistently achieved either the best or near-optimal solutions. In the pressure vessel problem, HCKEO obtained the lowest cost value, slightly outperforming advanced differential evolution variants and significantly surpassing GA and PSO, indicating its capability in handling variables with discrete-like behavior. Similar stability was observed in the spring design problem, where HCKEO matched the best-known solution under highly nonlinear stress and frequency constraints. In the speed reducer and welded beam problems, HCKEO demonstrated superior precision and feasibility, whereas several baseline algorithms failed to converge or violated constraints, highlighting its reliability in tightly constrained structural optimization tasks.

In network and system synthesis problems, such as the heat exchanger network and reactor network design, the presence of multiple nonlinear equality constraints poses significant challenges to conventional metaheuristics. The results show that GA and PSO frequently diverged or produced infeasible solutions with extremely large objective values, reflecting their limited constraint-handling capability. In contrast, HCKEO consistently obtained high-quality feasible solutions and achieved the best-known objective values, particularly in the reactor network problem, where mass balance equations must be strictly satisfied. These results indicate that the hybrid search strategy embedded in HCKEO is effective in navigating narrow feasible regions defined by coupled equality constraints.

For mechanical linkage and transmission problems, including the three-bar truss, step-cone pulley, robot gripper, and planetary gear train design, HCKEO again exhibited strong performance. In the three-bar truss problem, the algorithm successfully converged to the theoretical minimum located on the boundary of stress constraints. The step-cone pulley problem proved especially challenging for traditional algorithms, many of which returned excessively large objective values due to constraint violations. HCKEO, however, achieved the best feasible solution, demonstrating its ability to handle narrow feasible manifolds and complex geometric relationships. In the robot gripper and planetary gear train problems, which involve geometric feasibility conditions and mixed-integer decision variables, HCKEO attained the lowest objective values, indicating stable performance in discrete and hybrid search spaces.

Overall, the experimental evidence suggests that HCKEO exhibits strong constraint-handling capability, high numerical precision, and stable convergence behavior across diverse engineering problem categories. These characteristics make it a reliable and competitive optimization framework for complex real-world engineering design problems involving nonlinear, constrained, and mixed-variable search spaces.

\section{UAV Path Planning: Experimental Results and Analysis}\label{UAV}

In this study, the UAV path planning problem is formulated through a comprehensive cost function that integrates optimality criteria and operational constraints, as described in this section.

\subsection{Mathematical Model Definition}

To ensure safe and efficient flight operations, the UAV path planning task can be modeled as a constrained multi-objective optimization problem.
The objective is to determine an optimal trajectory that minimizes flight distance and energy consumption while satisfying environmental and physical constraints.
In this study, the UAV flight path is represented as a sequence of discrete waypoints in three-dimensional space, each corresponding to a feasible position that the vehicle can reach.
The mathematical model integrates four primary components:
(1) path optimality, which ensures the trajectory is as short as possible;
(2) safety and feasibility constraints, which prevent collisions with obstacles and maintain an adequate safety margin;
(3) altitude constraints, which keep the UAV within allowable operational height limits; and
(4) path smoothness, which guarantees dynamic feasibility by restricting excessive turning and climbing maneuvers.
The overall cost function is then formulated as a weighted combination of these components, serving as the optimization objective for the proposed algorithm.

\subsubsection{Path Optimality}

For practical UAV missions, the generated flight trajectory should be optimized according to specific performance criteria determined by the operational objectives.
In aerial photography, mapping, and surface inspection applications, minimizing the total path length is typically the primary objective.
Let the flight path of the $i$-th UAV be denoted as
$T_i = \{ T_{i1}, T_{i2}, \ldots, T_{in} \}$,
where each waypoint $T_{ij} = (x_{ij}, y_{ij}, z_{ij})$ represents a node in the search space.
The Euclidean distance between two consecutive nodes is expressed as
$\left\| T_{ij}T_{i,j+1} \right\|$.
Thus, the path length cost $F_1$ can be defined as

\begin{equation}\label{F1}
F_1(T_i) = \sum_{j=1}^{n-1} \left\| T_{ij}T_{i,j+1} \right\|.
\end{equation}

\subsubsection{Safety and Feasibility Constraints}

\begin{figure}[!htbp]
  \centering
  \includegraphics[scale=1]{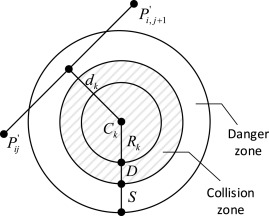}
  \caption{Determination of the threat cost.}
  \label{fig:threat_cost}
\end{figure}

In addition to optimality, the generated path must ensure safe and feasible UAV navigation by avoiding threats or obstacles within the operating space.
Let $K$ denote the set of all threats, each modeled as a cylinder whose projection is characterized by its center coordinate $C_k$ and radius $R_k$, as illustrated in Fig.~\ref{fig:threat_cost}.
For each path segment $\left\| T_{ij}T_{i,j+1} \right\|$, the corresponding threat cost depends on its distance $d_k$ to the threat center $C_k$.
By considering the UAV diameter $D$ and a safety margin $S$, the total threat cost $F_2$ is formulated as

\begin{equation}\label{F2}
\left\{
\begin{aligned}
F_2(T_i) &= \sum_{j=1}^{n-1} \sum_{k=1}^{K} O_k \left( T_{ij}T_{i,j+1} \right), \\[6pt]
O_k \left( T_{ij}T_{i,j+1} \right) &=
\begin{cases}
0, & d_k > S + D + R_k, \\[4pt]
(S + D + R_k) - d_k, & D + R_k < d_k \le S + D + R_k, \\[4pt]
\infty, & d_k \le D + R_k.
\end{cases}
\end{aligned}
\right.
\end{equation}
Here, $D$ depends on the UAV’s physical dimensions, whereas $S$ is related to mission type, environmental complexity, and GPS accuracy.
Typically, $S$ can range from tens of meters in static environments with strong GPS signals to hundreds of meters in dynamic or GPS-degraded environments.

\subsubsection{Altitude Constraint}

\begin{figure}[!ht]
  \centering
  \includegraphics[scale=1]{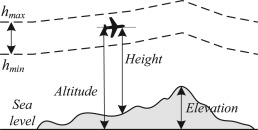}
  \caption{Illustration of the altitude cost computation.}
  \label{fig:altitude_cost}
\end{figure}

During UAV operation, the flight altitude is often restricted between two extrema, namely the minimum and maximum allowable heights.
For example, in aerial survey or search applications, the camera must operate within a specific range to guarantee the required image resolution and field of view.
Let $h_{\min}$ and $h_{\max}$ denote the minimum and maximum permissible heights, respectively.
The altitude cost associated with each waypoint $T_{ij}$ is defined as

\begin{equation}
H_{ij} =
\begin{cases}
\left| h_{ij} - \dfrac{h_{\max} + h_{\min}}{2} \right|, & \text{if } h_{\min} \le h_{ij} \le h_{\max}, \\[6pt]
\infty, & \text{otherwise,}
\end{cases}
\end{equation}
where $h_{ij}$ is the flight height relative to the ground, as illustrated in Fig.~\ref{fig:altitude_cost}.
This formulation maintains the flight altitude near the average height and penalizes out-of-range values.
The total altitude cost is then expressed as

\begin{equation}\label{F3}
F_3(T_i) = \sum_{j=1}^{n} H_{ij}.
\end{equation}

\subsubsection{Path Smoothness}

\begin{figure}[!htbp]
  \centering
  \includegraphics[scale=1]{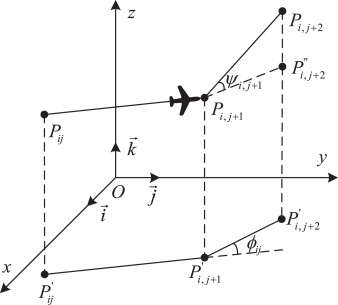}
  \caption{Calculation of turning and climbing angles.}
  \label{fig:turning_angle}
\end{figure}

Path smoothness evaluates the turning and climbing rates, which are crucial for generating dynamically feasible UAV trajectories.
As illustrated in Fig.~\ref{fig:turning_angle}, the turning angle $\phi_{ij}$ is defined as the angle between two consecutive path segments,
$T'_{ij}T'_{i,j+1}$ and $T'_{i,j+1}T'_{i,j+2}$, projected onto the horizontal plane.
Let $\vec{k}$ denote the unit vector along the $z$-axis. The projected vector can be obtained as

\begin{equation}
T'_{ij}T'_{i,j+1} = \vec{k} \times \left( P_{ij}\vec{P}_{i,j+1} \times \vec{k} \right).
\end{equation}
Hence, the turning angle is computed by

\begin{equation}
\phi_{ij} = \arctan \left(
\frac{\left\| \vec{P}_{ij}P'_{i,j+1} \times \vec{P}_{i,j+1}P'_{i,j+2} \right\|}
{\vec{P}_{ij}P'_{i,j+1} \cdot \vec{P}_{i,j+1}P'_{i,j+2}}
\right).
\end{equation}
Similarly, the climbing angle $\psi_{ij}$ represents the inclination between the 3D path segment $P_{ij}P_{i,j+1}$ and its horizontal projection $P'_{ij}P'_{i,j+1}$:

\begin{equation}
\psi_{ij} = \arctan \left(
\frac{z_{i,j+1} - z_{ij}}
{\left\| \vec{P}_{ij}P'_{i,j+1} \right\|}
\right).
\end{equation}
The smoothness cost is therefore defined as

\begin{equation}
F_4(T_i) = a_1 \sum_{j=1}^{n-2} \phi_{ij}
+ a_2 \sum_{j=1}^{n-1} \left| \psi_{ij} - \psi_{i,j-1} \right|,
\end{equation}
where $a_1$ and $a_2$ are weighting coefficients associated with the turning and climbing penalties, respectively.

\subsubsection{Overall Cost Function}

By integrating the optimality, safety, altitude, and smoothness considerations, the overall cost function for a given path $T_i$ can be expressed as

\begin{equation}\label{all_cost}
F(T_i) = \sum_{k=1}^{4} b_k F_k(T_i),
\end{equation}
where $b_k$ denotes the weight coefficient associated with each component.
This overall cost function serves as the objective to be minimized in the UAV path planning process.

\subsection{Scenario Setup}

\begin{figure}[!htbp]
  \centering
  \includegraphics[scale=0.65]{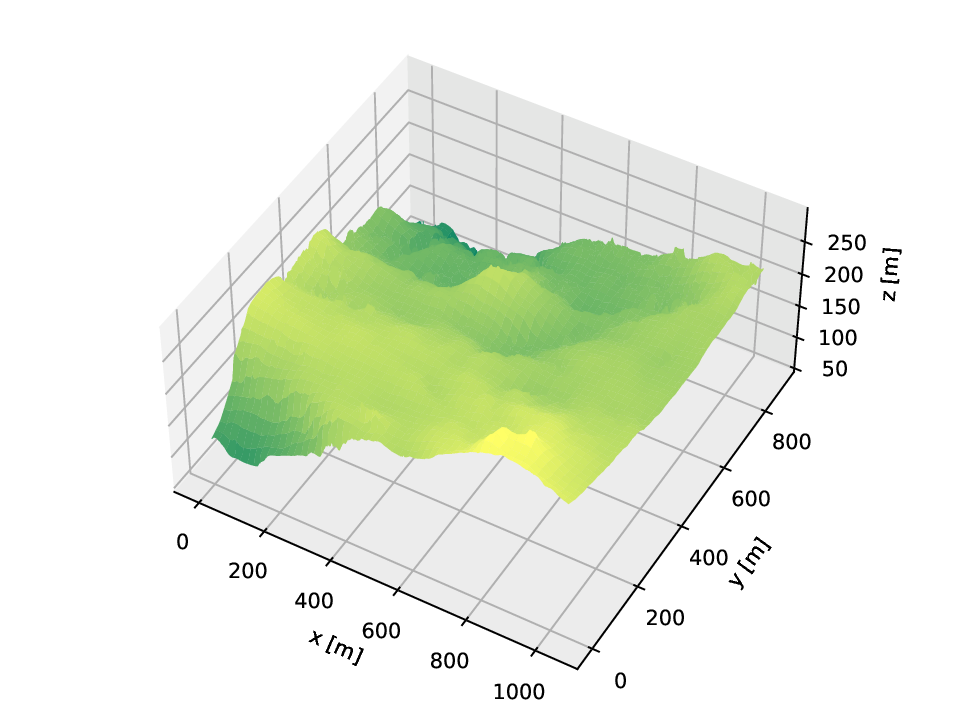}
  \caption{Calculation of turning and climbing angles.}
  \label{fig:3dmap}
\end{figure}

\begin{figure}[!htbp]
  \centering
  \includegraphics[scale=0.65]{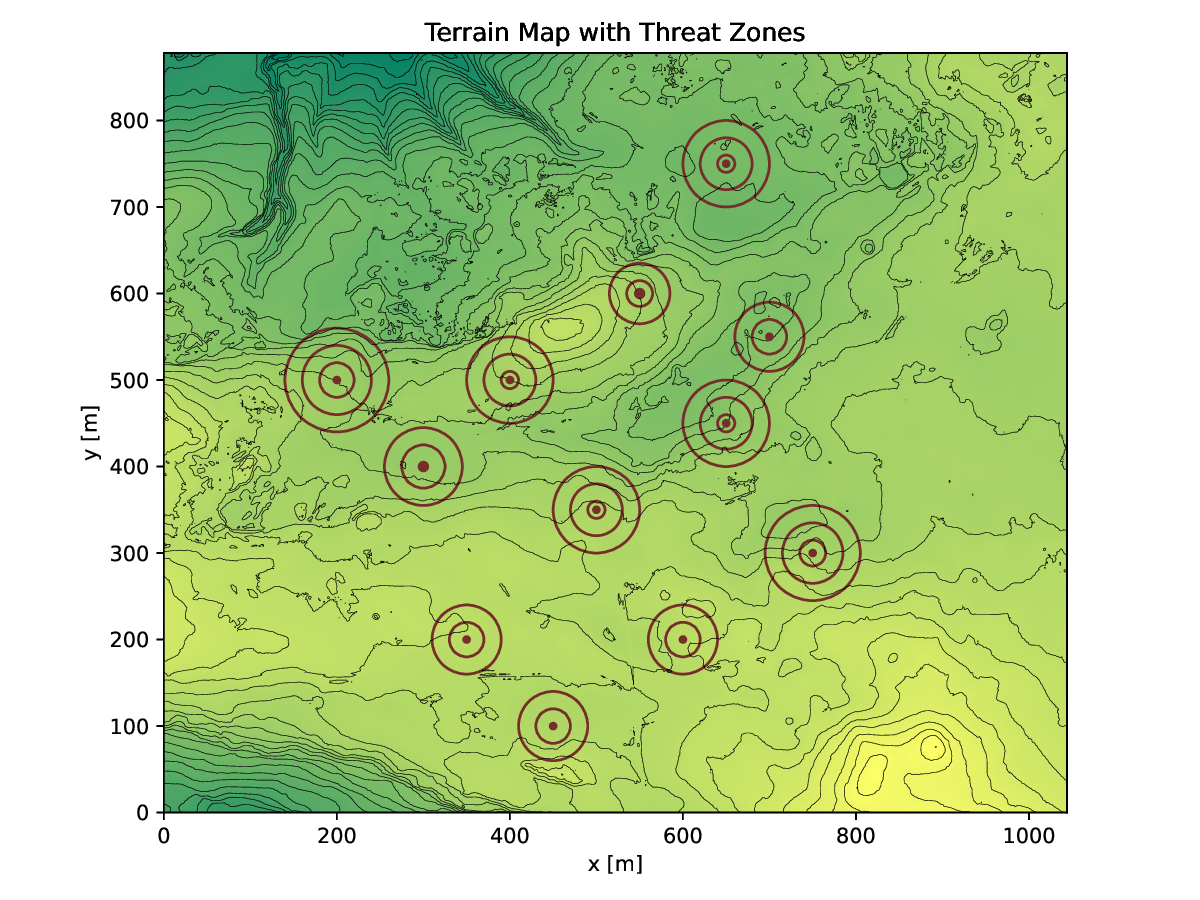}
  \caption{Calculation of turning and climbing angles.}
  \label{fig:sensors}
\end{figure}

The scenario used for evaluation is based on a real digital elevation model (DEM) map obtained from a LiDAR sensor \cite{Geoscience2015}, which is further complicated.
As shown in Fig \ref{fig:3dmap}, an area of Christmas Island in Australia was selected and 12 threat zones were added (as indicated by the red circular area in the Fig \ref{fig:sensors}).

\subsection{Result and Analysis}

To comprehensively evaluate the performance of the proposed HCKEO algorithm in UAV path planning, a comparative study was conducted against several representative metaheuristic algorithms, including GA, PSO, WSO, SAO, BKA, LSHADE, KEO and SPSO \cite{Manh2021} (a various of PSO for UAV path planning).
Each algorithm was executed under identical parameter settings and environmental conditions, with the weight coefficients $b_1$, $b_2$, $b_3$, and $b_4$ in Eq.~(\ref{all_cost}) set to 5, 1, 10, and 1, respectively.
In addition, the infinite penalty values defined in Eq.~(\ref{F2}) and Eq.~(\ref{F3}) were represented by $1\times10^{16}$ in the program implementation to handle constraint violations numerically.

\begin{table}[!ht]
  \centering
  \caption{Fitness values of the paths generated by HCKEO and other methods}\label{tab:path_result}
  \begin{tabular}{llllll}
    \toprule
    Algorithm & F1 & F2 & F3 & F4 & F  \\
    \midrule
    GA & 942.4542 & 30.4300 & 36.9316 & 0 & 51120.1691 \\
    PSO & \textbf{166.9197} & 1.0000E+17 & 160.4703 & 459.0027 & 1.0000E+17 \\
    SMO & 1330.2852 & \textbf{7.4906} & 252.0680 & 233.5092 & 9413.1057 \\
    WSO & 939.3808 & 11.3319 & 28.1909 & 0 & 4986.3736 \\
    SAO & 935.7081 & 1.0000E+17 & 0.7971 & 0 & 1.0000E+17 \\
    BKA & 1052.6512 & 20.7338 & 78.7665 & 0 & 6071.6549 \\
    QSSA & 1423.4461 & 1.0000E+17 & 227.4372 & 464.4064 & 7.0000E+17 \\
    LSHADE & 941.2683 & 17.9923 & 27.0719 & 0 & 4995.0523 \\
    KEO & 929.3484 & 1.0000E+17 & 0.5962 & 0 & 1.0000E+17 \\
    SPSO & 1038.8484 & 50.9854 & 353.7158 & 150.4522 & 8932.8383\\
    HCKEO & 931.3509 & 9.4698 & \textbf{1.2419E-05} & \textbf{0} & \textbf{4667.4820} \\
    \bottomrule
  \end{tabular}
\end{table}

Table~\ref{tab:path_result} summarizes the fitness components obtained by all tested algorithms.
It can be observed that HCKEO achieved the lowest overall cost ($F_5$) among all competitors, with a total value of 4667.4820, outperforming other metaheuristics by a considerable margin.
Although PSO obtained the smallest path length ($F_1$ = 166.9197), its safety cost ($F_2$) and overall cost ($F_5$) reached the maximum penalty ($1\times10^{17}$), indicating that the generated trajectory violated the threat avoidance constraint and thus was infeasible.
Similarly, SAO and KEO exhibited extremely low altitude and smoothness costs, but their paths intersected with threat zones, resulting in infeasible solutions.
In contrast, HCKEO successfully maintained a balance between path optimality, safety, and smoothness, yielding a feasible and collision-free trajectory with minimal constraint violations.
Compared with the basic KEO, the introduction of the L\'{e}vy long jump and horizontal–vertical crossover strategies significantly improved global exploration capability and convergence stability, as evidenced by the markedly reduced $F_2$ and $F_3$ values.

The specific path coordinates generated by HCKEO are listed in Table~\ref{tab:path}.
It can be seen that the UAV departs from $(200, 100, 150)$ and travels through ten intermediate waypoints before reaching the target point $(800, 800, 150)$.
Throughout the trajectory, the altitude remains nearly constant at approximately 150~m, satisfying the prescribed altitude constraint.
The path exhibits smooth directional transitions without sharp turns, indicating good dynamic feasibility for UAV motion control.
Furthermore, the intermediate waypoints are spatially distributed along a gradually ascending diagonal direction, effectively avoiding all predefined threat zones.
This demonstrates that HCKEO not only ensures safe and feasible flight but also minimizes unnecessary detours, achieving both efficiency and robustness.

\begin{figure}[!ht]
  \centering
  \includegraphics[scale=0.65]{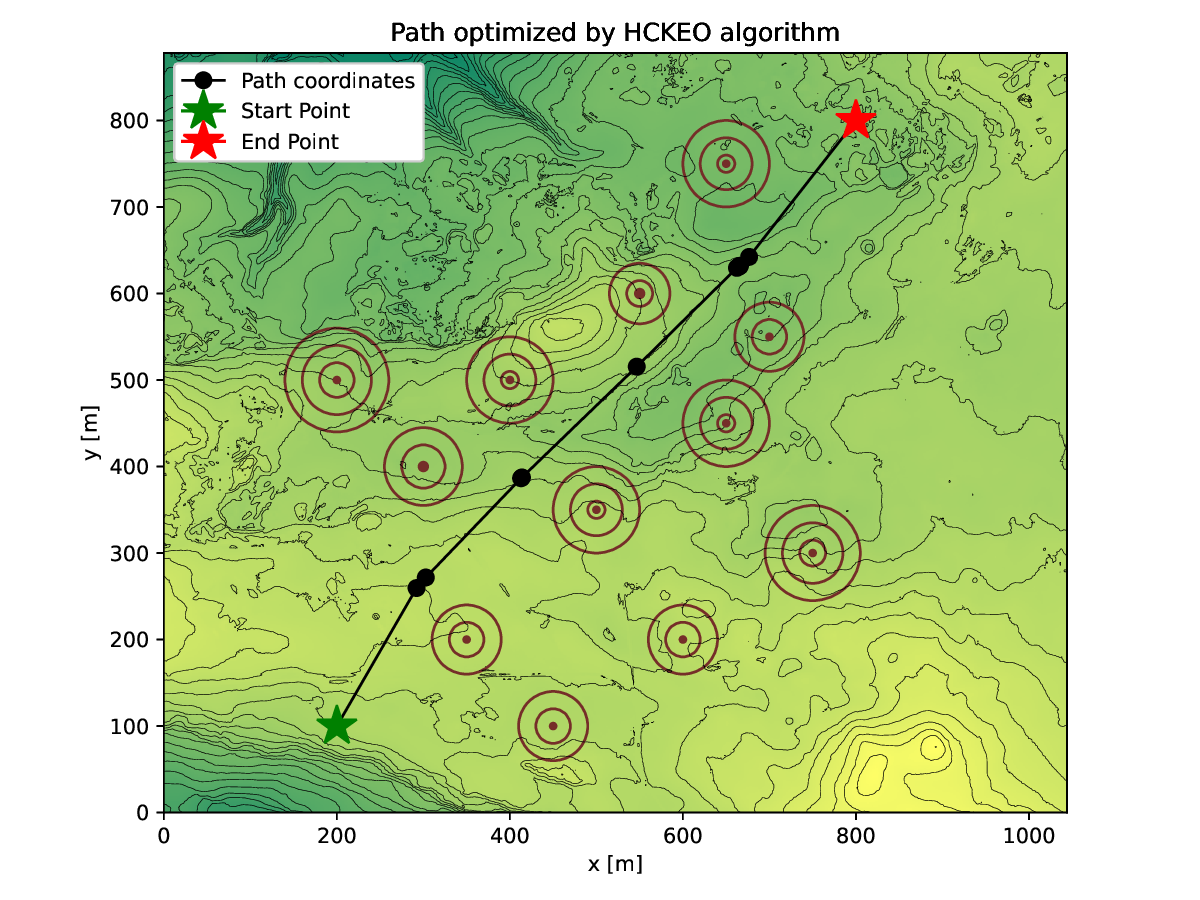}
  \caption{The path optimized by HCKEO.}
  \label{fig:opt_path}
\end{figure}

In summary, the proposed HCKEO demonstrates strong optimization capability in UAV path planning tasks.
By enhancing both global exploration and local exploitation through dual crossover strategies, it achieves superior solution quality, robust constraint handling, and smooth trajectory generation compared with other state-of-the-art metaheuristics.
The resulting feasible path effectively balances flight distance, safety, altitude stability, and smoothness, proving the algorithm’s potential for practical UAV navigation and mission planning applications.
The specific path coordinates obtained by HCKEO are shown in Table \ref{tab:path}.
And the path optimized by HCKEO is shown in Fig \ref{fig:opt_path}.
\begin{table}[!ht]
  \centering
  \caption{The path coordinates obtained by HCKEO}\label{tab:path}
  \begin{tabular}{lccc}
    \toprule
     & $x$ & $y$ & $z$   \\
    \midrule
    Starting point & 200.0000 & 100.0000 & 150.0000 \\
    1 & 292.1247 & 259.5185 & 150.0009 \\
    2 & 302.8170 & 271.7684 & 150.0007 \\
    3 & 412.9512 & 386.4542 & 149.8908 \\
    4 & 413.9689 & 387.4040 & 150.0027 \\
    5 & 546.6148 & 515.4863 & 149.9994 \\
    6 & 662.4394 & 629.6276 & 149.9921 \\
    7 & 665.1117 & 631.4551 & 149.9997 \\
    8 & 665.1143 & 631.4551 & 149.9999 \\
    9 & 665.9281 & 632.0249 & 150.0034 \\
    10 & 676.5317 & 642.3140 & 150.0000 \\
    End point & 800.0000 & 800.0000 & 150.0000 \\
    \bottomrule
  \end{tabular}
\end{table}

\section{Multi-backend Implementation and Computational Performance}\label{fealpy}
Metaheuristic algorithms have traditionally been implemented in high-level numerical environments such as MATLAB or NumPy-based Python, which provide convenient abstractions for vectorized computation and rapid algorithm prototyping \cite{Ilhem2013}.
With the growing demand for large-scale and high-dimensional optimization, heterogeneous computing platforms, especially GPU-accelerated architectures, are increasingly adopted to improve computational efficiency \cite{Crainic2024}.

In such environments, numerical acceleration is typically achieved through tensor-based computation.
Notably, the population of metaheuristic algorithms is naturally represented in matrix form, while tensors can be regarded as higher-order generalizations that naturally extend such representations \cite{Jean2019}.
As a result, tensor-based abstractions provide a natural and unified representation for population-level operations, particularly when extending algorithm execution across different computing backends and hardware platforms.

To address these requirements, the proposed HCKEO is implemented based on the FEALPy framework, which provides a unified backend abstraction layer for numerical computation.
This design enables the same optimization logic to be executed seamlessly across different numerical libraries and computing platforms, including NumPy and PyTorch, as well as CPU- and GPU-based architectures, without modifying the core algorithmic structure.
To the best of our knowledge, this is the first implementation of the KEO family within a unified tensor-based multi-backend framework.

Rather than emphasizing algorithmic competitiveness in this section, a simple benchmark function $Step$, defined as $f(x)=\sum_{i=1}^{D}(x_i+0.5)^2$, is deliberately adopted to isolate the effects of backend implementation.
This experimental design allows the numerical consistency and computational efficiency of HCKEO to be evaluated across different platforms, independent of problem-specific difficulty or search complexity.

Experiments in this case are conducted using identical algorithmic parameters and stopping criteria to ensure a fair and consistent comparison across different implementations.
The results, summarized in Table~\ref{tab:compare}, show that the proposed implementation achieves consistent solution quality across all tested platforms.
Moreover, significant computational speedups are observed on GPU-enabled backends, especially for large-scale and high-dimensional problem settings, highlighting the scalability and practical efficiency of the multi-backend HCKEO implementation.
It should be emphasized that these speedups primarily stem from backend-level parallelism rather than algorithmic modifications.

These results indicate that the proposed HCKEO is not only algorithmically effective but also implemented in an engineering-oriented manner, making it suitable for practical deployment in heterogeneous computing environments.

\begin{table}[H]
\centering
\caption{Performance comparison of the proposed algorithm on a test function
under different problem dimensions, population sizes, programming languages,
and computing platforms}
\label{tab:compare}
\begin{tabular}{cccc}
\toprule
Problem setting $(D, N)$ & Backend / Platform & Best fitness value & Runtime (s) \\
\midrule
\multirow{4}{*}{$D=100$, $N=10000$}
& Matlab           & $4.4772\times10^{-24}$ & 1221.2508 \\
& NumPy (CPU)      & $0$ & 53.4238 \\
& PyTorch (CPU)    & $0$ & 18.0221 \\
& PyTorch (GPU)    & $0$ & \textbf{1.3212} \\
\midrule
\multirow{4}{*}{$D=50$, $N=1000$}
& Matlab           & $3.7286\times10^{-31}$ & 29.9382 \\
& NumPy (CPU)      & $0$ & 1.2883 \\
& PyTorch (CPU)    & $0$  & 1.1620 \\
& PyTorch (GPU)    & $0$ & \textbf{1.1078} \\
\midrule
\multirow{4}{*}{$D=30$, $N=100$}
& Matlab           & $0$ & 2.9667 \\
& NumPy (CPU)      & $0$ & \textbf{0.1340} \\
& PyTorch (CPU)    & $0$ & 0.4913 \\
& PyTorch (GPU)    & $0$ & 1.1226 \\
\bottomrule
\end{tabular}
\end{table}

Listing \ref{lst:example} shows the script of this example.
It can be seen that seamless switching between computing backend and computing platform can be achieved with just one line of code.

\begin{lstlisting}[caption={Example of backend-agnostic HCKEO usage based on FEALPy}, label={lst:example}]
# import FEALPy backend manager
from fealpy.backend import backend_manager as bm
# import HCKEO algorithm
from fealpy.opt import initialize, opt_alg_options, HybridCrossoverKangarooEscapeOpt

# switch backend and computing device without changing algorithm code
bm.set_backend('pytorch')
bm.set_default_device('cuda')

# define objective function
def fobj(x):
    return bm.sum((x + 0.5) ** 2, axis=-1)

# parameter setting
pop_size = 10000
dim = 100
ub = 100 * bm.ones((dim,))
lb = -100 * bm.ones((dim,))

# initialization
x0 = initialize(pop_size, dim, ub, lb)
options = opt_alg_options(x0, fobj, (lb, ub), pop_size)
optimizer = HybridCrossoverKangarooEscapeOpt(options)

# run optimization
optimizer.run()
optimizer.print_optimal_result()
\end{lstlisting}

\section{Conclusion and Future Work}\label{conclusion}

This paper investigated the limitations of the original Kangaroo Escape Optimizer (KEO), particularly its insufficient local search capability and susceptibility to premature convergence when addressing complex and high-dimensional optimization problems.
To overcome these shortcomings, a novel Hybrid Crossover Kangaroo Escape Optimizer (HCKEO) was proposed.
By incorporating a LLC strategy and a HVC strategy, the proposed algorithm effectively enhances population diversity, strengthens global exploration, and improves the ability to escape from local optima.

The performance of HCKEO was comprehensively evaluated on the IEEE CEC2022 benchmark suite as well as a set of real-world constrained engineering design problems.
The experimental results demonstrate that HCKEO consistently outperforms nine state-of-the-art metaheuristic algorithms in terms of convergence speed, solution accuracy, and robustness.
Furthermore, HCKEO was successfully applied to UAV path planning in a complex digital elevation model (DEM) environment, where it generated smooth, safe, and collision-free trajectories while achieving a favorable balance between flight cost and operational safety.

Beyond the algorithmic contributions, HCKEO was developed within a multi-backend computational framework based on FEALPy, enabling unified execution across different numerical backends and heterogeneous computing platforms, including CPU and GPU environments.
This design facilitates fair performance comparison across platforms and enhances the practicality of the proposed algorithm for large-scale and engineering-oriented optimization tasks, without altering the core algorithmic structure.

Future work will focus on extending HCKEO in two directions.
First, adaptive mechanisms such as entropy-based parameter control and ensemble hybridization strategies will be investigated to further improve scalability and adaptability when dealing with highly nonlinear, dynamic, and constrained optimization problems.
Second, the multi-backend framework will be further exploited to explore efficient heterogeneous computing strategies, including task-level CPU--GPU cooperation and large-scale parallel optimization, thereby improving computational efficiency for real-world industrial applications.

\section*{CRediT authorship contribution statement}

H.B.L and H.W designed the algorithm, H.B.L conducted the experiments, H.W guided the analysis of the experiment results. H.B.L wrote the main manuscript text, H.W modified it. All authors reviewed the manuscript.

\section*{Declaration of competing interest}

The authors declare that they have no known competing financial interests or personal relationships that could have appeared to influence the work reported in this paper.

\section*{Ethical and informed consent for data used}
Our research is based on open source data, and there are no human subjects in the article, so ethical and informed consent are not applicable.

\section*{Acknowledgment}

This work was supported by the National Natural Science Foundation of China (Grant Nos. 12371410, 12261131501), the Postgraduate Scientific Research Innovation Project of Hunan Province (Grant No. CX20250944), the Postgraduate Scientific Research Innovation Project of Xiangtan University (Grant No. XDCX2025Y205), and the Construction of Innovative Provinces in Hunan Province (Grant No. 2021GK1010).

\section*{Data availability}

The data that support the findings of this study are available upon request.

\end{document}